\documentclass[11pt]{article}
\usepackage{amsmath}[1996/11/01]
\usepackage{amssymb,amsthm,amsxtra}
\usepackage{amscd}
\usepackage[margin=1in]{geometry}
\usepackage[title,titletoc]{appendix}

\def\[#1\]{\begin{equation}#1\end{equation}}
\makeatletter
\def\beq{%
   \relax\ifmmode
      \@badmath
   \else
      \ifvmode
         \nointerlineskip
         \makebox[.6\linewidth]%
      \fi
      $$%%$$ BRACE MATCH HACK
   \fi
}
\def\eeq{%
   \relax\ifmmode
      \ifinner
         \@badmath
      \else
         $$%%$$ BRACE MATCH HACK
      \fi
   \else
      \@badmath
   \fi
   \ignorespaces
}

\def\enddisplaymath{\eeq\global\@ignoretrue}
\makeatother

\newtheorem{thm}{Theorem}
\newtheorem{cor}[thm]{Corollary}
\newtheorem{lem}[thm]{Lemma}
\newtheorem{prop}[thm]{Proposition}
\newtheorem{conj}{Conjecture}
\theoremstyle{remark}
\newtheorem*{rem}{Remark}

\newtheorem{eg}{Example}
\theoremstyle{definition}
\newtheorem{defn}{Definition}

\numberwithin{equation}{section}
\numberwithin{thm}{section}
\numberwithin{eg}{section}
\numberwithin{defn}{section}

\newcommand{\F}{\mathbb F}
\renewcommand{\P}{\mathbb P}
\newcommand{\Q}{\mathbb Q}
\newcommand{\Z}{\mathbb Z}

\newcommand{\A}{\mathbb A}
\newcommand{\G}{\mathbb G}

\DeclareMathOperator{\HH}{HH}
\DeclareMathOperator{\Aut}{Aut}
\DeclareMathOperator{\GL}{GL}

\DeclareMathOperator{\SL}{SL}
\DeclareMathOperator{\Pic}{Pic}

\DeclareMathOperator{\NS}{NS}

\DeclareMathOperator{\Ext}{Ext}
\DeclareMathOperator{\Sym}{Sym}
\DeclareMathOperator{\End}{End}
\DeclareMathOperator{\Hom}{Hom}
\DeclareMathOperator{\Der}{Der}
\DeclareMathOperator{\Inn}{Inn}
\DeclareMathOperator{\Mat}{Mat}
\DeclareMathOperator{\Lie}{Lie}

\DeclareMathOperator{\rank}{rank}

\DeclareMathOperator{\ad}{ad}

\DeclareMathOperator{\Tr}{Tr}
\DeclareMathOperator{\Spec}{Spec}
\DeclareMathOperator{\Proj}{Proj}

\DeclareMathOperator{\im}{im}

\DeclareMathOperator{\coh}{coh}

\DeclareMathOperator{\gr}{gr}
\DeclareMathOperator{\Gr}{Gr}
\DeclareMathOperator{\OGr}{OGr}
\newcommand{\sO}{\mathcal O}

\DeclareMathOperator{\sHom}{{\mathcal H}\!{\it om}}
\DeclareMathOperator{\sEnd}{{\mathcal E}\!{\it nd}}

\DeclareMathOperator{\ord}{ord}

\begin{document}

\title{Filtered deformations of elliptic algebras}

\author{Eric M. Rains\\Department of Mathematics, California Institute of Technology}

\date{February 16, 2022}
\maketitle

\begin{abstract}
One of the difficulties in doing noncommutative projective geometry via
explicitly presented graded algebras is that it is usually quite difficult
to show flatness, as the Hilbert series is uncomputable in general.  If
the algebra has a regular central element, one can reduce to understanding
the (hopefully more tractable) quotient.  If the quotient is particularly
nice, one can proceed in reverse and find {\em all} algebras of which it is
the quotient by a regular central element (the {\em filtered deformations}
of the quotient).  We consider in detail the case that the quotient is an
{\em elliptic algebra} (the homogeneous endomorphism ring of a vector
bundle on an elliptic curve, possibly twisted by translation).  We
explicitly compute the family of filtered deformations in many cases and
give a (conjecturally exhaustive) construction of such deformations from
noncommutative del Pezzo surfaces.  In the process, we also give a number
of results on the classification of exceptional collections on del Pezzo
surfaces, which are new even in the commutative case.

\end{abstract}

\tableofcontents

\section{Introduction}

The origin of the present work was an attempt to produce elliptic analogues
of the algebras of \cite{EOR}, with a particular focus on understanding the
associated noncommutative analogues of del Pezzo surfaces.  (Note that for
our purposes, a ``(degenerate) del Pezzo surface'' is a smooth projective
surface with an {\em almost} ample anticanonical bundle; that is, the
anticanonical morphism is allowed to contract a finite configuration of
$-2$-curves.)  The algebras of \cite{EOR} were given by very simple
explicit presentations (they are generated by three elements that multiply
to 1 and have specified minimal polynomials), and are difficult to
generalize as such, but they also have (several) filtrations such that the
associated gradeds have nice {\em geometric} descriptions.  In particular,
the associated graded algebras are constructed from the
multiplication-by-$q$ map on a nodal curve of genus 1, and thus {\em their}
elliptic analogues are relatively easy to describe.  One is thus left with
the problem of understanding their ``filtered deformations'', i.e., for
which filtered algebras does one obtain the given associated graded.
Although the problem of classifying filtered deformations of a given graded
algebra is almost certainly not computable in general, it is still
relatively straightforward to do non-rigorous calculations and thus
formulate conjectures.  In particular, we (P. Etingof and the author) found
that the (putative) moduli spaces classifying such deformations were often
quite nice: not only were they frequently rational, they were in fact quite
frequently weighted projective spaces or products thereof.

Although the resulting conjectures were quite striking, it was not clear
how to prove them rigorously; similarly, although the resulting algebras
were clearly noncommutative analogues of del Pezzo surfaces, it was unclear
how to prove anything much about them.  As a result, the project was
shelved for quite some time.  In the interim, the author developed a
completely different approach to noncommutative rational surfaces
\cite{noncomm1,noncomm2,elldaha}, extending beyond the del Pezzo case and
giving a great deal of control over the corresponding representation theory
(i.e., the categories and derived categories of sheaves), making much of
the original motivation moot.  This is not entirely the case, however: the
work of \cite{ChekovL/MazzoccoM/RubstovV:2021} in particular suggests that
there is independent interest in having representations of such surfaces
via graded algebras.  The object of the present work is therefore not so
much to construct noncommutative del Pezzo surfaces as filtered
deformations, but rather to construct filtered deformations via
noncommutative del Pezzo surfaces.

As the title indicates, the algebras we are deforming are ``elliptic
algebras'', which are constructed as follows.  The simplest version is
simply the graded homogeneous coordinate ring of an elliptic curve: if $E$
is an elliptic curve and ${\cal L}$ is a line bundle, then there is a
corresponding graded algebra
\[
B_{E,{\cal L}} = \bigoplus_{n\in \Z} \Gamma(E;{\cal L}^n),
\]
with multiplication induced by the obvious map ${\cal L}^m\otimes {\cal
  L}^n\cong {\cal L}^{m+n}$.  This has a noncommutative analogue (the
twisted homogeneous coordinate rings of \cite{ArtinM/VandenBerghM:1990}),
in which we twist by a translation $\tau\in \Aut(E)$:
\[
B_{E,{\cal L},\tau} = \bigoplus_{n\in \Z} \Gamma(E;{\cal L}_n),
\]
where
\[
{\cal L}_1 = {\cal L};\quad
{\cal L}_n = {\cal L}\otimes (\tau^{-1})^* {\cal L}_{n-1},
\]
with multiplication induced by
\[
{\cal L}_m\otimes (\tau^{-m})^* {\cal L}_n
\cong
{\cal L}_{m+n}.
\]
A somewhat cleaner description is based on the observation that the functor
\[
\Psi:M\mapsto {\cal L}\otimes (\tau^{-1})^* M
\]
is actually an autoequivalence of $\coh(E)$.  (In fact, any autoequivalence
is of this form, except that $\tau$ need not be a translation.)  We may
thus instead define
\[
B_{\Psi} = \bigoplus_{n\in \Z} \Gamma(E;\Psi^n \sO_E),
\]
with multiplication given by
\[
x_m\cdot x_n = \Psi^m(x_n)\circ x_m
\]
where we view $x_n$ as a homomorphism $\sO_E\to \Psi^n\sO_E$.

Although the filtered deformations of these algebras are enough to cover
the examples of \cite{ChekovL/MazzoccoM/RubstovV:2021}, analogues of the
algebras of \cite{EOR} require a further generalization.\footnote{We should
note here that although the algebras of \cite{EOR} are analogues of double
affine Hecke algebras, the elliptic double affine Hecke algebras of
\cite{elldaha} are {\em not} in general filtered deformations of elliptic
algebras, even in rank 1.}  For the commutative version, there is a
generalization involving a vector bundle $V$, now with
\[
B_{E,{\cal L},V} = \bigoplus_{n\in \Z} \Hom(V,V\otimes {\cal L}^n).
\]
This also has a twisted version,
\[
B_{V,\Psi} := \bigoplus_{n \in \Z} \Hom(V,\Psi^n V).
\]
Note that unlike the instance $V=\sO_E$, this may fail to be nonnegatively
graded, but will always be bounded below.

The basic question we study in the present work is: What are the filtered
deformations of $B_{V,\Psi}$?  Although a full answer appears out of reach
for the moment, we do give some partial results that narrow things down.
In one direction, we show that when $\tau=1$, any filtered deformation
arises as the endomorphism ring of a vector bundle $V^+$ on a deformation
of the appropriate cone over $E$, with the property that $V^+$ is both
rigid ($\Ext^1(V^+,V^+)=0$) and unobstructed ($\Ext^2(V^+,V^+)=0$).  Thus
the problem essentially reduces to understanding such bundles on del Pezzo
surfaces.  (In principle, one could have a nontrivial filtered deformation
of $B_{V,{\cal L}}$ living on an undeformed cone, but experiment suggests
this is not actually possible.)  More generally, we classify such bundles
on {\em noncommutative} del Pezzo surfaces, and recover the phenomenology
alluded to earlier by showing that in many cases the resulting family of
filtered deformations is a weighted projective space (or a product thereof).

The main gaps in obtaining a full classification are that (a) we only
classify rigid, unobstructed bundles on {\em smooth} del Pezzo surfaces,
while in general for $\tau=1$ the center is the coordinate ring of a {\em
  singular} del Pezzo surface, (b) when $\dim(V)\ne 1$, the algebra is not
a deformation of a del Pezzo surface per se, but rather of an Azumaya
algebra on the del Pezzo surface (in particular, there is no a priori
reason to expect it to have a well-defined structure sheaf!), and (c)
although the results of \cite{noncomm1,noncomm2,elldaha} certainly
construct many examples of noncommutative del Pezzo surfaces (obtaining the
same family in multiple ways), and some of the results below strongly
suggest that these examples are exhaustive, actually proving this is beyond
the reach of present techniques.

Note that our focus in the present work has been on the general theory,
with little attempt to give explicit presentations for any of the filtered
deformations constructed below (especially for $\tau\ne 1$).  Those
interested in explicit presentations should instead view this work as a
guide for future investigations, by showing where one expects the full
family to be particularly nice, and letting one predict the degrees of the
generators and relations.

\medskip

The plan of the paper is as follows.  In Section 2, we give a general
discussion of filtered deformations, in particular not requiring the graded
algebra to be nonnegatively graded.  The main results for elliptic algebras
are (a) the moduli stack of filtered deformations is a closed substack of a
weighted projective space, and (b) if $\tau$ is torsion, then the center of
any filtered deformation of $B_{V,\Psi}$ is a filtered deformation of the
center of $B_{V,\Psi}$.  (This section also gives explicit computations of
the moduli stack in some cases with $\tau=1$, and briefly discusses some
interesting phenomena arising when viewing $B_{V,\Psi}$ itself as a
filtered deformation.)

Any actual {\em algorithm} for classifying filtered deformations requires
understanding not just the generators and relations but also the {\em
  syzygies} of the relations, so in Section 3, we give an algorithm for
computing resolutions of the algebras $B_{V,\Psi}$.  One particularly nice
feature is that qualitative information about the resolution (i.e., the
description of each term as a sum of projective modules) can be computed
directly from qualitative information about the algebra (the degree of
${\cal L}$ and the slopes and multiplicities of the constant-slope summands
of $V$).  In particular, one finds that for nearly every case in which $V$
is semistable, the algebra $B_{V,\Psi}$ is Koszul; in particular, it is a
quadratic algebra with syzygies only in degree 3, and thus it is
straightforward to write down equations cutting out the family of filtered
deformations.  This understanding of resolutions also enables us to
control the cohomological properties of the deformations, in particular
showing that the tautological sheaf associated to any filtered deformation
is rigid and unobstructed.

This suggests in general that we should study rigid, unobstructed
torsion-free sheaves on noncommutative del Pezzo surfaces, which is done in
Section 5 after first reviewing the definition of and main results on such
surfaces in Section 4.  It is fairly straightforward to show that any such
sheaf gives rise to a filtered deformation of an elliptic algebra, and thus
in the absence of an understanding of {\em all} filtered deformations, we
would at least like to understand all the deformations arising from del
Pezzo surfaces, or equivalently to understand all such sheaves.  In the
commutative case, these were studied in
\cite{KuleshovSA/OrlovDO:1994,KuleshovSA:1997}, and it was in particular
shown that any such bundle arises in a natural way from an exceptional
collection on the surface.  It is straightforward to extend their arguments
(especially in the simplified form we give below) to the noncommutative
case, and thus the problem of understanding rigid, unobstructed,
torsion-free sheaves reduces to that of understanding exceptional
collections.  In fact, the problem simplifies further: for the cases of
interest, any sequence of exceptional bundles with the desired numerical
data is automatically an exceptional collection, so one needs only classify
exceptional bundles.  These are, moreover, uniquely determined by their
numerical data, so the only real question is existence.  This turns out to
be an open and closed condition, so that one can immediately reduce to the
commutative case (and further to nondegenerate del Pezzo surfaces of degree
1, as well as to a related question on rational elliptic surfaces).
Although this existence question remains open in most cases, we are at
least able to settle it for all cases of rank at most 9, which is enough to
cover the most interesting examples for filtered deformations.

Finally, in Section 6, we describe how to use the above theory to compute
the corresponding moduli spaces as well as their images in the moduli
spaces of filtered deformations.  Since exceptional collections are unique
when they exist, their moduli spaces are finite covers of the moduli spaces
of surfaces (the covering arising from the fact that the numerical data
breaks symmetries), and the maps from such spaces to moduli spaces of
elliptic algebras are straightforward to compute.  One finds in general
that any family of algebras arising in this way is explicitly parametrized
by the quotient of a subgroup of $E^n$ by some Weyl group acting linearly
on $E^n$.  This partly explains why these spaces are often weighted
projective spaces, as results of \cite{LooijengaE:1976} give a number of
cases in which such quotients are weighted projective spaces; the remainder
of the explanation involves showing that the natural line bundle on the
moduli stack of filtered deformations is the same as the invariant line
bundle arising in \cite{LooijengaE:1976}.  Using this, we describe the
results of the classification in a large number of cases, as well as
explaining some of the consequences arising when the result is a weighted
projective space with generators of degree $>1$.  (For instance, Sklyanin's
original deformation \cite{SklyaninEK:1982} of $\P^3$ arises in our theory
from the fact that when $\deg({\cal L})=4$, the space of deformations of
$B_{\sO_E,\Psi}$ has two generators of degree 2, which become the two
central degree 2 elements of Sklyanin's algebra.)

{\bf Acknowledgements}.  A great many thanks are due to P. Etingof, not
only for suggesting several versions of the above problem (both the
original version in which $V$ was a sum of line bundles of degree 0 and
more recent versions in which $V$ is semistable) but also for a large
number of discussions of filtered deformations over the years, including
some recent pointers regarding Poisson structures on centers.  The author
would also like to thank the organizers of the workshop ``Hypergeometry,
Integrability and Lie Theory'' at the Lorentz Center in Leiden for
providing the necessary impetus to finish and write up the present work.
In addition, the author would like to thank the referee for extensive
comments which led to many improvements in the current version of this
paper.

\section{Filtered deformations}

We first set a convention: for our purposes, a graded algebra $B$ is a
$\Z$-graded algebra (over a commutative ring $k$) in the usual sense such
that (a) $B_m$ is a finitely generated $k$-module for all $m$, and (b)
$B_m=0$ for $m\ll 0$.  (The algebras will typically be algebras over a
field $k$ of characteristic 0, but will on occasion need to be algebras
over a field of finite characteristic or a discrete valuation ring.  We
will mainly give arguments in characteristic 0, indicating in remarks where
changes need to be made to include finite characteristic.)

Similarly, our convention for filtrations is that they are {\em ascending}
filtrations, i.e., an assignment of a finitely generated submodule $A_{\le
  m}\subset A$ to each integer $m$ such that (a) $A_{\le m}A_{\le n}\subset
A_{\le m+n}$, (b) $A_{\le m}\subset A_{\le n}$ for $m\le n$, (c) $A_{\le
  m}=0$ for $m\ll 0$, and (d) $\bigcup_n A_{\le n}=A$.

Given a filtered algebra, there are two natural ways to produce a graded
algebra: the {\em associated graded} $\gr A := \bigoplus_i (A_{\le i}/A_{\le
  i-1})$, and the {\em Rees algebra} $A^+ := \bigoplus_i A_{\le i} t^i$,
where $t$ is an auxiliary variable.  Similarly, a graded algebra may be
viewed as a filtered algebra with $B_{\le i}=\bigoplus_{j\le i} B_j$; in
that case, the associated graded is $B$ and the Rees algebra is $B[t]$.

\begin{defn}
   Let $B$ be a graded algebra.  A {\em filtered deformation} of $B$ is a
   filtered algebra $A$ equipped with an isomorphism $\gr A\cong B$.
\end{defn}

\begin{rem} Note that the isomorphism is part of the data.\end{rem}

The following fact will later allow us to reduce to the case of
nonnegatively graded algebras (for the algebras of interest, at least!).

\begin{prop}\label{prop:ids_lift}
  Let $B$ be a graded algebra.  For any homogeneous idempotent $e\in B$ and
  any filtered deformation $A$ of $B$, there is an idempotent $\hat{e}\in
  A$ with leading term $e$.
\end{prop}

\begin{proof}
  An idempotent has degree 0, and thus $\hat{e}$ should be in $A_{\le 0}$.
  This is a {\em subalgebra}, and our conventions ensure that it is
  finite-dimensional, and thus Artinian.  The ideal $A_{\le -1}$ is nil, so
  any idempotent of $A_{\le 0}/A_{\le -1}=B_0$ lifts to $A_{\le 0}$, and
  thus to $A$.
\end{proof}

There is a natural action of $\G_m$ on any graded algebra $B$ with $t\in \G_m$
acting by $t^i$ on $B_i$, and this acts on the space of filtered
deformations by composing the isomorphism $\gr A\cong B$ with the induced
automorphism of $B$.

Any filtered isomorphism $A\cong A'$ of filtered algebras induces an
isomorphism $\gr A\cong \gr A'$, and thus if $A$ and $A'$ are both filtered
deformations of $B$, the isomorphism induces an automorphism of $B$.  The
natural notion of equivalence on filtered deformations is an isomorphism
such that the induced automorphism is the identity; we also consider the
notion of a {\em weak} equivalence, in which the induced automorphism comes
from the action of $\G_m$.

One example of a filtered deformation is the filtered algebra associated to
$B$ itself; in general, we say that a filtered deformation is {\em trivial}
if it is equivalent to this deformation.

\begin{defn}
  The {\em moduli stack of filtered deformations} of $B$ is the stack
  classifying nontrivial filtered deformations of $B$ up to weak
  equivalence.
\end{defn}

\begin{rem}
  There is an even weaker notion of equivalence one may consider, in which
  one allows {\em all} automorphisms of $B$ rather than just $\G_m$.  This
  is not as well-behaved in families, however, and leads to more
  complicated moduli stacks (esp. in our examples below where the present
  definition gives a weighted projective space).
\end{rem}

\begin{prop}
  If $B$ is a finitely presented graded algebra over a field $k$, then the
  moduli stack of filtered deformations of $B$ is algebraic.
\end{prop}

\begin{proof}
  Fix a finite (homogeneous) presentation of $B$, and choose a basis of $B$
  consisting of monomials in the generators.  Given a filtered deformation
  $A$ of $B$, since $\gr A\cong B$, each generator $x_i$ of $B$ has a
  preimage inside $A_{\le \deg(x_i)}$, which is unique modulo $A_{\le
    \deg(x_i)-1}$.  Choose a lift $\hat{x}_i$ of each generator. This
  induces a basis of $A$ compatible with the filtration by replacing $x_i$
  with $\hat{x}_i$ in the chosen monomial basis of $B$.

  If we make a similar substitution in a relation $\rho_i$ of $B$, then the
  result will be an element of $A_{\le \deg(\rho_i)-1}$; the degree
  $\deg(\rho_i)$ term cancels since the relation holds in $\gr A$.  This
  element is then uniquely expressible in terms of our basis of $A$, and
  thus gives us a unique lift of each relation.

  In the other direction, if we start with a collection of deformed
  relations (i.e., adding lower-order terms from the chosen basis to each
  relation from the presentation of $B$), the condition that the
  corresponding finitely presented algebra is a filtered deformation of $B$
  is closed.  To see this, note that it is equivalent to ask that the
  chosen monomials in the generators of the filtered algebra be linearly
  independent.  It follows that a given collection of deformed relations
  gives a filtered deformation iff for any combination of deformed
  relations in which only the ``good'' monomials appear, the coefficients
  all vanish.  But this is clearly a collection of polynomial equations,
  and thus the filtered deformation locus is an intersection of closed
  subschemes, so is closed.

  This shows that the moduli stack of filtered deformations {\em with
    chosen lifts of the generators} is a scheme (an affine scheme, in
  fact).  Two such deformations are equivalent iff one can be obtained from
  the other by changing the lifts of the generators.  Each generator
  $\hat{x_i}$ can be changed by adding an element of $A_{\le \deg(x_i)-1}$,
  which is in turn given by a polynomial in the generators of lower degree.
  We thus see that if we order the generators $x_1,\dots,x_n$ by degree,
  then after lifting $x_1,\dots,x_{n-1}$, the set of possible lifts of
  $x_n$ is a torsor over a unipotent group scheme, so that forgetting the
  lift takes the quotient by that group, leaving the stack algebraic.  We
  can then continue inductively, and thus find that the moduli stack of
  filtered deformations up to equivalence is algebraic.

  Finally, we can remove the trivial deformations (which correspond to a
  single closed point of the stack) and quotient out by $\G_m$ to obtain
  the desired moduli stack, which is therefore still algebraic.
\end{proof}

\begin{rem}
  Note that in the above argument, if the degrees of the generators are
  $\le$ the degrees of the relations, then the moduli stack is a quotient
  of a scheme by a group of the form $\G_m\ltimes U$ where $U$ is
  unipotent.  The above induction is needed in general, however, since the
  coordinatization of the set of lifts of a generator depends on the
  relations of lower degree, which in turn depend on the lifts of
  generators of up to that degree.
\end{rem}

One caution here is that this argument is {\em highly} nonconstructive;
indeed, it seems likely that even the question of whether $B$ admits a
nontrivial filtered deformation is undecidable in general, even if we are
given a finite presentation of $B$.  (The situation improves dramatically
if we are also given generators of the module of syzygies.)

This moduli stack is nicest under the following condition.

\begin{prop}
  Suppose that the finitely presented graded algebra $B$ admits no
  homogeneous derivations of negative degree.  Then the moduli stack of
  filtered deformations of $B$ is a closed substack in a weighted
  projective space.
\end{prop}

\begin{proof}
  For each positive integer $d$, consider the effect of adding terms
  of degree $\deg(x_i)-d$ to each lifted generator $\hat{x}_i$.
  The resulting effect on the relations is trivial until one reaches
  degree $\deg(\rho_i)-d$, where it induces a linear transformation
  from the space of deformations of generators to the space of
  deformations of relations.  This linear transformation is independent
  of the choice of filtered deformation of $B$, and may thus be
  computed using the trivial filtered deformation.  In particular,
  we find that if we view the linear transformation as its own
  derivative, then the kernel is precisely the space of degree $-d$
  derivations of $B$, and thus by hypothesis the linear transformation
  is injective.

  We now choose for each positive integer $d$ a complement to the
  image of the linear transformation in the space of degree $-d$
  deformations of relations.  We then claim that every filtered
  deformation is uniquely equivalent to one in which the deformations
  of every degree lie in the chosen complement.  This may be seen
  by induction on $d$: if the deformations of relative degree $>-d$
  all lie in their respective subspaces, then (by injectivity and
  complementarity) there is a unique way to shift the generators
  by terms of relative degree $-d$ to make the relations lie in the
  complementary subspace in degree $-d$.  Since this does not affect
  the terms of higher degree, the result follows.

  It follows that the moduli stack of filtered deformations up to
  equivalence is a scheme, and is in fact $\Spec(M)$ for some (commutative)
  graded algebra $M$ (with grading coming from the $\G_m$ action).  The
  desired moduli stack is then obtained by removing the origin (the unique
  point corresponding to the trivial deformation) and quotienting by
  $\G_m$, and is thus closed in the weighted projective space corresponding
  to the homogeneous generators of $M$.
\end{proof}

\begin{rem}
  Here we view weighted projective spaces as stacks; if all coordinates of
  degree prime to $l$ vanish, then the result will be an orbifold point
  with stabilizer containing $\mu_l$.  We can recover the usual scheme
  versions by taking coarse moduli spaces.  Note that in either
  interpretation, weighted projective spaces are proper, and thus so are
  their closed substacks/subschemes.
\end{rem}

An important special case in which it is straightforward to control
derivations is when $B$ is a commutative homogeneous coordinate ring.

\begin{lem}\label{lem:derivs_on_hcr}
  Let $X$ be an irreducible smooth projective scheme of positive dimension
  with ample line bundle ${\cal L}$, with associated homogeneous coordinate
  ring $B:=\bigoplus_{i\ge 0}\Gamma(X;{\cal L}^i)$.  For any integer $k>0$,
  there is a short exact sequence
  \[
  0\to \Omega^k_X\to V_k\to \Omega^{k-1}_X\to 0
  \]
  such that there is an isomorphism
  \[
  \Der^k(B)\cong \bigoplus_{i\in \Z} \Hom(V_k,{\cal L}^i)
  \]
  of graded vector spaces.
\end{lem}

\begin{proof}
  If we remove the cone point from $\Spec(B)$, the resulting scheme is a
  $\G_m$-torsor $X^+$ over $X$, so that an element of $\Der^k(B)$ thus
  determines an element of $\Der^k(\sO_{X^+})$.  This is, in fact, a
  $\G_m$-equivariant isomorphism
  \[
  \Der^k(B)\cong \Der^k(\sO_{X^+}).
  \]
  Indeed, injectivity follows from the fact that $X^+$ still contains the
  generic point of $\Spec(B)$, while surjectivity follows by observing that
  any holomorphic section of $\sO_{X^+}$ decomposes via the $\G_m$ action
  as a sum of homogeneous elements, each of which induces a section of the
  corresponding power of ${\cal L}$ and thus gives a homogeneous element of
  $B$.

  Since $\pi:X^+\to X$ is smooth of dimension 1, there is a short exact
  sequence
  \[
  0\to \pi^*\Omega_X\to \Omega_{X^+}\to \Omega_{X^+/X}\to 0,
  \]
  with $\Omega_{X^+/X}$ a line bundle.  This line bundle is in fact
  canonically trivial, via the map
  \[
  \Omega_{X^+/X}\otimes \Lie(\G_m)\to \sO_{X^+}
  \]
  coming from the $\G_m$-torsor structure.  We further find that the
  induced map
  \[
  \Gamma(X^+;\sO_{X^+})\to H^1(X^+;\pi^*\Omega_X)
  \cong H^1(X;\pi_*\pi^*\Omega_X)
  \cong H^1(X;\bigoplus_i \Omega_X\otimes {\cal L}^i)
  \]
  is homogeneous of degree 0, and thus corresponds to an element of
  $H^1(X;\Omega_X)$.  In other words, there is a short exact sequence
  \[
  0\to \Omega_X\to V_1\to \sO_X\to 0
  \]
  such that $\Omega_{X^+}\cong \pi^*V_1$.  Taking exterior powers
  gives a similar exact sequence
  \[
  0\to \Omega^k_X\to V_k\to \Omega^{k-1}_X\to 0
  \]
  for $V_k:=\wedge^k V_1$.  We then find
  \begin{align}
  \Der^k(\sO_{X^+})
  &\cong
  \Gamma(X^+;\sHom(\pi^* V_k,\sO_{X^+}))\notag\\
  &\cong
  \Gamma(X;\pi_*\sHom(\pi^* V_k,\pi^*\sO_X))\notag\\
  &\cong
  \Gamma(X;\sHom(V_k,\pi_*\pi^*\sO_X))\notag\\
  &\cong
  \bigoplus_{i\in \Z} \Hom(V_k,{\cal L}^i)
  \end{align}
  as required.
\end{proof}

\begin{rem}
  Note that in this generality, $B$ can fail to be finitely generated.
\end{rem}

\begin{cor}
  Let $C$ be a smooth projective curve of positive genus, ${\cal L}$ an
  ample line bundle on $C$, and $B$ the corresponding homogeneous
  coordinate ring.  Then any homogeneous element of $\Der^k(B)$ of negative
  degree vanishes.
\end{cor}

\begin{proof}
  For $d>0$, the space of homogeneous elements of $\Der^k(B)$ of degree
  $-d$ is
  \[
  \Hom(V_k,{\cal L}^{-d})
  \]
  where $V_k$ is the sheaf from the Lemma, so that for $k=1$ there
  is an exact sequence
  \[
  0\to \Hom(\sO_C,{\cal L}^{-d})\to \Hom(V_1,{\cal L}^{-d})\to
  \Hom(\Omega_C,{\cal L}^{-d}),
  \]
  with both spaces on either side of the one of interest vanishing by degree
  considerations.  Vanishing for $k=2$ follows similarly from
  $V_2\cong \Omega_C$, while for $k>2$, $V_k=0$.
\end{proof}

\begin{cor}
  With $B$ as above, the moduli stack of filtered deformations is a closed
  substack of a weighted projective space, and every filtered deformation
  of $B$ is commutative.
\end{cor}

\begin{proof}
  The absence of derivations of negative degree gives the first claim.  For
  the second claim, suppose $A$ were a noncommutative filtered deformation
  of $B$, and view the Rees algebra of $A$ as a (flat) $k[t]$-module.  This
  gives a 1-parameter family of noncommutative algebras degenerating when
  $t=0$ to the commutative algebra $B$, and thus induces a nonzero Poisson
  bracket on $B$.  Since computing the Poisson bracket involves dividing
  commutators by an appropriate positive power of $t$, this Poisson bracket
  on $B$ is homogeneous of negative degree.  Since a Poisson bracket is in
  particular a $2$-derivation, it must therefore vanish, giving a
  contradiction.
\end{proof}

This makes it straightforward to compute the moduli space of filtered
deformations of commutative elliptic algebras.

\begin{eg}
Consider the case of an elliptic curve $E$, with ample line bundle ${\cal
  L}=\sO_E(\infty)$ of degree 1 coming from the point at infinity.  The
corresponding commutative graded algebra is generated by elements $z$, $x$,
$y$ of degree $1$, $2$, $3$ respectively, satisfying a single relation of
degree 6:
\[
y^2 + a_1xyz + a_3yz^3 = x^3+a_2 x^2z^2 + a_4 xz^4 + a_6z^6.
\]
By the corollary, any filtered deformation of this algebra is commutative,
and thus simply adds lower degree terms to the single relation; and since
there is a single commutative relation, there are no further conditions to
give a flat deformation.  To get equivalence, we must quotient by changes
of variable of the form
\[
(y,x,z)\mapsto (y+b_1 x + c_1 z^2 + b_2 z + b_3,x+d_1 z + c_2, z+e_1).
\]
This group acts faithfully, and thus to obtain the structure of the moduli
stack, we need simply remove four generators of degree 1, two of degree 2,
and one of degree 3.  We thus see that in this case, the moduli stack is
the weighted projective space $\P^{[122334456]}$ (i.e., with generators of
degrees $1,2,2,3,3,4,4,5,6$).  (These are the standard labels of the affine
Dynkin diagram of type $E_8$, as explained below).

A similar calculation works if we instead take ${\cal L}$ to have degree
$2$, $3$, or $4$; in the first two cases, the cone is a hypersurface, while
for $4$, it is a complete intersection, and thus in no case are there any
syzygies.  The result in each case is thus again a weighted projective
space ($\P^{[11222334]}$ ($E_7$); $\P^{[1112223]}$ ($E_6$); and
$\P^{[111122]}$ ($D_5$) respectively).

For higher degrees, there are syzygies, so some actual calculation is
required, but one can still easily compute the corresponding moduli space
since there is no need to consider noncommutative deformations.  One finds
that the space is $\P^4$ for degree 5 and $\P^1\times \P^2$ for degree 6.
For degrees 7,8,9, the structure is somewhat more complicated (but see
below), while for degree $>9$, there are no nontrivial deformations
\cite[\S9]{PinkhamHC:1974}.  For degree $\ge 3$, it was shown op.~cit.~that
any nontrivial deformation is rational and anticanonical (albeit possibly
singular), i.e., a del Pezzo surface of the given degree.  For degrees $3$
and $4$, this can also be seen by comparing the moduli stack as computed
above to the substack coming from del Pezzo surfaces, and the same argument
establishes the corresponding claim for degrees $1$ and $2$.  (This can
also be seen more directly: that the surface is Fano follows by adjunction
from the fact that the curve at infinity has genus 1; that it is rational
reduces (via successive blowups of points at infinity) to the degree 1
case, where it follows from the theory of elliptic surfaces.)
\end{eg}

\medskip

Note that the hypotheses above are quite strong,
and in particular exclude the cases in which $E$ becomes singular.  These
are, of course, still quite interesting, but harder to deal with
geometrically.  (By contrast, it is often quite straightforward to write
down explicit presentations in those cases!)  As an example of what can go
wrong, one has the following family in which the commutative deformation
fails to be a del Pezzo surface in the usual sense (and one obtains
deformations for ${\cal L}$ of arbitrarily large degree).

\begin{eg}\label{eg:weird_nodal}
  For $d\ge 3$, let $Z$ be the weighted projective space with generators $x$,
  $y$, $z$, $w$ of degrees $d-2$, $d-2$, $1$, $1$ respectively, and let $X_{ab}$
  be the family of hypersurfaces of degree $d$ with equation
  \[
  zw x + (az^2+bw^2)y = z^d+w^d.
  \]
  The curve $y=0$ in $X_{ab}$ is isomorphic to the standard nodal cubic,
  and its image under the embedding corresponding to the invertible sheaf
  $O_Z(d-2)$ is a degree $d$ model in $\P^{d-1}$; to be precise, it is the
  image of $\P^1$ under the map $(z^d+w^d,zw^{d-1},\dots,z^{d-1}w)$ that
  identifies the points $(0,1)$ and $(1,0)$.  In other words, the $(d-2)$nd
  Veronese of the homogeneous coordinate ring of $X_{ab}$ is the Rees
  algebra of a filtered deformation of that nodal cubic.  It is clear that
  this filtered deformation is nontrivial unless $a=b=0$, and thus gives us
  a $\P^1$ of such deformations for all $d$.  When $d=10$, one can verify
  by a direct Magma computation that the reduced subscheme of the full
  deformation scheme is $\P^1$, and thus all deformations arise in this
  way.  The presence of such deformations arises from the fact that the
  cone $X_{00}$ is not normal, and thus the filtered deformations
  themselves need not be normal (and, indeed, are not: $X_{ab}$ is singular
  along the curve $z=w=0$).  These are examples of the non-normal del Pezzo
  surfaces studied in \cite{ReidM:1994}.
\end{eg}

\medskip

Now, with $E$ still an elliptic curve and ${\cal L}$ an ample line bundle,
we can also obtain noncommutative graded algebras by taking
\[
B=\bigoplus_i \Hom(V,V\otimes {\cal L}^i)
\]
for $V$ a vector bundle on $E$.  For simplicity, we assume that $B$ is
nonnegatively graded; we will see below that this is not truly a
constraint.  This algebra is finite over its center, the usual homogeneous
coordinate ring, and is in fact the ring of global sections of a sheaf of
(Azumaya) algebras ${\cal B}:=\sEnd(\pi^*V)$ on the $\G_m$-torsor $E^+$
associated to ${\cal L}$.

For any associative algebra $A$, there is a natural map $\Der(A)\to
\Der(Z(A))$ coming from the fact that any automorphism (and thus any
infinitesimal automorphism) preserves the center.

\begin{lem}\label{lem:all_derivs_inner}
  A derivation of $B$ that is trivial on $Z(B)$ is inner.
\end{lem}

\begin{proof}
  A derivation of $B$ which is trivial over $Z(B)$ may be viewed as an
  $\sO_{\Spec(Z(B))}$-linear derivation on the corresponding sheaf of
  algebras, and thus induces an $\sO_{E^+}$-linear derivation on ${\cal B}$
  after removing the cone point, or equivalently a holomorphic family of
  derivations on the fibers of ${\cal B}$.  Since those fibers are matrix
  algebras over fields, their derivations are necessarily inner.  In
  particular, we may locally represent the given derivation as the
  commutator with a {\em traceless} matrix, and since this representation
  is canonical, it extends to a global such representation.  As in the
  proof of Lemma \ref{lem:derivs_on_hcr}, a global section of ${\cal B}$ is
  an element of $B$, and thus the original derivation of $B$ is inner.
\end{proof}

\begin{rem}
  This argument fails in finite characteristic, since it uses the
  fact that the map $\Mat_n\to \Inn\Mat_n$ is split by the trace,
  which fails when $n$ is a multiple of the characteristic.  Of
  course, any global splitting will do, so it in particular suffices
  to have {\em some} idempotent in $\End(V)$ with invertible trace,
  or equivalently for $V$ to have a summand with invertible rank.
  (Though even this appears to be stronger than needed for the
  conclusion to hold.)
\end{rem}

\begin{rem}
  An equivalent statement is that the natural map $\HH^1(B)\to \HH^1(Z(B))$
  is injective.
\end{rem}

\begin{cor}\label{cor:no_neg_derivs}
  The algebra $B$ has no homogeneous derivations of negative degree.
\end{cor}

\begin{proof}
  We have already seen that $Z(B)$ has no such derivations, and thus such a
  derivation of $B$ is trivial on $Z(B)$, so inner, and thus itself
  vanishes since $B$ is nonnegatively graded.
\end{proof}

\begin{rem}
  Note that $B$ is Noetherian since ${\cal B}$ is Noetherian, and is in
  particular finitely presented.  We may thus conclude that the moduli
  stack of filtered deformations of $B$ is a closed substack of a weighted
  projective space.
\end{rem}

The fact that filtered deformations remain commutative when $V=\sO_E$ also
has an analogue: the center of any filtered deformation of $B$ is a
filtered deformation of $Z(B)$.  To see this, we use ideas of Hayashi
\cite{HayashiT:1988} (see also \cite[\S 1]{KacV/RadulA:2000} for an
explicit discussion), in somewhat more formal terms.  We again view the
Rees algebra $A^+$ as a flat $k[t]$-algebra.  Call an element $a\in A^+$ of
this algebra {\em $m$-central} if $[a,A^+]\subset t^m A^+$, or equivalently
if its image in $A^+/t^m A^+$ is central; in particular, an element is
$1$-central iff its image in $B$ is central.  Now, given an $m$-central
element $a$, we may define a derivation $D_{a,m}(x):=t^{-m}[a,x]$ on $A^+$.
Since this is $k[t]$-linear, the image of $D_{a,m}(x)$ mod $t$ depends only
on the image of $x$ mod $t$, so that $D_{a,m}$ induces a derivation of $B$.
We also find
\[
D_{a+t^m b,m}(x) = D_{a,m}(x)+[b,x],
\]
and thus see that there is an $m+1$-central element of the form $a+t^m b$
iff $D_{a,m}$ is an inner derivation.  In other words, the obstruction to
lifting an element of $Z(A^+/t^m A^+)$ to an element of $Z(A^+/t^{m+1}A^+)$
is an element of $\Der(B)/\Inn(B)=\HH^1(B)$.  Note that since the product
of $m$-central elements is $m$-central and
$D_{ab,m}(x)=aD_{b,m}(x)+D_{a,m}(x)b$, the obstruction map
\[
Z(A^+/t^m A^+)\to \HH^1(B)
\]
is itself a derivation.

Now, let $m>0$ be an integer, and suppose that every central element of
$A^+/t^{m-1}A^+$ lifts to an $m$-central element of $A$ (this is vacuously
true for $m=1$).  If $a_1$, $a_2$ are $m$-central elements that reduce to
the same element of $B$, then
\[
D_{a_1,m}(x)-D_{a_2,m}(x) = t^{-m} [a_1-a_2,x] =
t^{1-m}[t^{-1}(a_1-a_2),x]
=
D_{t^{-1}(a_1-a_2),m-1}(x).
\]
Now, $t^{-1}(a_1-a_2)$ is $(m-1)$-central, and thus agrees modulo $t^{m-1}$
with an $m$-central element, so that $D_{t^{-1}(a_1-a_2),m-1}(x)$ is inner,
implying that
\[
D_{a_1,m}(x)-D_{a_2,m}(x)
\]
is inner.  We thus find that the obstruction map factors through a
derivation
\[
Z(B)\to \HH^1(B),
\]
which since it divides by $t^m$ is homogeneous of degree $-m$.  (As pointed
out in \cite[Ex. 2.17(4)]{EtingofP/SchedlerT:2010}, this map appears in a
spectral sequence for computing the Hochschild cohomology of the filtered
deformation.)

For any $a\in Z(B)$, the induced outer derivation of $B$ induces a
derivation of $Z(B)$, giving a 2-derivation of $Z(B)$ of negative degree,
which vanishes.  Since a derivation of $B$ vanishing on $Z(B)$ is inner, we
conclude that the obstruction map for degree $m$ vanishes.  We thus find
the following.

\begin{prop}\label{prop:center_lifts}
  Let $E$ be an elliptic curve, ${\cal L}$ an ample line bundle on $E$, and
  $V$ a vector bundle such that $\Hom(V,V\otimes {\cal L}^{-1})=0$.  The
  center of any filtered deformation $A$ of $\bigoplus_i \Hom(V,V\otimes
  {\cal L}^i)$ is a filtered deformation of $\bigoplus_i \Gamma({\cal
    L}^i)$.
\end{prop}

\begin{proof}
  Let $\bar{x}\in Z(B)$ be a homogeneous element, and consider the
  sequence of $m$-central lifts $x_1$, $x_2$, $x_3$,\dots, $x_m$,\dots
  guaranteed by the above discussion.  Then $x_{i+1}-x_i\in t^i
  A^+$, so that the image of $x_{i+1}-x_i$ in $A$ (i.e., setting
  $t=1$) has degree $\le \deg(\bar{x})-i$.  In particular, we find
  that any sequence of $m$-central lifts is eventually constant,
  and thus gives a lift of $\bar{x}$ which is $m$-central for all
  $m$, thus central.  This shows that $\gr Z(A)\supset Z(B)$, while
  the other inclusion is trivial, so that the claim follows.
\end{proof}

As with the commutative case, one can use these results to enable explicit
computations of the moduli stacks.  One issue is that it is more difficult
to write down an explicit presentation of $B$, especially if one wishes to
deal with all $E$ at once.  A useful trick for this purpose involves the
observation that $B$ itself has regular central elements of degree 1, so
may be viewed as a filtered deformation, and the resulting algebra turns
out to be independent of $E$.  To be precise, if $w\in Z(B)$ is a degree
$1$ element corresponding to $d=\deg({\cal L})$ distinct points of $E$,
i.e., corresponding to an isomorphism ${\cal L}\cong {\cal L}(D)$ for a
reduced divisor $D$, then the quotient $C:=B/\langle w\rangle$ has Hilbert
series of the form $\dim(B_0)+(n-\dim(B_0))t+nt^2+nt^3+\cdots$ where
$n=d\rank\sEnd(V)$.  If $\tilde{C}$ denotes the algebra $\Gamma(D;\sEnd(V))$,
then every homogeneous component of $C$ may be identified with a subspace
of $\tilde{C}$, and thus in particular dimension considerations tell us
that $C$ agrees with $\tilde{C}$ in all degrees $>1$.  (Note that the
identification with $\tilde{C}$ depends on a choice of trivialization of
${\cal L}(D)|_D$, or in other words on a system of uniformizers at the
points of $D$.  By adjunction, such a system is determined by a choice of
holomorphic differential $\omega$ on $E$, so there is an essentially
canonical choice of identification.)  In degree $0$, we of course have $C_0
= B_0 = \End(V)\subset \Gamma(D;\sEnd(V))$, so it remains only to determine
$C_1$.  There is a natural pairing $C_0\times C_1\to \Gamma({\cal L}(D))$
given by $(x,y)\mapsto \Tr(xy)$, which agrees with the
$\Gamma(\sO_D)$-valued trace pairing on $\tilde{C}$ if we take residues
(with the same choice of holomorphic differential).  Since the value is a
global function, the sum of residues vanishes, and thus any element of
$C_1$ must be orthogonal to $C_0$ under the sum-of-traces pairing on
$\tilde{C}$.  This cuts out a space of the correct dimension, so completely
determines $C_1$.

In other words, $C$ is determined by the algebra $\tilde{C}$ and its
subalgebra $B_0\subset \tilde{C}$:
\[
C_n =
\begin{cases} B_0 & n=0\\
    \{x:x\in \tilde{C}\mid \Tr(x B_0)=0\} & n=1\\ \tilde{C} & n>1
\end{cases}.
\]
Here $\tilde{C}$ is (geometrically) just a sum of $d$ copies of
$\Mat_{\dim(V)}$, while $B_0=\End(V)$ depends only on the gross structure
of $V$ (the slopes and multiplicities of the indecomposable summands, along
with which summands have isomorphic stable constituents).  It is then
straightforward to determine a presentation of $C$, either by using the
results for $B$ below or by using the fact that $C$ has regular central
elements (and in particular has a canonical regular central element of
degree 2, namely $1\in \tilde{C}=C_2$), and the quotient by such an element
is finite-dimensional.

\begin{eg}\label{eg:mu0r1}
  It is interesting to consider the commutative case from this perspective.
  Here the algebra $\tilde{C}$ is just $k^d$, while $B_0$ is the diagonal
  copy of $k$, so that for $d>2$, $C$ is the embedded coordinate ring of
  $d$ points in general position in $\P^{d-2}$.  (For $d=1$, it is
  $k[y,x]/(y^2-x^3)$ with $\deg(x)=2$, $\deg(y)=3$, and for $d=2$, it is
  $k[y,x]/(y^2-x^2y)$ with $\deg(x)=1$, $\deg(y)=2$.)  The commutative
  filtered deformations of such algebras were studied in
  \cite{LekiliY/PolishchukA:2019}, and in particular shown to always give a
  (possibly degenerate) genus 1 curve on which the marked degree $d$
  divisor is ample.  Moreover, the structure of the moduli stack was
  determined in low degree, giving the weighted projective spaces
  $\P^{[46]},\P^{[234]},\P^{[1223]},\P^{[11122]}$ in degrees $1,2,3,4$,
  $\P^4$ in degree $5$, and the Grassmannian $\Gr(2,5)$ in degree $6$.  Not
  only does this let us write down general versions of $B$ in low degree,
  but it in fact lets us compute the moduli stack of deformations of $B$: a
  filtered deformation of $B$ is essentially just a line through the
  corresponding point of the moduli stack of filtered deformations of $C$.
  (The description is somewhat more involved when the moduli space for $C$
  has orbifold points: in general, it corresponds to a homomorphism from
  the graded algebra parametrizing deformations of $C$ to $k[x,y]$ such
  that the induced rational map from $\P^1$ is nonconstant and $\infty$
  maps to the point associated to $B$.)  This recovers the results for
  $d\le 4$ above, as well as giving a more conceptual argument for $d=5,6$.
  (There are issues for $d=1,2$ in small characteristic, but these can be
  resolved by waiting until the last step to kill the unipotent freedom.)
  Interestingly, for $d>1$, if we quotient by a further regular element of
  degree 1, the deformation spaces remain nice:
  $\P^{[4]},\P^{[23]},\P^{[1122]},\P^6$ respectively for $2\le d\le 5$,
  while for $d=6$ the stack of deformations is the isotropic Grassmannian
  $\OGr^+(5,10)$.  (The arguments are, sadly, not conceptual, but simply
  involve explicitly computing the deformation space, and then, after
  guessing based on the Hilbert series, computing the explicit isomorphism
  for $d=6$.  It would be particularly interesting to have a conceptual
  construction of the resulting sheaf of filtered algebras over
  $\OGr^+(5,10)$!)
\end{eg}

\begin{eg}\label{eg:d1stab}
  Consider the case that $V$ is a stable bundle of rank $2$ and degree $1$
  and the ample bundle is ${\cal L}(\infty)$.  Then $B$ has Hilbert series
  $1+4t/(1-t)^2$, while $Z(B)$ has Hilbert series $1+t/(1-t)^2$.  In this
  case, we get a particularly nice algebra $D$ if we quotient $C$ by the
  regular central element of degree $2$: $D$ is simply the exterior algebra
  in three generators!  The moduli stack of filtered deformations of that
  exterior algebra is then straightforward to compute, and one finds that
  the space of filtered deformations is the weighted projective space
  $\P^{[222222]}$: any filtered deformation of the exterior algebra in
  three generators is equivalent to a Clifford algebra.  (There is a $\P^5$
  worth of Clifford algebras, but Clifford algebras have an additional
  automorphism acting as $-1\in \G_m$ on the associated graded.)  Since a
  filtered deformation $A$ of $B$ has center a filtered deformation of
  $Z(B)$, the quotient of the Rees algebra $A^+$ by its two degree 1
  elements and one degree 2 element is again the exterior algebra, and thus
  we may view the filtered deformation of $B$ as a family of filtered
  deformations of the exterior algebra.  We thus see that it is given by a
  map from the weighted projective space $\P^{[112]}$ to the moduli stack
  of Clifford algebras, or more precisely by a collection of six
  homogeneous polynomials of degree 2 in those generators.  Of course, the
  specializations of those polynomials to $t=0$ are determined by $B$ (that
  is, together with the specific choices of degree 1 and 2 elements), and
  thus there are only two free parameters for each of the polynomials, one
  of degree 1 and one of degree 2.  Of those, one each is eliminated by
  translations of the degree 2 central element of $B$, and one more of
  degree 1 is eliminated by the translations of the degree 1 central
  element.  We thus conclude that for {\em any} elliptic curve and
  semistable bundle of rank 2, degree 1, the moduli stack of filtered
  deformations of the corresponding graded algebra is the weighted
  projective space $\P^{[111122222]}$ ($D_8$).  (One caveat: the above
  calculation assumes $2$ is invertible, both because the bundle has even
  rank, breaking the proof of Lemma \ref{lem:all_derivs_inner}, and because
  the exterior algebra over $\bar\F_2$ has derivations of negative degree.
  One can, however, verify by a more complicated calculation that the
  claims remain true in characteristic 2.)
\end{eg}

\begin{rem}
  For $V$ stable of slope $1/2$ and ${\cal L}$ of degree $d$, the algebra
  $C$ is equal to $\Mat_2^d$ in degree $>1$, and the subspace of $d$-tuples
  of total trace $0$ in degree $1$, and in particular its center is the
  algebra considered in Example \ref{eg:mu0r1}.  There is an action of
  $\GL_2^d$ on the corresponding deformation space (of filtered
  deformations respecting the center), and we find (computationally) the
  following descriptions after taking that into account:
  \begin{align}
    d=1:{} & S^4(V_1)/(\G_m/\langle \pm 1\rangle)\\
    d=2:{} & (S^2(V_1)\otimes S^2(V_2))/\G_m\\
    d=3:{} & \!\Gr(2,V_1\otimes V_2\otimes V_3),
  \end{align}
  where $V_1$,\dots are the fundamental representations of the respective
  copies of $\GL_2$.  In each case, this presumably has a geometric modular
  interpretation in terms of rank 2 vector bundles (or torsion-free
  sheaves) over the center, but it appears that there is an alternate
  interpretation in terms of hyperelliptic curves of genus 1 with
  additional structure.  Indeed, for $d=1$, $d=2$, the given representation
  of $\GL_2^d$ is precisely the representation used in
  \cite{BhargavaM/HoW:2016} in connection with $2$-Selmer groups of
  elliptic curves with $d-1$ marked points, while for $d=3$, the space
  considered there is $V_1\otimes V_2\otimes V_3\otimes V_4$, and thus our
  space is the GIT quotient of their space by the ``extra'' copy of
  $\GL_2$.  Similarly, for the case $\rank(V)=3$, $d=1$, the deformation
  space is two disjoint copies of $S^3(V_1)/\G_m$ (with $V_1$ the
  fundamental representation of $\GL_3$), corresponding to cubic plane
  curves (and thus related to $3$-Selmer groups).  (One obtains two copies
  as there are two possible slopes modulo 1 with denominator $3$.)  The
  remaining representations considered in \cite{BhargavaM/HoW:2016} thus
  suggest three more deformation spaces associated to stable bundles:
  \begin{align}
    \rank(V)=3,\, d=2:{} &\!\Gr(3,V_1\otimes V_2)\\
    \rank(V)=4,\, d=1:{} &\!\Gr(2,S^2(V_1))\\
    \rank(V)=5,\, d=1:{} &\!\Gr(5,S^2(V_5));
  \end{align}
  in the first two cases, one expects two components of that form, while in
  the third case, one expects a total of four components, two of which are
  of the given form.  Note that in each case, the resulting prediction for
  the deformation spaces of the resulting elliptic algebras agrees with the
  computation in Section 5 below.  It would be interesting to know if there
  are any other instances with nice deformation space.
\end{rem}

\begin{eg}\label{eg:d1r2}
  Similarly, if $V$ is a sum of two nonisomorphic degree 0 line bundles,
  then $C$ is generated (over $B_0\cong k^2$) by elements $F$ of degree $1$
  and $G$ of degree 2 interacting with the two idempotents of $B_0$ as $e_1
  F = F e_2$, $e_1 G = G e_2$, and satisfying relations $FG=GF$, $G^2=F^4$.
  (These both correspond to the matrix $\begin{pmatrix} 0 & 1 \\ 1 &
    0\end{pmatrix}\in \tilde{C}$.)  Any filtered deformation is equivalent
    to one in which the relations involving the idempotents remain
    unchanged, and two such deformations are equivalent iff they are
    related by a shift $G\mapsto G+b_1 F$.  (We can only add lower terms
    that interact the same way with the idempotents.)  An explicit
    calculation with syzygies of low degree tells us that the universal
    deformation of this algebra has the form
\begin{align}
  FG-GF &= a_3(e_1-e_2),\\
  G^2 &= F^4 + a_2 F^2 + a_4.
\end{align}
  In particular, $B$ corresponds to a point of the weighted projective
  space $\P^{[234]}$, and a filtered deformation of $B$ extends that point
  to a map from $\P^1$, from which it follows easily that the moduli stack
  of filtered deformations of $B$ is a weighted projective space with
  degrees $\P^{[11222334]}$ (again $E_7$).  (This calculation again assumes
  $2$ is invertible, but can be extended to work over $\Z$ by waiting until
  the last step to kill the unipotent symmetry.)
\end{eg}

\begin{rem}
  Just as $A$ may be viewed as a line of filtered deformations of $C$, we
  may obtain filtered deformations of $A^+$ as either a line of filtered
  deformations of $B$ or a plane of filtered deformations of $C$.  (This
  imposes the additional condition that $Z(A^+)$ lifts, which is not
  necessarily automatic.)  This at the very least gives rise to a family of
  (algebras on) Fano $3$-folds (embedded in such a way that $-K_X=\sO(2)$),
  and we could hope to obtain noncommutative Fano 3-folds by replacing $B$
  by a twisted elliptic algebra.  (Indeed, the Sklyanin algebra
  \cite{SklyaninEK:1982} may be obtained by a variation on this idea, see
  below.)  This can, of course, be iterated as long as enough degrees of
  freedom remain to allow the deformations to be nontrivial.  When the
  moduli stack of deformations of $B$ is a weighted projective space, a
  particularly interesting variation is the {\em universal} deformation,
  in which one adjoins every generator of the weighted projective space to
  $B$.  In the commutative case, the resulting scheme is independent of
  $E$, and is isomorphic to $\P^{[11222333445]},
  \P^{[1111222233]},\P^{[111111222]},\P^7,\Gr(2,5)$ for $1\le d\le 5$, so
  that we obtain noncommutative deformations of those spaces by twisting
  $B$.  (It would be interesting to see if the other cases with nice moduli
  stacks give rise to recognizable algebras in the untwisted case.)
\end{rem}

\medskip

The claims about the center can be further generalized.  Consider an
elliptic algebra of the form $B_{V,\Psi}:=\bigoplus_i \Hom(V,\Psi^i V)$,
where $\Psi$ is an autoequivalence of $\coh(E)$ of the form $M\mapsto {\cal
  L}\otimes (\phi^{-1})^*M$ with $\deg({\cal L})>0$.  If $\phi$ has finite
order, so that $\Psi^d$ is twisting by a line bundle, then $B$ has a large
center, namely the homogeneous coordinate ring of the quotient curve
$E/\langle \phi\rangle$ relative to the norm of ${\cal L}$, with degrees
multiplied by the order of $\phi$.  When $\phi$ is a translation by a
$d$-torsion point $q$, $E':=E/\langle \phi\rangle$ is again an elliptic
curve, and thus its homogeneous coordinate ring has no derivations of
negative degree.  (Note that this could fail when $\phi$ is not a
translation, since then $E'$ has genus 0.  It would be interesting to
understand the filtered deformations in such cases, but they fall outside
the scope of the present work.)

Let $E^{\prime +}$ be the $\G_m$-torsor associated to the norm line bundle
$N_{E/E'}({\cal L})$ on $E'$, modified by pulling back the action of $\G_m$
through the $d$-th power map.  (This is to correct for the fact that the
degrees in $Z(B)$ are all multiples of $d$, and in particular the degree
$dm$ elements are sections of $N_{E/E'}({\cal L})^m$.)  Then $B$ is the
ring of global sections of a corresponding sheaf of algebras ${\cal B}$ on
$E^{\prime +}$.  This sheaf is $\G_m$-equivariant in a somewhat complicated
way, but the (geometric) fibers are again just matrix algebras.  Indeed,
the $\mu_d$-invariant subalgebra of ${\cal B}$ is the direct image of the
endomorphism ring of a vector bundle on a $\G_m$-torsor over $E$, so is
locally a matrix algebra over a cyclic degree $d$ cover of $E^{\prime +}$.
Relative to this interpretation, the other congruence classes correspond to
spaces of linear transformations which are semilinear w.r.to the
appropriate power of $\phi$.  We thus find that ${\cal B}$ is locally a
matrix algebra over a {\em cyclic} algebra, and is thus an Azumaya algebra.

It is then straightforward to extend the above arguments to show that $B$
has no negative derivations.  We can actually extend this to the case that
$\tau$ has infinite order, but this will require us to control derivations
in finite characteristic, where the argument of Lemma
\ref{lem:all_derivs_inner} can fail.  Luckily, if we modify the conclusion
somewhat, we obtain a result that holds in general.

\begin{lem}\label{lem:derivs_cong}
  Let $B$ be an elliptic algebra (possibly in finite characteristic),
  twisted by a torsion point of exact order $d$.  Then any derivation of
  $B$ of degree not a multiple of $d$ is inner.
\end{lem}

\begin{proof}
  Let $D$ be a derivation of $B$ of degree $\delta$, with $\delta$ not a
  multiple of $d$.  Since every element of $Z(B)$ has degree a multiple of
  $d$, it follows immediately that $D$ is trivial on $Z(B)$.  We may then
  argue as in the proof of Lemma \ref{lem:all_derivs_inner} that $D$ is
  locally inner; in other words, there is a collection of representations
  $D=\ad(z_i^{-1}x_i)$ with $z_i\in Z(B)$, $x_i\in B$, such that the open
  subsets $U_i:=[z_i\ne 0]$ cover $E'$.

  Now, since $\ad(z_i^{-1}x_i)=\ad(z_j^{-1}x_j)$, the system of pairs must
  also satisfy $z_jx_i-z_ix_j\in Z(B)$ for all $i,j$.  Since
  $\deg(z_i^{-1}x_i)=\delta$, we see that
  \[
  \deg(z_jx_i-z_ix_j)\equiv \delta\pmod{d},
  \]
  which (again using the fact that elements $Z(B)$ have degree a multiple
  of $d$) forces $z_jx_i=z_ix_j$.  But this implies that the local
  representations glue to give a global representation, so that $D$ is
  inner as required.
\end{proof}

\begin{lem}\label{lem:no_neg_derivs_general}
  Let $B$ an elliptic algebra twisted by a translation $\tau$.  If $B$ has
  a nonzero derivation of degree $-d$ with $d>0$, then $\tau^d=1$.
\end{lem}

\begin{proof}
  Suppose otherwise, so that $B$ has such a derivation, but $\tau^d\ne 1$.
  If the field of definition $K$ is infinitely generated, then we can
  replace it by the field generated by (a) the coefficients of some
  Weierstrass equation for $E$, (b) the coordinates of the point by which
  $\tau$ is a translation, and (c) the coordinates of the points
  corresponding to the determinants of stable constituents of $V$, so that
  we may assume $K$ is finitely generated.  (See \cite{AtiyahMF:1957} or
  the next section for a discussion of the structure of $V$.)

  If $K$ is finite, then $\tau$ has finite order, with $d$ not a multiple
  of $\ord(\tau)$, and thus Lemma \ref{lem:derivs_cong} gives a
  contradiction.  Otherwise, let $v$ be a valuation of $K$, with associated
  discrete valuation ring $R$ having residue field $k$ such that (a) $E_k$
  is smooth, (b) $\tau_k^d\ne 1$, and (c) the nonisomorphic stable
  constituents of $V$ have nonisomorphic reductions.  (Such a place exists
  because we are imposing finitely many open conditions.)  Then
  $B_{V_R,\Psi_R}$ is a flat extension of $B_{V,\Psi}$, so that the
  dimension of the space of derivations of degree $-d$ is upper
  semicontinuous, implying that $B_k$ itself has a nonzero derivation of
  degree $-d$.  Since $\dim(k)=\dim(K)-1$, the claim follows by induction
  on the dimension.
\end{proof}

\begin{rem}
  A similar reduction shows that any elliptic algebra is Noetherian, by a
  direct generalization of the proof of
  \cite[Thm. 8.3]{ArtinM/TateJ/VandenBerghM:1990}.
\end{rem}

This is enough to establish the following for an arbitrary elliptic algebra
in characteristic 0.

\begin{thm}
  Let $B$ be an elliptic algebra over a field of characteristic 0, twisted
  by an arbitrary translation.  Then the moduli space of filtered
  deformations of $B$ is a closed substack of a weighted projective space.
\end{thm}

\begin{proof}
  If $\tau$ has finite order, then the arguments of Lemma
  \ref{lem:all_derivs_inner} and Corollary \ref{cor:no_neg_derivs} carry
  over mutatis mutandis to show that $B$ has no negative derivations.
  Otherwise, $\ord(\tau)=\infty$, and we may simply replace Lemma
  \ref{lem:all_derivs_inner} by Lemma \ref{lem:no_neg_derivs_general}.
\end{proof}

The argument of Proposition \ref{prop:center_lifts} also carries over.

\begin{thm}
Let $B$ be an elliptic algebra of characteristic 0 twisted by translation
by a torsion point.  The center of any filtered deformation of $B$ is a
filtered deformation of $Z(B)$.
\end{thm}

\section{Resolutions of elliptic algebras}

Let $E/k$ be an elliptic curve over an algebraically closed field
$k$, let $V$ be a vector bundle on $E$, and let $\Psi:\coh(E)\to
\coh(E)$ be an autoequivalence of the form $M\mapsto {\cal L}\otimes
(\tau^{-1})^* M$ with $\tau\in E\subset \Aut(E)$ and ${\cal L}$ an
ample line bundle of degree $d>0$.  (In fact, most of the arguments
below work even if $\tau$ is not a translation, since they only
depend on how $\Psi$ affects slopes; the only exception is in the
explicit calculations of minimal resolutions, where minimality also
depends on how it affects determinants.) As mentioned above, these
ingredients determine a graded algebra $B_{V,\Psi}$ such that for
$n\in \Z$,
\[
(B_{V,\Psi})_n = \Hom(V,\Psi^n V),
\]
with the obvious composition.  Since ${\cal L}$ is ample, this
algebra is bounded below, and we are interested in its filtered
deformations.

Any summand of $V$ determines an idempotent of $B_{V,\Psi}$, which lifts
(nonuniquely) to any filtered deformation by Proposition
\ref{prop:ids_lift}.  It follows immediately that applying a power of
$\Psi$ to a summand of $V$ does not change the space of filtered
deformations---it simply shifts the degrees of the corresponding
off-diagonal blocks.  Since any vector bundle on an elliptic curve is a sum
of semistable bundles (see \cite{AtiyahMF:1957}, which also completely
determines the structure of indecomposable bundles) and $\Psi$ adds $d$ to
the slope of a bundle, this freedom allows us to assume that every summand
of $V$ has slope in $[0,d)$ (or, more generally, in any desired half-open
  interval of length $d$).  Since a nonzero morphism between semistable
  bundles forces the slope of the codomain to be at least the slope of the
  domain, this constraint on $V$ implies that $(B_{V,\Psi})_n=0$ for $n<0$.

This, of course, is not quite canonical; if $V$ has $c$ different slopes
(modulo $d$), then there are $c$ inequivalent choices of nonnegative
gradings on $B_{V,\Psi}$.  The most natural way to resolve this involves
replacing $B_{V,\Psi}$ by a category.  For each slope $\mu\in \Q$, let
$V_{[\mu]}$ denote the direct sum of all slope $\mu$ indecomposable
summands of the modules $\Psi^n V$, so that the construction of the
previous paragraph involves replacing $V$ by $\bigoplus_{\mu\in
  [\alpha,\alpha+d)} V_{[\mu]}$.  The slopes for which $V_{[\mu]}\ne 0$
  form a discrete subset of $\Q$, a union of finitely many cosets of $d\Z$,
  and thus we may choose an ordered bijection between $\Z$ and the set of
  such slopes.  If $i\mapsto \mu_i$ is that bijection, then we may define
  $V_i:=V_{[\mu_i]}$.  We then let ${\cal B}_{V,\Psi}$ be the subcategory
  of $\coh E$ determined by $V_i$; i.e., it is the category with objects
  indexed by integers and $\Hom$ spaces given by
\[
{\cal B}_{V,\Psi}(i,j) = \Hom(V_i,V_j).
\]
(This is a nonnegatively graded $\Z$-algebra in the sense of
\cite{BondalAI/PolishchukAE:1993}.)  It is then straightforward to see that
the category of modules over this category (i.e., contravariant functors to
the category of vector spaces) is isomorphic to the category of graded
$B_{V,\Psi}$-modules.  Note that since $\Psi$ induces an autoequivalence
shifting the objects by $c$, this category induces a graded algebra which
in degree $d$ is
\[
\bigoplus_{1\le i,j\le c} {\cal B}_{V,\Psi}(i,j+d).
\]
This graded algebra agrees with the na\"{i}ve grading on $B_{V,\Psi}$, up
to rescaling the grading by $c$ and shifting the grading on off-diagonal
blocks relative to the natural idempotents, and thus gives the same space
of filtered deformations.  (The only caveat is that the additional terms
added to the relations must remain in the same congruence class modulo $c$,
but this is actually forced by consistency with the degree 0 involutions.)
With this in mind, where $B_{V,\Psi}$ appears below, we will always take
this canonical grading.

The theory of minimal resolutions, though most familiar in the case of
connected graded algebras, applies more generally (following
\cite{EilenbergS:1956}) to nonnegatively graded algebras which are
finite-dimensional in every degree.  In particular, any $B_{V,\Psi}$-module
has such a minimal resolution, and thus the same applies to ${\cal
  B}_{V,\Psi}$.  The one caveat here is that when $B_0$ is not semisimple,
the minimal resolution of a finitely graded module need not be convergent;
that is, there may be degrees such that the restriction of the complex to
those degrees gives an infinite resolution.

\begin{lem}
  Let $B$ be a nonnegatively graded $k$-algebra such that $B_0$ is
  finite-dimensional.  Then the minimal resolution of $B_0$ as a graded
  $B$-module is convergent iff each $B_n$ has finite homological dimension
  as a $B_0$-module.  Moreover, if this holds, then for any bounded below
  graded $B$-module $M$, $M$ has convergent minimal resolution iff each
  $M_n$ has finite homological dimension.
\end{lem}

\begin{proof}
  First suppose that $B_n$ and $M_n$ all have finite homological dimension,
  and WLOG assume that $M_n=0$ for $n<0$ and $M_0\ne 0$.  Then the
  restriction of the minimal resolution of $M$ to degree 0 gives a
  (minimal!) projective resolution of $M_0$ as a $B_0$-module, which thus
  has finite length.  In particular, projective modules of degree 0 appear
  in only finitely many terms, and thus we may consider the partial
  resolution truncated at the last point where they occur.  In each degree,
  the terms of the resulting complex have finite homological dimension
  (being homogeneous components of $B$ or $M$) and thus, since the complex
  is exact except at the very left, the kernel also has finite homological
  dimension.  This gives us a new graded module $M'$ satisfying the same
  hypotheses but with a larger lower bound on its degrees, and thus the
  result follows by induction.  This gives the ``if'' part of the second
  claim, and thus (taking $M=B_0)$ of the first claim as well.

  Now, supposing still that the components $B_n$ all have finite
  homological dimension, let $M$ be a bounded below module with at least
  one homogeneous component of infinite homological dimension.  Then the
  same argument of the previous paragraph allows us to reduce to the case
  that $M_n=0$ for $n<0$ but $M_0$ has infinite homological dimension.
  (Here we use the fact that in a partial resolution by modules of finite
  homological dimension, the kernel has finite dimension iff the cokernel
  does.)  But then the minimal resolution in degree 0 is a projective
  $B_0$-module resolution of $M_0$, which is infinite by assumption,
  preventing convergence.

  Finally, if $B$ has a homogeneous component of infinite homological
  dimension, let $n$ be the minimal such degree.  Then the degree $n$ part
  of the minimal resolution of $B_0$ is a resolution of $B_n$ by sums of
  summands of components $B_i$ for $i<n$, and thus of terms having finite
  homological dimension.  Since $B_n$ has infinite homological dimension,
  this forces the resolution to have infinite length.
\end{proof}

This property holds for the canonical grading on $B_{V,\Psi}$, and in fact
we have something stronger.

\begin{lem}
  Suppose that $V$, $W$ are indecomposable bundles with $\mu(V)<\mu(W)$.
  Then $\Hom(V,W)$ is free as a
  $\End(V)\otimes_k\End(W)^{\text{op}}$-module.
\end{lem}

\begin{proof}
  We have $\End(V)\cong k[x]/(x^n)$ and $\End(W)\cong k[y]/(y^m)$ for
  suitable positive integers $m$ and $n$, and thus we need to show that
  $\Hom(V,W)$ is a free $k[x,y]/(x^n,y^m)$-module, or equivalently that
  \[
  \dim_k\Hom(V,W) = mn \dim_k(\Hom_{\End(V)}(k,\Hom(V,W)\otimes_{\End(W)}k)),
  \]
  or in other words
  \[
  \dim_k\Hom(V,W) = mn \dim_k(\Hom(V\otimes_{\End(V)}k,W\otimes_{\End(W)}k)).
  \]
  Since $V\otimes_{\End(V)}k$ is the stable constituent of $V$, and
  similarly for $W\otimes_{\End(W)}k$, this certainly applies at the level
  of Euler characteristics, and the difference in slopes forces $\Ext^1$ to
  vanish (by duality).
\end{proof}

\begin{lem}
  Suppose that $V$, $W$ are semistable bundles with $\mu(V)<\mu(W)$.  Then
  $\Hom(V,W)$ is projective as a
  $\End(V)\otimes_k\End(W)^{\text{op}}$-module.
\end{lem}

\begin{proof}
  Since each of $V$ and $W$ splits naturally as a direct sum over
  nonisomorphic stable constituents, we may assume that each is built from
  a single stable bundle.  In particular, $V$ has a unique isomorphism
  class $V_{\max}$ of indecomposable summands of maximal rank, and one has
  \[
  V\cong V_{\max}\otimes_{\End(V_{\max})}M_V
  \]
  for some faithful $\End(V_{\max})$-module $M_V$, with a similar
  description applying to $W$.  We then find that
  \[
  \Hom(V,W)\cong
  \Hom_{\End(V_{\max})}(M_V,\Hom(V_{\max},W_{\max})\otimes_{\End(W_{\max})}M_W).
  \]
  Since $\Hom(V_{\max},W_{\max})$ is free, it remains only to show that
  \[
  \Hom_{\End(V_{\max})}(M_V,\End(V_{\max}))\otimes_k M_W
  \]
  is projective as an $\End(V)\otimes_k \End(W)^{\text{op}}$-module.

  Here $M_W$ is a faithful module over $\End(W_{\max})\cong k[x]/(x^n)$,
  and (since $x$ is central in $\End(W)$) we have $\End(W)\cong
  \End_{k[x]/(x^n)}(M_W)$.  In particular, $M_W\cong \bigoplus_i
  k[x]/(x^{m_i})$ with $n=m_1\ge m_2\ge\cdots$, and thus there is an
  idempotent cutting out a copy of $k[x]/(x^n)$, so that the corresponding
  submodule of $\End_{k[x]/(x^n)}(M_W)$ is $M_W$, implying projectivity.
  The claim for $V$ is analogous (relative now to the left $\End(V)$-module
  structure on $\End(V)$).
\end{proof}

\begin{prop}
  The homogeneous components of $B_{V,\Psi}$ are projective modules over
  $(B_{V,\Psi})_0$.
\end{prop}

\begin{proof}
  This reduces to showing that the Hom spaces ${\cal B}_{V,\Psi}(i,j)$ are
  projective over the relevant endomorphism ring.  If $j=i$, this is immediate,
  while if $j>i$, we have in fact shown that it is projective over the
  tensor product of endomorphism rings.
\end{proof}

We thus see that minimal resolutions of $B_{V,\Psi}$-modules (or,
equivalently, ${\cal B}_{V,\Psi}$-modules) are convergent as long as they
are bounded below and have homogeneous components of finite homological
dimension.  Let $P_i$ denote the projective ${\cal B}_{V,\Psi}$-module
given by the contravariant functor
\[
j\mapsto {\cal B}_{V,\Psi}(j,i),
\]
and let $S_i$ be the module agreeing with $P_i$ in degree $i$ but $0$ in
all other degrees.  Since $S_i(i)=\End(V_i)$, $S_i$ has a convergent
projective resolution.

An important property of the $\Z$-algebras ${\cal B}_{V,\Psi}$ is that they
are Gorenstein in a suitable sense.

\begin{prop}
  The space
  \[
  \Ext^p(S_i,P_j)
  \]
  vanishes unless $i=j$ and $p=2$, in which case there is a canonical
  isomorphism
  \[
  \Ext^2(S_i,P_i)\cong\Ext^1(V_i,V_i).
  \]
\end{prop}

\begin{proof}
  The module $S_i$ has a convergent projective resolution
  \[
  \cdots\to Q_{i+2}\to Q_{i+1}\to P_i\to S_i.
  \]
  For each $j>i$, projectivity of $Q_j$ implies that it is the module of
  global sections of a vector bundle $W_j$ on $E$ (this is true for any
  summand of any $P_j$) and thus the complex
  \[
  \cdots\to W_{i+2}\to W_{i+1}\to V_i
  \]
  of sheaves on $E$ is exact.  Now, $R\Hom(S_i,P_j)$ is represented by the
  complex
  \[
  \Hom(P_i,P_j)\to \Hom(Q_{i+1},P_j)\to \Hom(Q_{i+2},P_j)\to\cdots,
  \]
  or equivalently
  \[
  \Hom(V_i,V_j)\to \Hom(W_{i+1},V_j)\to \Hom(W_{i+2},V_j)\to\cdots.
  \]
  Since the double complex
  \[
  R\Hom(V_i,V_j)\to R\Hom(W_{i+1},V_j)\to R\Hom(W_{i+2},V_j)\to\cdots
  \]
  represents $R\Hom(0,V_j)$, it is acyclic, which, since every column has
  magnitude contained in $[0,1]$ (as a $R\Hom$ between sheaves on $E$),
  implies a quasi-isomorphism between the complex of $\Hom$ spaces and a
  shift of the complex
  \[
  \Ext^1(V_i,V_j)\to \Ext^1(W_{i+1},V_j)\to \Ext^1(W_{i+2},V_j)\to\cdots ,
  \]
  or (after fixing an isomorphism $H^1(\sO_E)\cong k$)
  \[
  \Hom(V_j,V_i)^*\to \Hom(V_j,W_{i+1})^*\to \Hom(V_j,W_{i+2})^*\to\cdots .
  \]  
  This is dual to the complex
  \[
  \cdots \to\Hom(V_j,W_{i+2})\to \Hom(V_j,W_{i+1})\to \Hom(V_j,V_i),
  \]
  or equivalently
  \[
  \cdots \to\Hom(P_j,Q_{i+2})\to \Hom(P_j,Q_{i+1})\to \Hom(P_j,P_i),
  \]  
  which (up to a shift) represents $R\Hom(P_j,S_i)\cong \Hom(P_j,S_i)$.
  This vanishes unless $i=j$, when it equals $\End(V_i)$.

  We thus conclude that $\Ext^p(S_i,P_j)=0$ unless $i=j$ and (taking into
  account the shift) $p=2$, in which case one has a canonical isomorphism
  \[
  \Ext^2(S_i,P_j)\cong \Hom_k(\End(V_i),H^1(\sO_E))\cong \Ext^1(V_i,V_i)
  \]
  as required.
\end{proof}

Now, let ${\cal B}^+$ be the $\Z$-algebra obtained from ${\cal B}_{V,\Psi}$
by adjoining a central element $t$ of degree $c$.  In other words,
\[
{\cal B}^+(i,j) = \bigoplus_{l\ge 0} {\cal B}_{V,\Psi}(i,j-cl),
\]
with $t\in {\cal B}^+(i,i+c)$ the element corresponding to $1\in {\cal
  B}_{V,\Psi}(i,i)$.  Let $P^+_i$ be the natural family of projective
${\cal B}^+$-modules, with associated modules $S^+_i$, and note that we
still have ${\cal B}^+(i,i)=\Hom(V_i,V_i)$.

\begin{cor}
  The space
  \[
  \Ext^p(S^+_i,P^+_j)
  \]
  vanishes unless $j=i-c$ and $p=3$, in which case there is a canonical
  isomorphism
  \[
  \Ext^3(S^+_i,P^+_{i-c})\cong\Ext^1(V_i,V_i).
  \]
\end{cor}

\begin{proof}
  Since $P^+_j$ restricts to $\bigoplus_{l\ge 0} P_{j-lc}$, we have
  \[
  \Ext^p(S_i[t],P^+_j)\cong \bigoplus_{l\ge 0} \Ext^p(S_i,P_{j-lc}),
  \]
  where $S_i[t]$ is the induced module.  We may thus use the short
  exact sequence
  \[
  0\to S_i[t][-c]\to S_i[t]\to S^+_i\to 0
  \]
  to obtain a long exact sequence
  \[
  \cdots\to
  \Ext^p(S^+_i,P^+_j)\to 
  \bigoplus_{l\ge 0} \Ext^p(S_i,P_{j-lc})\to
  \bigoplus_{l\ge -1} \Ext^p(S_i,P_{j-lc}) \to\cdots
  \]
  This immediately gives an isomorphism
  \[
  \Ext^p(S^+_i,P^+_j)\cong \Ext^{p-1}(S_i,P_{j+c}),
  \]
  so the formula reduces to the analogous formula for ${\cal B}$.
\end{proof}

Given a graded algebra (or $\Z$-algebra) $A$, there is an associated
category $\Proj(A)$, the quotient of the category of all graded modules by
the subcategory of ``torsion'' modules (in which every element generates a
finite-dimensional submodule).  The $\Ext$ groups in that category can then
be computed as limits of $\Ext$ groups in the category of graded modules,
per \cite{ArtinM/ZhangJJ:1994}, subject to a certain technical condition
(``property $\chi$'') which is satisfied for $B_{V,\Psi}$ by Thm 4.5
op.~cit. and thus for the Rees algebra of any filtered deformation by
Thm. 8.8 op.~cit.  Indeed, in our setting, the only nontrivial hypothesis
between the two theorems is that $\Psi$ is ``ample'' relative to $V$: for
any coherent $M$, there exists $n$ such that $\Ext^p(\Psi^{-n}V,M)=0$ for
$p>0$ and $\Hom(\Psi^{-n}V,M)\otimes_k \Psi^{-n}V\to M$ is surjective.
Both facts reduce immediately to the case that both $V$ and $M$ are
indecomposable, with vanishing of $\Ext^1$ an immediate consequence via
duality of the eventual inequality in slopes, and global generation proved
in Corollary \ref{cor:slope_global_gen} below.

\begin{thm}
  Let $\bar{P}^+_i$ denote the object of $\Proj({\cal B}^+)$ corresponding
  to $P^+_i$.  Then $\Ext^p(\bar{P}^+_i,\bar{P}^+_j)=0$ unless $p=0$ and
  $j\ge i$ or $p=2$ and $j\le i-c$, with
  \[
  \Hom(\bar{P}^+_i,\bar{P}^+_j) \cong \Hom(V_i,V_j)
  \]
  and
  \[
  \Ext^2(\bar{P}^+_i,\bar{P}^+_j) \cong \Ext^1(V_i,V_{j+c}).
  \]
\end{thm}

\begin{proof}
  We have
  \[
  \Ext^p(\bar{P}^+_i,\bar{P}^+_j)\cong \lim_{n\to\infty}
   \Ext^p((P^+_i)_{\ge n},P^+_j).
  \]
  Now consider the short exact sequence
  \[
  0\to (P^+_i)_{\ge n}\to P^+_i\to P^+_i/(P^+_i)_{\ge n}\to 0.
  \]
  Since $P^+_i$ is projective, $R\Hom(P^+_i,P^+_j)=\Hom(V_i,V_j)$
  concentrated in degree 0.  Moreover, since the degree $l$ component of
  $P^+_i$ is a free $\End(V_l)$-module, so $P^+_i/(P^+_i)_{\ge n}$ is a
  finite extension of such modules, we find that for $n>j-c$,
  \[
  R\Hom(P^+_i/(P^+_i)_{\ge n},P^+_j)
  \cong
  \Ext^1(V_i,V_{j+c})
  \]
  concentrated in degree 3.  The desired claim follows immediately.
\end{proof}

For the next claim, we revert to one of the na\"ive nonnegative gradings of
$B_{V,\Psi}$.

\begin{cor}\label{cor:fildef_is_rigid_unobstructed}
  Let $A^+$ be the Rees algebra of a filtered deformation of $B_{V,\Psi}$,
  and let $\sO(l)$ be the objects of $\Proj(A^+)$ corresponding to shifts
  of $A^+$.  Then for $l\le m$, $\Ext^p(\sO(l),\sO(m))=0$ for $p\ne 0$ and
  for $l>m$, $\Ext^p(\sO(l),\sO(m))=0$ for $p\ne 2$.  Moreover, there is a
  canonical duality
  \[
  \Ext^2(\sO(l),\sO(m))\cong \Hom_k(\Hom(\sO(m),\sO(l-1)),H^1(\sO_E)).
  \]
  In particular, $\sO(0)$ has no higher endomoprhisms, and if $V$ is
  stable, then $\sO(0)$ is exceptional.
\end{cor}

\begin{proof}
  The $\G_m$ action on filtered deformations gives rise to a family of
  graded algebras with generic fiber isomorphic to $A^+$ and special fiber
  associated to the trivial deformation $B^+\cong B[t]$.  By
  semicontinuity, it suffices to prove vanishing when $A^+$ is trivial,
  where it reduces to computing
  \[
  R\Hom(\bigoplus_{0\le
    i<c} \bar{P}^+_{i-lc},
  \bigoplus_{0\le
    i<c} \bar{P}^+_{i-mc})
  \]
  in $\Proj({\cal B}^+)$.  Note that since each case has precisely one
  nonvanishing cohomology space, the corresponding sheaves on
  $\A^1=\G_m\cup \{0\}$ are flat.

  For the duality between nonvanishing $\Ext$ spaces, we note that
  $\coh(E)\cong \Proj(B_{V,\Psi})$ embeds in $\Proj(A^+)$ as the
  subcategory on which the natural transformation ${-}(-1)\to \text{id}$
  associated to the central element of $A^+$ vanishes.  This embedding
  factors into a pair $i_*i^*$ of adjoint functors with $i^*:\Proj(A^+)\to
  \coh(E)$ exact, so that we have distinguished triangles
  \[
  M(-1)\to M\to Ri_*i^*M\to 
  \]
  for any $M$ and thus distinguished triangles
  \[
  R\Hom(M,N(-1))\to R\Hom(M,N)\to R\Hom_E(i^*M,i^*N)\to.
  \]
  In particular, taking $M=N=\sO(l)$ gives a canonical isomorphism
  \[
  \Ext^2(\sO(l),\sO(l-1))\cong \Ext^1(\Psi^lV,\Psi^lV)
  \]
  and thus a canonical trace
  \[
  \Ext^2(\sO(l),\sO(l-1))\to H^1(\sO_E)(\cong k).
  \]
  This induces pairings in the usual way, and it remains only to show that
  the pairings are perfect.  But this follows from another semicontinuity
  argument together with the fact they they are perfect in $\Proj(B^+)$.
\end{proof}

\begin{rem}
  Note that the trace form produces two pairings:
  \[
  \Hom(\sO(l),\sO(m))\otimes \Ext^2(\sO(m),\sO(l-1))\to
  \Ext^2(\sO(l),\sO(l-1))\to k
  \]
  and
  \begin{align}
    &\Ext^2(\sO(m),\sO(l-1))\otimes \Hom(\sO(l),\sO(m))\notag\\
  {}\cong {}&\Ext^2(\sO(m),\sO(l-1))\otimes \Hom(\sO(l-1),\sO(m-1))
  \to
  \Ext^2(\sO(m),\sO(m-1))
  \to k.
  \end{align}
  The corresponding pairings in $\Proj(B^+)$ agree (since the analogous
  pairings agree in $\Proj(B_{V,\Psi})\cong \coh(E)$), and thus they differ
  by a compatible family of automorphisms of $\Hom(\sO(l),\sO(m))$ agreeing
  with 1 on the associated graded, or in other words by an automorphism of
  $A^+$ acting as $1$ on the associated graded.  But this implies that the
  pairings agree: if this automorphism were nontrivial, it would induce a
  nonzero derivation of $B_{V,\Psi}$ of negative degree.  It follows that
  (up to a choice of holomorphic differential on $E$) ${-}(-1)[2]$ is a
  Serre functor on $\Proj(A^+)$, since the relevant duality and
  compatibility relations hold on generators.  (In particular,
  $S[-2]={-}(-1)$ is an abelian autoequivalence, so that $\Proj(A^+)$ is
  Gorenstein.)
\end{rem}
  
\medskip

If the components of the slope decomposition of $V$ have commutative
endomorphism rings, we can be more precise about the structure of the
minimal resolution.  To understand this, we first recall some facts about
derived autoequivalences of elliptic curves.  In addition to the abelian
autoequivalences, which consist of automorphisms of $E$ along with
translations and twists by line bundles, there is an additional
autoequivalence $\Phi_{\sO_E}$ which fits into a distinguished triangle
\[
\sO_E\otimes_k R\Hom(\sO_E,M)\to M\to \Phi_{\sO_E} M\to.
\]
(We may either use Fourier-Mukai theory or a suitable dg-enhancement to see
that this cone can be made functorial.)  The inverse is also a cone:
\[
\Phi_{\sO_E}^{-1}M\to M\to \Hom_k(\Hom(M,\sO_E),\sO_E)\to,
\]
which can be seen by observing that
\[
R\Hom(\sO_E,\Phi_{\sO_E}M)\cong \Ext^1(\sO_E,\sO_E)\otimes R\Hom(\sO_E,M)
\cong
R\Hom(M,\sO_E)^*
\]
Twisting by the line bundle $\sO_E(z)$ (where $z$ is the identity) has a
similar description:
\[
\sO_z\otimes_k R\Hom(\sO_z,M)\to M\to M(z)\to
\]
so that we may call it $\Phi_{\sO_z}$.  (These functors generate a central
extension of $\SL_2(\Z)$.)  The following is then standard (essentially
Atiyah's classification of vector bundles \cite{AtiyahMF:1957}):

\begin{prop}\label{prop:atiyah}
  Any semistable sheaf on $E$ is the image of a torsion sheaf under a
  composition of the functors $\Phi_{\sO_E}$ and $\Phi_{\sO_z}$ and their
  inverses.
\end{prop}

\begin{rem}
  If the semistable sheaf $M$ is not already torsion, then we twist by a
  suitable power of $\sO_E(z)$ to put its slope in $(0,1]$ and observe that
applying $\Phi_{\sO_E}$ gives a sheaf of smaller rank, so that iterating
eventually produces a torsion sheaf.
\end{rem}

\begin{rem}
  Since $\Phi_{\sO_E}$ and $\Phi_{\sO_z}$ both preserve the degree 0
  divisor class $c_1(M) - \deg(M)[z]$ and a stable sheaf is uniquely
  determined by its determinant, the torsion sheaf that results is uniquely
  determined, even though the derived equivalence getting there is not.
\end{rem}

\begin{rem}
  Note that since the action of $\Psi_{\sO_E}$ and $\Phi_{\sO_z}$ preserves
  $\gcd(\rank(M),\deg(M))$, we conclude that $M$ has
  $\gcd(\rank(M),\deg(M))$ stable constituents (counted with multiplicity),
  by reduction to the torsion case.
\end{rem}

\begin{lem}
Suppose that $M\in \coh(E)$ is such that $\End(M)$ is commutative.  Then
$M$ is semistable, and each indecomposable summand of $M$ is S-equivalent
to a power of a different isomorphism class of stable sheaves.
\end{lem}

\begin{proof}
  Let $M=\bigoplus M_i$ be the decomposition of $M$ into indecomposable
  summands.  We clearly have $\Hom(M_i,M_j)=0$ for $i\ne j$, and thus
  $\mu(M_i)=\mu(M_j)$.  Since an indecomposable sheaf is semistable, it
  follows immediately that $M$ is semistable.  Moreover, each
  indecomposable sheaf is an iterated extension of a unique stable
  constituent, and the vanishing of $\Hom(M_i,M_j)$ forces the stable
  constituents to be nonisomorphic.
\end{proof}

\begin{cor}
  A sheaf $M$ on $E$ has commutative endomorphism ring iff it is the image
  of the structure sheaf of an effective divisor $D_M$ under an element of
  the group of derived autoequivalences generated by $\Phi_{\sO_E}$ and
  $\Phi_{\sO_z}$.
\end{cor}

\begin{proof}
  Since the structure sheaf of an effective divisor is commutative, the
  image of such a sheaf under a derived autoequivalence also has
  commutative endomorphism ring.  Conversely, if $M$ has commutative
  endomorphism ring, then the torsion sheaf it is derived equivalent to
  must also have commutative endomorphism ring, so is the structure sheaf
  of a $0$-dimensional subscheme as required.
\end{proof}

\begin{rem}
  With this in mind, we call sheaves on $E$ with commutative endomorphism
  ring ``divisorial''.
\end{rem}

Now, given a divisorial sheaf $M$, consider the functor $\Phi_M$ on
$D^b\coh(E)$ fitting into the distinguished triangle
\[
M\otimes^{\bf L}_{\End(M)}R\Hom(M,N)\to N\to \Phi_M(N)\to.
\]
(Again, we either use a Fourier-Mukai kernel or a dg-enhancement to ensure
that this functor is well-defined.)

\begin{prop}\label{prop:phiM_preserves_shifted_sheaves}
  If $M$ is divisorial, then the functor $\Phi_M$ is a derived
  autoequivalence, and if $N$ is a semistable sheaf then either $\Phi_M(N)$
  or $\Phi_M(N)[-1]$ is a semistable sheaf.
\end{prop}

\begin{proof}
  The conjugate of the functor $\Phi_M$ under a derived autoequivalence has
  the same form, and thus to see that it is a derived autoequivalence, it
  suffices to consider the case that $M$ is the structure sheaf of a
  divisor.  But in that case, $\Phi_{\sO_D}N\cong N(D)$, which is clearly
  an autoequivalence.  It also follows from this that the functor
  $M\otimes_{\End(M)}{-}$ is exact, and thus we may take the underived
  tensor product above.

  For any sheaf $N$, we then have an exact sequence
  \[
  0\to h^{-1}\Phi_M(N)\to M\otimes_{\End(M)}\Hom(M,N)\to N\to
  h^0\Phi_M(N)\to M\otimes_{\End(M)}\Ext^1(M,N)\to 0
  \]
  with $h^p\Phi_M(N)=0$ for $p\notin\{-1,0\}$.  If $N$ is stable, then it
  is simple ($\End(N)\cong k$), and the same applies to $\Phi_M(N)$,
  forcing $\Phi_M(N)$ to be a shift of a simple (and thus stable!) sheaf.
  (Since $\Ext^2$ vanishes in $\coh(E)$, any complex in $D^b_{\coh}(E)$ is
  formal, and thus its endomorphism ring contains at least one idempotent
  for each nonvanishing cohomology sheaf.)  Since there are only two
  possible shifts, we can distinguish whether it is a sheaf or a shift from
  its class in the {\em numerical} Grothendieck group, so that the claim
  immediately extends to semistable sheaves.
\end{proof}

\begin{cor}
  If the divisorial sheaf $M$ has stable constituents $M_1$,\dots,$M_m$
  (with multiplicity, and in any order), then $\Phi_M\cong
  \Phi_{M_1}\circ\cdots\circ \Phi_{M_m}$.  Moreover, $\Phi_M$ depends only
  on $\rank(M)$ and $\det(M)$.
\end{cor}

\begin{proof}
  When $M$ is torsion, the claims are immediate; here the second claim
  simply says that the autoequivalence $N\mapsto N(D)$ depends only on the
  isomorphism class of the line bundle $\sO_E(D)$.  The claim in general
  follows by conjugating by the appropriate element of the group generated
  by $\Phi_{\sO_E}$ and $\Phi_{\sO_z}$.
\end{proof}

\begin{rem}
  Note that the claim as stated descends to smooth genus 1 curves over
  non-closed fields, where the functor $\Phi_{\sO_z}$ is not defined.
\end{rem}

The second claim can be restated as saying that $\Phi_M$ only depends on
the class of $M$ in $K_0(E)$.  If $M$ is stable, then we clearly have
\[
[\Phi_M(N)] = [N]-\chi(M,N)[M],
\]
where $\chi(M,N):=\sum_{i\in \Z} (-1)^i\dim\Ext^i(M,N)$ is the Mukai
pairing on $K_0$.  More generally, if $M$ has $m=\gcd(\rank(M),\deg(M))$
stable constituents, then we can use the factorization along with the fact
that $\chi(M_i,M_j)=0$ (since they are stable sheaves of the same slope) to
find that
\[
[\Phi_M(N)] = [N]-m^{-1}\chi(M,N)[M].
\]
Since
\[
\chi(M,N) = \rank(M)\deg(N)-\deg(M)\rank(N),
\]
this lets us explicitly compute the action of $\Phi_M$ on the Grothendieck
group.

\begin{cor}
  Suppose that the divisorial sheaf $M$ has rank $r$, degree $d$, and
  $m=\gcd(r,d)$ stable constituents.  Then the autoequivalence $\Phi_M$
  acts on the Grothendieck group of $E$ by
  \begin{align}
  \rank(\Phi_M(N)) &= \frac{m+dr}{m}\rank(N) - \frac{r^2}{m}\deg(N)\notag\\
  \deg(\Phi_M(N)) &= \frac{d^2}{m}\rank(N) + \frac{m-dr}{m}\deg(N)
  \end{align}
  and
  \[
  \det(\Phi_M(N)) \cong \det(N)\otimes \det(M)^{-(r\deg(N)-d\rank(N))/m}.
  \]
\end{cor}

\begin{rem}
  As a function of $m$ with $r/m$, $d/m$ fixed, this gives a homomorphism
  from $\Z$ to $\SL_2(\Z)$; in particular, the inverse is just the image
  under negating $m$, $r$, $d$.
\end{rem}

In the following statement, and below, the convention is that the slope of
a (nonzero) torsion sheaf is {\em positive} infinity.  (This is compatible
with the definition of semistability: a semistable sheaf of positive rank
has no torsion subsheaf but plenty of torsion quotients.)

\begin{cor}\label{cor:slope_global_gen}
  Let $V$ be a divisorial vector bundle with $\rank(V)=r$, $\deg(V)=d$.
  Then for any semistable sheaf $M$, the morphism
  \[
  V\otimes_{\End(V)} \Hom(V,M)\to M
  \]
  is surjective if
  \[
  \mu(M)>\mu(V)+\frac{\gcd(r,d)}{r^2},
  \]
  and is otherwise injective.  In either case, the (co)kernel is
  semistable, of slope $\le \mu(V)-\frac{\gcd(r,d)}{r^2}$ iff the morphism is
  surjective.
\end{cor}

\begin{proof}
  The strict lower bound on the slope of $M$ ensures that
  $\rank(\Phi_V(M))<0$, and thus $\Phi_V(M)[-1]$ is a semistable sheaf
  fitting into a short exact sequence
  \[
  0\to \Phi_V(M)[-1]\to V\otimes_{\End(V)} \Hom(V,M)\to M\to 0.
  \]
  If the opposite inequality holds, then $\rank(\Phi_V(M))>0$, so that
  $\Phi_V(M)$ is a semistable sheaf and the morphism is injective.
  Finally, if equality holds, then $\rank(M)>0$ and eliminating $\deg(M)$
  gives (assuming $V$ has $m$ stable constituents)
  \[
  \deg(\Phi_V(M))=\frac{m}{r^2} \rank(M),
  \]
  which since $\rank(\Phi_V(M))=0$ again forces $\Phi_V(M)$ to be a
  semistable sheaf.

  If $\mu(M)>\mu(V)$, the map vanishes, so the cokernel is certainly
  semistable.  When $\mu(M)<\mu(V)$, the (co)kernel is a shift of
  $\Phi_V(M)$, so is semistable by Proposition
  \ref{prop:phiM_preserves_shifted_sheaves}; when equality holds (forcing
  injectivity), this fails, but the morphism is between semistable bundles
  of the same slope, and thus the cokernel must also have that slope (which
  also agrees with the slope of $\Phi_V(M)$).

  Since $\mu(\Phi_V(M))$ is a linear fractional function of $\mu(M)$, the
  surjective region in terms of $\mu(\Phi_V(M))$ is the image of the
  interval $(\mu(V)+\frac{\gcd(r,d)}{r^2},\infty]$ under that LFT, giving
  the stated result.
\end{proof}

\begin{prop}
  Let $V$ be a divisorial vector bundle on $E$, and
  suppose $M$ is a coherent sheaf such that $(\Phi_V\Psi)^n M[-n]$ is a
  sheaf for all $n$.  Then $M$ has a natural resolution of the form
  \[
  \cdots\to \Psi^{-n}V\otimes_{\End(V)}M_n\to
   \cdots\to \Psi^{-1}V\otimes_{\End(V)}M_1\to M \to 0
  \]
  for suitable $\End(V)$-modules $M_i$, and the corresponding complex of
  $B_{V,\Psi}$-modules is exact in all positive $B$-degrees.
\end{prop}

\begin{proof}
  For the first claim, we first refine to say that there is an exact
  sequence
  \[
  0\to \Psi^{-n}(\Phi_V\Psi)^n M[-n]
   \to \Psi^{-n}V\otimes_{\End(V)}M_n\to
   \cdots\to \Psi^{-1}V\otimes_{\End(V)}M_1\to M \to 0.
  \]
  for all $n$.  This is straightforward: it is trivial for $n=0$, and we
  may increase $n$ using the short exact sequence
  \begin{align}
  0\to \Psi^{-n-1}(\Phi_V\Psi)^{n+1}M[-n-1]
   &\to \Psi^{-n-1}V\otimes_{\End(V)} \Hom(V,\Psi(\Phi_V\Psi)^n M[-n])\notag\\
   &\to \Psi^{-n}(\Phi_V\Psi)^n M[-n] \to 0
  \end{align}
  obtained by applying $\Psi^{-n-1}$ to the short exact sequence associated
  to $\Phi_V(\Psi(\Phi_V\Psi)^nM[-n])$.

  To see that this also gives a resolution of the (truncated!) $B$-module
  associated to $M$, observe that every term (except possibly $M$ itself)
  of the sequence for $n$ is acyclic for $R\Hom(V,\Psi^{n+1}{-})$.  For the
  intermediate terms, this is immediate from the fact that
  $\mu(\Psi^iV)>\mu(V)$, while for the first term it follows from the fact
  that $(\Phi_V\Psi)^{n+1}M[-n-1]$ is a sheaf.  In particular, the spectral
  sequence computing the (vanishing!) hypercohomology collapses
  immediately, and thus all homology groups vanish as required.
\end{proof}

This is particularly powerful for the following reason: the condition that
$(\Phi_V\Psi)^n M[-n]$ should be a sheaf is essentially {\em numerical} in
nature.  Indeed, it follows easily (by considering indecomposable summands
of $M$) that this holds iff it holds for every stable constituent of $M$.
But for a stable sheaf, we know that $(\Phi_V\Psi)^n M[-n]$ is a sheaf iff
$\rank((\Phi_V\Psi)^j M[-j])>0$ for $1\le j\le n$, iff
$(-1)^j\rank((\Phi_V\Psi)^j M)>0$ for $1\le j\le n$.

\begin{cor}
  With $V$, $M$ as above, suppose that $M$ is a semistable bundle.  Then
  the above conclusion holds iff the power series
  \[
  \frac{1 + (r^2\mu(M)+(m-dr))t/m}
       {1-(\deg({\cal L}) r^2/m-2)t+t^2}
  \]
  has positive coefficients, where $r=\rank(V)$, $d=\deg(V)$, and
  $m=\gcd(r,d)$ is its number of stable constituents.
\end{cor}

\begin{proof}
  We need $(-1)^j \rank((\Phi_V\Psi)^j M)>0$ for $j\ge 1$, or equivalently
  that the power series
  \[
  \sum_{j\ge 1} (-t)^j \rank((\Phi_V\Psi)^j M)
  \]
  has positive coefficients.  Since $\Phi_V\Psi$ acts linearly on the
  Grothendieck group, we may rephrase this as saying that all nonconstant
  coefficients of
  \[
  \rank(\sum_{j\ge 0} (-t)^j (\Phi_V\Psi)^j [M])
  =
  \rank((1+t\Phi_V\Psi)^{-1} [M])
  \]
  are positive.  The claim then follows by inverting the relevant $2\times
  2$ matrix and dividing by $\rank(M)$ to rewrite it in terms of the slope.
\end{proof}

\begin{rem}
  In the case that $M$ is torsion, we obtain a resolution iff
  \[
  \frac{1}
       {1-(\deg({\cal L}) r^2/m-2)t+t^2}
  \]
  has positive coefficients, or equivalently when
  \[
  \deg({\cal L}) r^2\ge 4m.
  \]
\end{rem}

\begin{lem}
  If $\tau\ge 2$, then the power series $(1+\alpha t)/(1-\tau t+t^2)$ has
  positive coefficients iff $\alpha\ge -(\tau+\sqrt{\tau^2-4})/2$.
\end{lem}

\begin{proof}
  We may write $\tau = \lambda+1/\lambda$ with $\lambda=(\tau+\sqrt{\tau^2-4})/2\ge 1$.  We then find that
  \[
  \frac{1-\lambda t}{1-(\lambda+1/\lambda)t+t^2} = \frac{1}{1-t/\lambda}
  \]
  has positive coefficients, as (apart from the constant term) does
  $\frac{t}{1-(\lambda+1/\lambda)t+t^2}$.  Moreover, the coefficients of the latter are
  unbounded, while the coefficients of the former go to 0, and thus the
  coefficients of
  \[
  \frac{1-(\lambda+\epsilon)t}{1-(\lambda+1/\lambda)t+t^2}
  \]
  are either always positive (if $\epsilon\le 0$) or eventually negative
  (if $\epsilon>0$), implying the desired result.
\end{proof}

\begin{rem}
  If $\tau\in \{-1,0,1\}$, positivity fails for any $\alpha$.
\end{rem}

\begin{thm}\label{thm:often_koszul}
  Let $V$ be a divisorial vector bundle on $E$ with rank $r$, degree $d$
  and $m=\gcd(r,d)$ stable constituents.  Then $B_{V,\Psi}$ is Koszul
  unless $m=r$ (so $r|d$) and $\deg({\cal L})r\le 3$.
\end{thm}

\begin{proof}
  If we apply the corollary to the case $M=V$, we find that the relevant
  power series is
  \[
  \frac{1+t}{1-(\deg({\cal L}) r^2/m-2)t+t^2},
  \]
  which has positive coefficients as long as the denominator has positive
  real roots.  This gives the condition
  \[
  \deg({\cal L}) (r/m)^2 m \ge 4
  \]
  which holds unless $r|d$ and $\deg({\cal L})r\le 3$.

  The resulting complex of $B_{V,\Psi}$-modules is exact in all positive
  degrees, and thus gives a {\em free} resolution of $(B_{V,\Psi})_0$.
  Since this resolution is manifestly linear, the claim follows.
\end{proof}

Something similar applies to minimal resolutions in the category ${\cal
  B}_{V,\Psi}$ more generally, again with the assumptions that all objects
are divisorial.  We may in fact generalize somewhat: we choose for each
$i\in \Z$ a divisorial sheaf $M_i$.  If $\mu(M_i)$ is an increasing
function of $i$, this gives a category as before, but even if not, we can
still define a {\em dg}-category:
\[
  {\cal B}_{\vec{M}}(i,j)
  =\begin{cases}
  0, & i>j\\
  \End(M_i), &i=j\\
  R\Hom(M_i,M_j) & i<j.
  \end{cases}
\]
(Here, of course, $R\Hom(M_i,M_j)$ should be replaced by the appropriate
complex.)  Let $P_i$ be the associated free modules (which are no longer
projective, per se), and $S_i$ be defined analogously to the $B_{V,\Psi}$
case.

Suppose now that $N$ is a semistable sheaf, and consider the
truncated module:
\[
M_{N,\ge 0}:i\mapsto \begin{cases} 0 & i>0\\
  R\Hom(M_i,N) & i\le 0.
\end{cases}
\]
over this dg-category.  The computation of the minimal resolution of
$M_{N,\ge 0}$ begins with the morphism
\[
P_0\otimes_{\End(M_0)}R\Hom(M_0,N)\to M_{N,\ge 0},
\]
apart from replacing $R\Hom(M_0,N)$ itself by a minimal resolution (as an
$\End(M_0)$-module).  Let $M'$ be the corresponding cone, so that we have
for each $d\in \Z$ a distinguished triangle
\[
R\Hom(P_{-d},P_0)\otimes_{\End(M_0)}R\Hom(M_0,N)\to R\Hom(P_{-d},M_{N,\ge
  0})\to M'(-d)\to.
\]
For $d<0$, the first two terms vanish, while for $d=0$ the map becomes
\[
\End(M_0)\otimes_{\End(M_0)}R\Hom(M_0,N)\cong R\Hom(M_0,M_{N,\ge
  0}),
\]
so that $M'(-d)=0$ for $d\le 0$, while for $d>0$ we may rewrite the
distinguished triangle as
\[
R\Hom(M_{-d},M_0)\otimes_{\End(M_0)}R\Hom(M_0,N)\to R\Hom(M_{-d},N)\to
M'(-d)\to.
\]
If we pull the tensor product inside the first $R\Hom$, we see that the
first morphism is the image under $R\Hom(M_{-d},{-})$ of the natural
morphism
\[
M_0\otimes_{\End(M_0)}R\Hom(M_0,N)\to N,
\]
and thus we find that
\[
M'(-d)\cong R\Hom(M_{-d},\Phi_{M_0}N).
\]
In other words, we have a distinguished triangle
\[
P_0\otimes_{\End(M_0)}R\Hom(M_0,N)\to M_{N,\ge 0} \to M_{\Phi_{M_0}N,\ge 1}\to,
\]
from which it is straightforward to inductively construct the minimal
resolution of $M_{N,\ge 0}$.  In particular, we find that the associated
graded of the resulting filtered complex is the direct sum
\begin{align}
&P_0\otimes_{\End(M_0)} R\Hom(M_0,N)\notag\\
\oplus&
P_1\otimes_{\End(M_1)} R\Hom(M_1,\Phi_{M_0}N)\notag\\
\oplus&
P_2\otimes_{\End(M_2)} R\Hom(M_2,\Phi_{M_1}\Phi_{M_0}N)\notag\\
\oplus&
\cdots,
\end{align}
again with each $R\Hom$ complex replaced by a suitable minimal resolution.

The one caveat above is that it is not a priori obvious that $R\Hom(M_0,N)$
admits a finite projective resolution.  If $\mu(N)>\mu(M_0)$, then this is
just $\Hom(M_0,N)$, which we have shown is already free, and similarly if
$\mu(N)<\mu(M_0)$, we know by duality that $\Ext^1(M_0,N)$ is free.  Thus
only the case that equality holds is an issue.  We may reduce to the case
that $M_0$, $N$ are both structure sheaves of divisors, and from that
immediately reduce to the case that they are structure sheaves of jets
supported at the same point, and may then pass to the complete local ring.
In other words, we need to know that
\[
R\Hom_{k[[x]]}(k[[x]]/(x^n),k[[x]]/(x^m))
\]
is represented by a finite complex of free $k[[x]]/(x^n)$-modules.  Using
the two-term injective resolution
\[
0\to k[[x]]/(x^m)\to k((x))/x^m k[[x]]\to k((x))/k[[x]]\to 0,
\]
we obtain a two-term free resolution, given by the morphism
\[
x^m:k[[x]]/(x^n)\to k[[x]]/(x^n).
\]
(This morphism vanishes on $k[[x]]/(x)$, so this is the minimal
resolution.)  It follows more generally that if $N$ with $\mu(N)=\mu(M_0)$
is also divisorial, then $R\Hom(M_0,N)$ is represented by a two-term
complex in which both terms are the same direct summand of $\End(M_0)$.
Note that we can also obtain a (non-minimal) free resolution in this case,
by taking the direct sum of this complex with a complex having morphism $1$
on a complementary direct summand.

We thus obtain the following.

\begin{prop}
  The minimal resolution of the module $S_i$ is given by a filtered complex
  such that the associated graded is
  \[
  P_i
  \oplus
  P_{i-1}\otimes R\Hom(M_{i-1},M_i)[1]
  \oplus
  P_{i-2}\otimes R\Hom(M_{i-2},\Phi_{M_{i-1}}M_i)[1]
  \oplus
  \cdots.
  \]
  In particular, each $P_i$ either appears in a single cohomological
  degree, tensored with a free module, or in two consecutive degrees,
  tensored with the same direct summand of $\End(M_i)$.
\end{prop}

In particular, we can compute the structure of the minimal resolutions from
purely numerical data.  Indeed, the object
\[
\Phi_{M_{j+1}}\cdots \Phi_{M_{i-1}} M_i
\]
is always a shift of a sheaf (by induction using Proposition
\ref{prop:phiM_preserves_shifted_sheaves}), with the amount of shift
monotone in $j$ and increasing by at most 1 when we decrease $j$.
Moreover, the jumps occur precisely when the pair $(\rank,\deg)$ switches
to the other side of $(0,0)$ in lexicographic ordering, so depend only on
the image of $(\rank(M_i),\deg(M_i))$ under the relevant product of
elements of $\SL_2(\Z)$.  The support of
\[
R\Hom(M_j,\Phi_{M_{j+1}}\cdots \Phi_{M_{i-1}} M_i)
\]
is then determined by adding $0$, $1$, or $[0,1]$ to that shift, depending
on how the two slopes compare.

The above structure of the minimal resolution has the following
consequence, using the fact that $R\Hom(P_i,S_j)=0$ unless $i=j$, when
$R\Hom(P_i,S_j)=\End(M_j)$.

\begin{prop}
  One has
  \[
  R\Hom(S_i,S_j)\cong
  \begin{cases}
    0, & i<j\\
    \End(M_i), & i=j\\
    R\Hom_{\End(M_j)}(
    R\Hom(M_j,\Phi_{M_{j+1}}\cdots\Phi_{M_{i-1}}M_i),\End(M_j))[-1], & i>j.
  \end{cases}
  \]
\end{prop}

\begin{rem}
  If we define a new sequence of (shifts of) divisorial sheaves by
  \[
  M'_i = \Phi_{M_0}\cdots \Phi_{M_{i-1}}M_i
  \]
  for $i\ge 0$ and
  \[
  M'_i = \Phi_{M_{-1}}^{-1}\cdots \Phi_{M_i}^{-1}M_i
  \]
  for $i\le 0$, then one has a non-canonical isomorphism
  \[
  R\Hom(M'_i,M'_j)
  \cong
  R\Hom_{\End(M'_j)}(R\Hom(M'_j,M'_i),\End(M'_j))[-1],
  \]
  which together with the canonical isomorphism
  \[
  R\Hom(M'_j,M'_i)
  \cong
  R\Hom(M_j,\Phi_{M_{j+1}}\cdots\Phi_{M_{i-1}}M_i)
  \]
  for $i\ge j$ produces an isomorphism
  \[
  R\Hom(S_i,S_j)\cong
  \begin{cases}
    0, & i<j\\
    \End(M'_i), & i=j\\
    R\Hom(M'_i,M'_j), & i>j.
  \end{cases}
  \]
  Presumably this quasi-isomorphism of complexes can be made compatible
  with composition.  This is certainly true in any of the cases to which
  Theorem \ref{thm:often_koszul} above applies, when the dg-category with
  Hom complexes $R\Hom(S_i,S_j)$ is essentially just a $\Z$-algebra version
  of the Koszul dual.  The general conjecture is thus that if we relax
  $B_{\vec{M}}$ to allow (cohomological) {\em shifts} of divisorial
  sheaves, then the Koszul dual of $B_{\vec{M}}$ is a dg-category of the
  same form, and if $B_{\vec{M}}$ is periodic up to derived autoequivalence
  (letting us turn the $\Z$-dg-algebra structure into a dg-algebra
  structure with an additional grading), then so is its Koszul dual.
\end{rem}

\begin{eg}
Consider the case $V=\sO_E$, $\deg({\cal L})=1$ of our $B_{V,\Psi}$
construction.  This is already not generated in degree 1, since
$\Phi_{V_{-1}}(V_0)$ is a sheaf, the structure sheaf of a point.  We do,
however, get a short exact sequence
\[
0\to \Phi_{V_{-2}}\Phi_{V_{-1}}(V_0)[\pm 1]\to V_{-2}\oplus V_{-1}\to
V_0\to 0
\]
and find that
\[
\Phi_{V_{-2}}\Phi_{V_{-1}}(V_0)[\pm 1]\cong V_{-2}\otimes V_{-1}.
\]
In particular, it has slope $-3$, so the standard resolution requires us to
add the generic morphism from $V_{-3}$ to $V_0$.  The resulting complex
\[
V_{-3}\to V_{-3}\oplus V_{-2}\oplus V_{-1}\to V_0
\]
represents a shift of $\Phi_{V_{-3}}\Phi_{V_{-2}}\Phi_{V_{-1}}V_0$, and the
corresponding complex of ${\cal B}$-modules represents the truncation of
this module to degree $\ge 4$ (i.e., to the functor taking $V_i$ to 0 if
$i>-4$).  Apart from a shift by $-3$, this is numerically the same as the
case we started with, and we thus obtain a free resolution of $k$ in the
following periodic form:
\[
\cdots
\to
   [-12,-11,-10,-9]\to [-9,-8,-7,-6]\to [-6,-5,-4,-3]\to [-3,-2,-1]\to [0],
\]
where $[i,j,\dots,k]$ denotes $P_i\oplus P_j\oplus\cdots\oplus P_k$.  If
$\tau^2=\text{id}$, then this is minimal, and otherwise the minimal
resolution has shape
\[
\cdots
\to [-14,-13,-12] \to [-12,-11,-10]\to [-8,-7,-6]\to [-6,-5,-4]\to [-2,-1]\to [0].
\]

Similarly, for the case $V=\sO_E$, $\deg({\cal L})=2$, the resolution has
shape
\[
\cdots \to [-8,-7,-7,-6]\to [-6,-5,-5,-4]\to [-4,-3,-3,-2]\to [-2,-1,-1]\to [0]\to 0,
\]
reducing to
\[
\cdots \to [-9,-9,-8] \to [-8,-7,-7]\to [-5,-5,-4]\to [-4,-3,-3]\to [-1,-1]\to [0]\to 0
\]
unless $\tau^2=\text{id}$.  For $\deg({\cal L})=3$, the resolution has
shape
\begin{align}
\cdots&\to [-7,-7,-7,-6]\to [-6,-5,-5,-5]\notag\\
&\to [-4,-4,-4,-3]\to [-3,-2,-2,-2]\to
      [-1-1,-1]\to [0]\to 0
\end{align}
and is always minimal.
\end{eg}

An interesting common property of these three cases is that, in each case,
if we omit the highest-degree relation, the resulting algebra still has a
nice resolution, and in particular has the same Hilbert series as a free
polynomial ring in the generators, with the highest-degree relation
corresponding to a {\em central} element of this ring. (This will follow
once we understand the filtered deformations.  In each case, the graded
algebra of deformation parameters is a free polynomial ring with one
generator of the same degree as the highest-degree relation.  So if we set
the other deformation parameters to 0, we get a deformation in which the
lower-degree relations are undeformed while the top-degree relation has a
(nonzero) scalar added.  The commutator of the deformed relation with any
generator is thus expressible in terms of the undeformed relations, giving
the requisite centrality statement, and letting us recover the Hilbert
series.)

The other two cases in which the algebra fails to be Koszul have
$\deg({\cal L})=1$ and $V$ a slope 0 divisorial bundle with rank $2$ or
$3$.  (This is generically a sum of line bundles of degree 0, except that
when the line bundles become isomorphic, they should be replaced by the
appropriate self-extension.)  Here the free resolution has the same shape
as in the previous two cases:
\[
\cdots \to [-8,-7,-7,-6]\to [-6,-5,-5,-4]\to [-4,-3,-3,-2]\to [-2,-1,-1]\to [0]\to 0
\]
or
\begin{align}
\cdots&\to [-7,-7,-7,-6]\to [-6,-5,-5,-5]\notag\\
&\to [-4,-4,-4,-3]\to [-3,-2,-2,-2]\to
      [-1-1,-1]\to [0]\to 0;
\end{align}
the main difference is that minimality can fail in more (and more
complicated) ways.  Again, in both of these cases, there is a deformation
with a single parameter of degree 4 or 3 as appropriate, thus expressing
the algebra as a central quotient of another algebra with well-behaved
Hilbert series.

\medskip

One question suggested by the above arguments is: given that we can compute
resolutions of $B_{V,\Psi}$ purely in terms of sheaves on $E$, can we do
something similar for infinitesimal deformations?  One could settle a
number of cases of Conjecture \ref{conj:main} below (and in many ways the
most interesting cases) if one could compute the negative degree part of
$\HH^2(B_{V,\Psi})$.  Certainly, having an explicit resolution of
$B_{V,\Psi}$ makes it feasible to compute this group in special cases, but
one would like to do so in general (or at least by hand).

\section{Elliptic noncommutative del Pezzo surfaces}

We have seen in Corollary \ref{cor:fildef_is_rigid_unobstructed} that any
filtered deformation of an elliptic algebra $B_{V,\Psi}$ is associated to a
sheaf on the corresponding compactification (suitably defined) with no
higher endomorphisms.  This suggests a way to construct such deformations
by starting with a suitable sheaf on a (noncommutative) surface.  Since (by
the same Corollary) the compactification is Gorenstein with dualizing sheaf
$\sO(-1)$, this surface should be del Pezzo in a suitable sense.  In
general, of course, this would be a quite broad sense, since the cone over
$B_{V,\Psi}$ satisfies the same conditions, but this at least suggests that
we should obtain a large class of deformations by taking more
familiar-looking del Pezzo surfaces.  Of course, if the del Pezzo surfaces
are {\em too} familiar (i.e., commutative), then this will constrain $\Psi$
to be the tensor product with a line bundle, and thus fail to explain why
various properties of the moduli space of filtered deformations extend to
more general $\Psi$.

Luckily, there is a theory of noncommutative surfaces available
\cite{noncomm2}, including a large subclass coming from elliptic curves
\cite{noncomm1}.  This gives two classes of direct constructions, both of
which are rather involved.  Note that in contrast to the commutative
theory, there is no direct geometric interpretation of a noncommutative
surface; what one actually studies in general is a category of
(quasi-)coherent sheaves on such a surface (which in our case is a flat
deformation of the category of (quasi-)coherent sheaves on a commutative
surface).  It turns out that the corresponding {\em derived} category is
quite simple to describe.  For simplicity, we will restrict our attention
to coherent sheaves (and bounded complexes with coherent cohomology); in
particular, here and below, all sheaves on noncommutative surfaces will be
assumed coherent.  Note that since the noncommutative surfaces we consider
are Noetherian \cite[Prop. 6.1]{noncomm2}, coherent sheaves are always
Noetherian.

In general, if $X$ is a smooth projective rational surface, then the
structure sheaf $\sO_X$ is exceptional, and thus determines a
semiorthogonal decomposition $(\sO_X^\perp,\sO_X)$ of $D^b_{\coh}(X)$.  If
$X$ has an anticanonical curve $Q$ (e.g., if $X$ is del Pezzo), with
corresponding morphism $i:Q\to X$, then one has a functorial distinguished
triangle
\[
  {-}\otimes \omega_X\to \text{id}\to i_*i^*\to
\]
and thus in particular a distinguished triangle
\[
R\Hom(M,N\otimes \omega_X)\to R\Hom(M,N)\to R\Hom_Q(M|^L_Q,N|^L_Q)\to
\]
for any objects $M$, $N\in D^b_{\coh}(X)$.  (One caveat is that the natural
transformation ${-}\otimes\omega_X\to\text{id}$ is only determined up to a
scalar by the choice of $Q$, with that scalar corresponding to a choice of
isomorphism $H^1(\sO_E)\cong H^2(\omega_X)\cong k$.) In particular, if
$R\Hom(N,M)=0$, then Serre duality gives $R\Hom(M,N\otimes \omega_X)=0$ and
thus $R\Hom(M,N)\cong R\Hom_Q(M|^L_Q,N^L|_Q)$.  As a special case, if $M\in
\sO_X^\perp$, then $R\Hom(M,\sO_X)\cong R\Hom_Q(M|^L_Q,\sO_Q)$.

The key observation is that this lets one reconstruct $D^b_{\coh}X$ from the
(dg-)subcategory $\sO_X^\perp$ and the restriction-to-$Q$ functor on that
subcategory.  Moreover, the dg-functor
\[
M\mapsto R\Hom_Q(M|^L_Q,\sO_Q)
\]
admits a natural flat deformation: given an invertible sheaf $q$ in the
identity component of $\Pic(X)$, we may define $D^b_{\coh}(X_q)$ to be the
dg-category obtained by ``gluing'' \cite{OrlovD:2016} $\sO_X^\perp$ to $\sO_X$
via the functor
\[
M\mapsto R\Hom_Q(M|^L_Q,q)
  \cong  R\Hom_Q(M|^L_Q\otimes q^{-1},\sO_Q).
\]
(That this is the derived category of a noncommutative surface is shown in
\cite[Thm. 4.30]{noncomm2}.)  Note that since the forward maps in the
semiorthogonal decomposition come from $Q$, the functor $Li^*$ immediately
extends to $D^b_{\coh}(X_q)$, with $Li^*\sO_{X_q}\cong \sO_Q$ and the
restriction on $\sO_X^\perp$ twisted by $q^{-1}$.

In particular, we see that the Grothendieck group of $D^b_{\coh}(X_q)$ is
independent of $q$ (and largely independent of $X$ itself); this lets one
define the rank, first Chern class, and Euler characteristic of a general
object (by taking the corresponding homorphisms on $K_0(X)$).  Moreover,
since this is a flat deformation, the Mukai pairing $\chi([M],[N]):=\chi
R\Hom(M,N)$ is also independent of $q$, so that one immediately has an
analogue of Hirzebruch-Riemann-Roch \cite[Cor. 7.13]{noncomm2} for objects
in $D^b_{\coh}(X_q)$:
\[
\chi(M,N)
=
-\rank(M)\rank(N)
+\rank(M)\chi(N)+\rank(N)\chi(M)
-c_1(M)\cdot (c_1(N)+\rank(N)Q).
\]
We also obtain Serre duality, i.e., the existence of a (canonical)
dg-functor $S$ with a (suitably coherent) isomorphism
\[
R\Hom(M,N)\cong R\Hom(R\Hom(N,SM),k).
\]
(This follows from standard facts about the interactions of Serre functors
with semiorthogonal decompositions.)  Moreover, $Li^*$ has a right adjoint
$Ri_*$ fitting into a functorial distinguished triangle
\[
\theta \to \text{id}\to Ri_*Li^*\to
\]
with $\theta:=S[-2]$ (the analogue of twisting by the anticanonical
bundle).  (This distinguished triangle is again not quite canonical, as the
natural transformation $\theta\to\text{id}$ is only determined up to scalar
multiple.)  There is also \cite[Prop. 5.1]{noncomm2} an analogue of
Cohen-Macaulay duality, an equivalence $\ad_q:(D^b_{\coh}
X_q)^{\text{op}}\cong D^b_{\coh}(X_{q^{-1}})$ satisfying
$\ad_{q^{-1}}\ad_Q\cong \text{id}$, acting on $K_0(X)$ as
$R\Hom({-},\omega_X)$.

When $Q$ is smooth, the description of the derived category can be
simplified further by observing that $\sO_X^\perp$ admits a full
exceptional collection.  Indeed, in that case we can nearly always write
$X$ as an iterated blowup of $\P^2$:
\[
X\cong X_m\to X_{m-1}\to\cdots\to X_{-1}\cong \P^2
\]
where each map $X_i\to X_{i-1}$ blows down a single $-1$-curve.  (We take
the sequence to end at $-1$ rather than $0$ to be compatible with the
general case discussed in \cite{noncomm2}, in which one must typically
consider iterated blowups of ruled surfaces.)  Then we have a full
exceptional collection on $X$ of the form
\[
(\sO_{e_m}(-1),\sO_{e_{m-1}}(-1),\dots,\sO_{e_0}(-1),\sO_{\P^2}(-2),\sO_{\P^2}(-1),\sO_{\P^2}),
\]
where we omit the pullback functors from $X_i$ to $X_m$.  (The only case
not covered by iterated blowups of $\P^2$ is when $X$ is an even Hirzebruch
surface $F_0$ or $F_2$, in which case we have a full exceptional collection
of the form
\[
(\sO_X(-s-f),\sO_X(-s),\sO_X(-f),\sO_X),
\]
where $f$ is the divisor class of a fiber of the ruling, while $s$ is the
unique divisor class satisfying $s\cdot f=1$, $s^2=0$.  The analogue of the
following discussion for that case is straightforward.)

We then see that $D^b_{\coh}(X_q)$ may be obtained from the sequence of
restrictions of the above sheaves to $Q$: it is the category of perfect
modules over a triangular dg-algebra obtained from
\begin{align}
R\End(&
\sO_{e_m}(-1)|_Q\oplus
\sO_{e_{m-1}}(-1)|_Q\oplus
\cdots\oplus
\sO_{e_0}(-1)|_Q\oplus{}\notag\\
&
\sO_{\P^2}(-2)|_Q\otimes q^{-1}\oplus
\sO_{\P^2}(-1)|_Q\otimes q^{-1}\oplus
\sO_Q
)
\end{align}
by omitting all ``backwards'' morphisms and all nontrivial endomorphisms of
summands.  Since each summand is a simple coherent sheaf on $Q$, it is
determined up to isomorphism by its class in $K_0(Q)$, and thus we see that
$D^b_{\coh}(X_q)$ is determined by the restriction morphism $K_0(X)\to
K_0(Q)$.  In fact, for any morphism $\phi:K_0(X)\to K_0(Q)$ satisfying
\[
\rank(\phi(M))=\rank(M),\qquad
\deg(\phi(M))=c_1(M)\cdot Q,\qquad
\phi([\sO_X])=[\sO_Q],
\]
there is a corresponding sequence of simple modules on $Q$ and thus a
corresponding dg-category.  Moreover, if we set
\[
q=
\det(\phi([\sO_{\P^2}(-2)]))
\otimes \det(\phi([\sO_{\P^2}(-1)]))^{-2},
\]
then the dg-category has the form $D^b_{\coh}(X_q)$.  (Indeed, the surface $X$
is obtained from this data by first embedding $Q$ in $\P^2$ using the line
bundle $(\det(\phi([\sO_{\P^2}(-1)]))\otimes q)^{-1}\in \Pic^3(Q)$, then
sequentially blowing up the points corresponding to
$\det(\phi([\sO_{e_i}(-1)]))\in \Pic^1(Q)$.)

If we apply $\theta$ to the sum of exceptional objects before restricting
to $Q$, we of course obtain the same endomorphism dg-algebra, and thus
(since the summands generate $D^b_{\coh}(Q)$) there is necessarily a
derived autoequivalence $\xi$ of $Q$ such that $(\theta M)|^L_Q\cong
\xi(M|^L_Q)$.  Moreover, since the action of $\theta$ on $K_0(X_q)$ follows
from the commutative case by flatness, we may deduce the action of $\xi$ on
$K_0(Q)$.  We thus find that $\xi$ has the form $\Psi^{-1}$ where $\Psi$ is
an {\em abelian} autoequivalence of the kind considered above.

A choice of blowdown structure induces a basis the N\'eron-Severi lattice
of $X$ of the form $h,e_0,\dots,e_m$, with intersection form given by
\[
h^2=1,\quad h\cdot e_i=0,\quad e_i\cdot e_j=-\delta_{ij},
\]
with $Q=3h-e_0-\cdots-e_m=-K_X$.  (In the even Hirzebruch case, this
becomes $s^2=f^2=0$ and $s\cdot f=1$ with $Q=2s+2f$.)  Changing the
blowdown structure by swapping two consecutive blowups acts on the basis as
the reflection in $e_i-e_{i+1}$, while a quadratic transformation in the
first three points acts as the reflection in $h-e_0-e_1-e_2$; in either
case, the result corresponds to a blowdown structure iff the corresponding
$-2$-class is ineffective.

One can show that any two blowdown structures are related by a sequence of
such operations: see, e.g., \cite[Prop. 8.10]{noncomm2} for the analogous
result for arbitrary blowups of noncommutative ruled surfaces.
(Analogously, for even Hirzebruch surfaces, changing the marked ruling acts
as the reflection in $s-f$.)

These reflections are the simple reflections of a natural Coxeter group
action.  The point is that since $Q^2>0$, the Hodge index theorem tells us
that the lattice $Q^\perp:=\{v:v\in \Lambda|v\cdot Q=0\}$ is negative
definite relative to the intersection form, so negating the pairing gives a
positive definite integral lattice, which is in particular acted on by
reflections in any vectors of square norm $2$.  For $Q^2=9-m\le 6$, the
lattice is in fact {\em generated} by such vectors, and one recovers
familiar root lattices:
\[
E_8,E_7,E_6,D_5,A_4,A_2A_1
\]
(of rank $m=9-Q^2$).  (For $m=2$, one obtains the unique rank 2 even
lattice of determinant 7, while for $m=1$, the rank 1 lattice depends on
the parity of the Hirzebruch surface $X$.)  In the resulting finite root
system, we can take the positive roots to be those which are
lexicographically positive, in which case the simple roots are precisely
$h-e_0-e_1-e_2$, $e_0-e_1$,\dots  One then finds that a simple root is
effective iff its image in $\Pic^0(Q)$ is trivial.

\medskip

Given this description of the derived category, we can then recover the
actual category of coherent sheaves $\coh(X_q)$ from a suitable description
of the $t$-structure.  The simplest description involves first definining
``line bundles'', inductively defined \cite[Defn. 5.1]{noncomm2} by
starting with $\sO_{X_{-1,q}}(d)$ for $d\in \{-2,-1,0\}$ (i.e., the three
terms in the semiorthogonal decomposition) and taking the closure under
$\theta^{\pm 1}$ and pullback from $X_{l-1,w}$ to $X_{l,q}$.  (For the even
Hirzebruch surface case, we note that the subcategories generated by the
first two and last two exceptional sheaves in the above exceptional
collection generate copies of $D^b_{\coh}(\P^1)$, and we take as our line
bundles the images under powers of $\theta$ of line bundles on $\P^1$ in
those two subcategories.)  We caution that the resulting collection of
sheaves depends on the choice of blowdown structure (and the choice of
ruling when $X\cong \P^1\times \P^1$) in general.  (The precise dependence
is that we may only apply a reflection in a simple root when that simple
root is ineffective on the {\em noncommutative} surface, which happens
unless the corresponding class in $\Pic^0(Q)$ is a power of $q$.  In
particular, each simple root is {\em generically} ineffective, but is
effective over a dense collection of hypersurfaces in the moduli space of
surfaces.)

Then \cite[Lem. 5.7]{noncomm2}, an object $M\in D^b_{\coh}(X)$ is in
$D^b_{\coh}(X)^{\ge 0}$ iff $R\Hom(L,M)\in D^b_{\coh}(k)^{\ge 0}$ for every
line bundle $L$; this fully determines the $t$-structure, since a
$t$-structure is determined by its nonnegative part.  Moreover, line
bundles are themselves in the heart of this $t$-structure, and in fact
generate the heart.  Furthermore, nonzero morphisms between line bundles
are injective \cite[Prop. 7.21]{noncomm2}.  (The construction of
\cite{noncomm1} is via a representation of the category of line bundles in
which the morphisms are expressed as elliptic difference operators, and
such operators cannot be zero divisors.)  Furthermore, by the discussion
following \cite[Defn. 5.1]{noncomm2}, there is for any $D$ an abelian
equivalence of the form $\coh(X_q)\cong \coh(X'_q)$ acting on $K_0$ as
${-}(D)$ and taking $\sO_{X_q}(D')$ to $\sO_{X'_q}(D'+D)$ for all $D'$.
Note that since twisting by a line bundle changes the surface, there are
far fewer graded algebras associated to a noncommutative surface than in
the commutative case, and even fewer maps to noncommutative projective
spaces.

In particular, we define an {\em elliptic noncommutative (degenerate) del
  Pezzo surface} of degree $d>0$ to be an abelian category constructed from
a map $K_0(X)\to K_0(Q)$ as above, where $Q$ is a smooth genus 1 curve and
$K_0(X)$ is the Grothendieck group of a del Pezzo surface.  (It will be
convenient also to fix the scale of the natural transformation
$\alpha:\theta\to\text{id}$.)  When $q=\sO_Q$, such a surface corresponds
to a pair $(X,Q)$, where $X$ is a usual (smooth!) degenerate del Pezzo
surface and $Q$ is a smooth anticanonical curve.  Note that since these are
analogues of {\em smooth} surfaces, they are generally not Fano; in the
commutative case, the anticanonical embedding contracts $-2$-curves, and
the situation is worse in the noncommutative case.  (In general, if there
is {\em some} commutative del Pezzo surface with a $-2$-curve of class $D$,
then there is a pure 1-dimensional sheaf of class $[\sO_D(l)]$ on $X_q$ iff
$\phi[\sO_D(l)]=0$, in which case $Q$ fails to be ample.  And we have
already seen that this happens for a dense set of $q$!)

We may then use our understanding of the Grothendieck group to extend
various concepts from the commutative case: a sheaf is {\em torsion} if it
has rank 0 (note that every sheaf has nonnegative rank), a divisor class is
{\em effective} if it is the first Chern class of a torsion sheaf (and
\cite[Prop. 7.19]{noncomm2} if both $D$ and $-D$ are effective, then
$D=0$), a divisor class is {\em nef} if it is in the dual cone to the
monoid of effective divisors, and ample if it is in the interior of the nef
cone, with \cite[Thm. 9.10]{noncomm2} ample divisors satisfying
\[
\Ext^p(\sO_X(-nD),M)=0
\]
for $p>0$ and $\sO_X(-nD)\otimes \Hom(\sO_X(-nD),M)\to M$ surjective for
$n\gg 0$.  Note that in our case, $Q$ is always nef: the effective monoid
is generated by components of $Q$ along with those $-1$- and $-2$-classes
which are effective \cite[Thm. 9.2]{noncomm2}, and thus since $Q$ is
irreducible with $Q^2\ge 0$, it is nef.

Given an ample divisor, there is a corresponding notion of Hilbert
polynomial ($p_D(M):=\chi(\sO_X(-nD),M)$, which we note is $0$ iff $M=0$
and is otherwise positive for $n\gg 0$), which one may use to define
(semi)stable sheaves.  One then finds \cite[Prop. 11.34]{noncomm1} that the
corresponding semistable moduli spaces are proper.  (One caveat, however:
although it is likely that semistable moduli spaces are projective in
general, this has only been proved for $\rank\le 1$.)

In general, there are relatively few concepts for commutative surfaces that
do not immediately carry over to these noncommutative surfaces.  The most
basic issue is that the {\em support} of a torsion sheaf is not
well-defined in general.  However, since the Chern class {\em is}
well-defined, we can often work around this (e.g., replacing the condition
that a sheaf is supported on an exceptional curve $e$ by the condition that its
Chern class is a multiple of the divisor class $e$).  A related issue is
that noncommutative surfaces have relatively few points; indeed, on a
typical noncommutative surface, the only $0$-dimensional sheaves (i.e.,
with vanishing rank and first Chern class) are in the Serre subcategory
generated by $0$-dimensional sheaves on $Q$.  We have already mentioned the
dearth of embeddings in projective space, which breaks arguments involving
taking hyperplane sections.  In addition, the category of sheaves on a
noncommutative surface is no longer a closed monoidal category (i.e., there
is no internal Hom or tensor product).

Finally, and most significant for our purposes, although {\em line bundles}
make sense on a noncommutative surface, there is no fully satisfying
definition of {\em vector bundles}.  For many purposes (and indeed below)
it suffices to consider torsion-free sheaves.  One could also consider {\em
  reflexive} sheaves, i.e., sheaves $M$ on $X_q$ such that $\ad_q M$ is a
sheaf on $X_{q^{-1}}$.  (As in the commutative case, such sheaves are
necessarily torsion-free.)  This has a number of pathologies in general
(e.g., there can be reflexive sheaves of rank 1, trivial Chern class, and
abitrarily small Euler characteristic), but for present purposes those can
be ignored, as they do not arise for {\em rigid} reflexive sheaves (i.e.,
with $\Ext^1(M,M)=0$).  In any event, although one can show that the
sheaves we construct below are reflexive, we will not need this fact.

Note that when $M$ is torsion-free, $\alpha:\theta M\to M$ is injective, as
the kernel is supported on $Q$, so torsion.  In particular, we see that
$M|_Q$ is a sheaf.  (When $M$ is reflexive, we have $(\ad M)|_Q\cong
\sHom(M|_Q,\sO_Q)$ and thus $M|_Q$ is a vector bundle.)

\section{Filtered deformations from del Pezzo surfaces}

We return to our considerations of filtered deformations.  Let $X$ be a
noncommutative elliptic del Pezzo surface, with chosen anticanonical curve
$Q$, adjoint functors $i_*:\coh(Q)\to \coh(X)$, $i^*:\coh(X)\to \coh(Q)$,
and a choice of natural transformation $\alpha:\theta\to \text{id}$ giving
the usual functorial distinguished triangle
\[
\theta\to \text{id}\to i_*Li^*\to.
\]
Note that the choice of alpha can be encapsulated in the observation that
there is a canonical natural transformation $\theta\otimes H^1(\sO_Q)\to
\text{id}$.  (In principle, one could develop much of the theory below even
when $Q$ is singular, but we assume $Q$ smooth for simplicity.  This, of
course, is the only case that can produce filtered deformations of {\em
  elliptic} algebras, and is also (see Example \ref{eg:weird_nodal} above)
the only case in which we can hope to obtain {\em all} filtered
deformations using such surfaces.  In addition, the argument that the
moduli stack is proper breaks down when $Q$ is singular, though the claim
appears to still be valid; if so, one can obtain algebras in degenerate
cases by taking limits.)  Recall that $\theta^{-1}$ induces an
autoequivalence of $\coh(Q)$ of the form $\Psi$ considered above, with
associated line bundle of degree $Q^2$.

By Corollary \ref{cor:fildef_is_rigid_unobstructed}, we expect there to be
a relation between filtered deformations on elliptic algebras and
homogeneous endomorphism rings of sheaves on del Pezzo surfaces with no
higher endomorphisms, and this is indeed the case.

\begin{prop}\label{prop:fildef_from_sheaf}
  Let $(X,Q)$ be a noncommutative elliptic del Pezzo surface, let $M\in
  \coh(X)$ be a torsion-free sheaf with $\Ext^1(M,M)=\Ext^2(M,M)=0$ such
  that the slopes of $M|_Q$ are constrained to a half-open interval of
  length $Q^2$, and consider the graded algebra
  \[
  A^+_M := \bigoplus_i \Hom(M,\theta^{-i}M).
  \]
  This algebra is nonnegatively graded, the element
  $\alpha_{\theta^{-1}M}:M\to \theta^{-1}M$ of degree 1 is regular and
  central, and the quotient by $\alpha_{\theta^{-1}M}$ is
  $B_{M|_Q,\theta^{-1}}$.
\end{prop}

\begin{proof}
  That $\alpha_{\theta^{-1}M}$ is regular and central in $A^+_M$ follows
  immediately from the facts that it is injective and comes from a natural
  transformation.  For each $i\in \Z$, $\alpha$ induces a distinguished
  triangle
  \[
  R\Hom(M,\theta^{1-i}M)\to R\Hom(M,\theta^{-i}M)\to
  R\Hom(M|_Q,\theta^{-i}M|_Q)\to.
  \]
  For $i=0$, we may use duality to compute $R\Hom(M,\theta M)$, and thus
  (using the hypotheses on $M$) we obtain isomorphisms
  \[
  \End(M)\cong \End(M|_Q),\quad
  \Ext^2(M,\theta M)\cong \Ext^1(M|_Q,M|_Q).
  \]
  Moreover, since $\alpha_{\theta^iM}$ is injective for all $i$, the
  assumption that $\Ext^2(M,M)=0$ implies that $\Ext^2(M,\theta^i M)=0$ for
  all $i\le 0$, and thus that $\Hom(M,\theta^i M)=0$ for $i\ge 1$.  In
  addition, the condition on $M|_Q$ ensures that $\Ext^1(M|_Q,\theta^i
  M|_Q)=0$ for $i>0$.  We thus inductively find that $\Ext^1(M,\theta^i
  M)=0$ for $i\le 0$, giving us the desired short exact sequence
  \[
  0\to \Hom(M,\theta^{1-i}M)\to \Hom(M,\theta^{-i}M)\to
  \Hom(M|_Q,\theta^{-i}M|_Q)\to 0
  \]
  for all $i$.
\end{proof}

\begin{rem}
  Of course, the quotient by $\alpha_{\theta^{-1}M}$ is independent of the
  choice of scaling on $\alpha_{\theta^{-1}M}$.  In addition, although
  $\theta^{-1}$ appears in the definition of $A^+_M$, we can safely omit it
  from the notation because $\theta$ is canonical (in both senses of the
  word).
\end{rem}

Note that if we instead quotient by $\alpha_{\theta^{-1}M}-1$, then the
result will be a filtered algebra, and the Proposition tells us that the
associated graded algebra is $B_{M|_Q,\theta^{-1}}$.  In other words, the
quotient is precisely a filtered deformation of the type we are looking
for.  The filtered deformation depends only mildly on the choice of scale
on $\alpha$: the action of $\G_m$ by rescaling $\alpha$ agrees with the
standard action on filtered deformations.

We are thus left with the question of classifying such sheaves.  In the
commutative case, sheaves with $\Ext^1(M,M)=0$ on degenerate del Pezzo
surfaces were studied in \cite{KuleshovSA:1997}; though the classification
was only partial and some of the arguments do not apply in the
noncommutative setting, we can use the ideas therein with some supplements
to get a fairly complete understanding in general.

A useful simplification comes from the following fact.

\begin{prop}
  Suppose $M$, $N$ are sheaves on a noncommutative del Pezzo surface
  $(X,Q)$ such that $\Ext^1(N,M)=0$ and $N$ is torsion-free.  Then
  $\Hom(M,N)\to \Hom(M|_Q,N|_Q)$ is surjective.
\end{prop}

\begin{proof}
  We have the exact sequence
  \[
  0\to \Hom(M,\theta N)\to \Hom(M,N)\to \Hom(M|_Q,N|_Q)\to \Ext^1(M,\theta
  N),
  \]
  but duality gives
  \[
  \Ext^1(M,\theta N)\cong \Ext^1(N,M)^*=0.
  \]
\end{proof}

\begin{cor}
  If $M$ is a torsion-free sheaf with $\Ext^1(M,M)=0$, then $\End(M)\to
  \End(M|_Q)$ is surjective, and any idempotent in $\End(M|_Q)$ lifts to
  $\End(M)$.
\end{cor}

\begin{proof}
  Surjectivity is just the case $N=M$ above, and thus $\End(M|_Q)$ is the
  quotient of $\End(M)$ by the ideal $\alpha_M\circ \Hom(M,\theta M)$.
  Since idempotents lift through nil ideals, it remains to show that any
  element $\beta$ of the ideal is nilpotent.

  Since $X$ is Noetherian and $M$ is coherent, the ascending chain
  $\ker(\beta^i)$ of subsheaves is eventually constant, so that the
  descending chain $\im(\beta^i)\subset M$ is also eventually constant,
  stabilizing at $N$.  Since $\beta|_N$ is an automorphism and $\beta$
  factors through $\alpha_M$, we see that $\alpha_N$ must be surjective,
  and thus $\rank(N)=0$.  But $M$ is torsion-free, and thus $N=0$, so that
  $\beta$ is nilpotent as required.
\end{proof}

In general, we can write the sheaf $M|_Q$ as a direct sum $M|_Q =
\bigoplus_{\mu\in \Q\cup\{\infty\}} V_\mu$ where each $V_\mu$ is a
semistable sheaf of slope $\mu$.  By the Corollary, when $M$ is
rigid and torsion-free, this direct sum decomposition lifts
(nonuniquely!) to $M$ and thus $M$ itself is a sum of torsion-free
sheaves $M_\mu$ such that $(M_\mu)|_Q=V_\mu$.  (Note that
$\rank(M_\infty)=\rank(V_\infty)=0$, and thus $M_\infty=0$.)  And,
of course, each $M_\mu$ itself splits as a direct sum of sheaves
with {\em indecomposable} restriction to $Q$.

The key observation of \cite{KuleshovSA:1997} is that such a summand
$M_\mu$ is (slope) semistable in a suitable sense.

\begin{prop}
  Let $M$ be a torsion-free sheaf on $X$ such that $M|_Q$ is semistable.
  Then any nontrivial subsheaf $N\subset M$ satisfies
  \[
  \frac{c_1(N)\cdot Q}{\rank(N)}\le \frac{c_1(M)\cdot Q}{\rank(M)}.
  \]
  If $M|_Q$ is stable, equality implies $N=M$.
\end{prop}

\begin{proof}
  We may assume that $M/N$ is torsion-free, as otherwise replacing $N$ by
  the saturated subsheaf can only increase $c_1(N)\cdot Q$ without
  increasing $\rank(N)$.  (Here we use the fact that $Q$ is nef!)  But then
  (since $(M/N)|_Q$ is a sheaf) we have a short exact sequence
  \[
  0\to N|_Q\to M|_Q\to (M/N)|_Q\to 0
  \]
  and thus by assumption have $\mu(N|_Q)\le \mu(M|_Q)$.  The first claim then
  follows by observing that $\rank(N|_Q)=\rank(N)$ and
  $\deg(N|_Q)=c_1(N)\cdot Q$.  Moreover, stability of $M|_Q$ implies that
  $\mu(M|_Q)$ has denominator $\rank(M|_Q)$, and thus $\mu(N|_Q)=\mu(M|_Q)$
  only if $\rank(N)=\rank(M)$, only if $N=M$.
\end{proof}

Define $\mu_Q(M):=(c_1(M)\cdot Q)/\rank(M)=\mu(M|_Q)$, and say that the
torsion-free sheaf $M$ is $\mu_Q$-semistable if $\mu_Q(N)\le \mu_Q(M)$ for
nonzero $N\subset M$ (and similarly for $\mu_Q$-stable).  In particular, we
have just shown that (semi)stability of $M|_Q$ implies
$\mu_Q$-(semi)stability of $M$.

\begin{cor}
  If $M$, $N$ are $\mu_Q$-semistable torsion-free sheaves, then
  $\Hom(M,N)=0$ if $\mu_Q(M)>\mu_Q(N)$ and $\Ext^2(M,N)=0$ if
  $\mu_Q(M)<\mu_Q(N)+Q^2$ (and in particular if $M=N$).
\end{cor}

\begin{proof}
  If $\Hom(M,N)\ne 0$, then the image $I$ of such a homomorphism satisfies
  \[
  \mu_Q(M)\le \mu_Q(I)\le \mu_Q(N),
  \]
  since it is a quotient of $M$ and a subsheaf of $N$.  By duality, if
  $\Ext^2(M,N)\ne 0$, then $\Hom(N,\theta M)\ne 0$, and thus $\mu_Q(N)\le
  \mu_Q(\theta M)=\mu_Q(M)-Q^2$.
\end{proof}

Recall that an exceptional collection in a triangulated category is a
sequence $(E_1,\dots,E_n)$ of objects such that $R\Hom(E_i,E_i)\cong k$ and
$R\Hom(E_i,E_j)=0$ for $i>j$.  In the commutative case, the following is a
special case of \cite[Lem. 2.4.4]{KuleshovSA:1997}.  (The proof below is
somewhat more elementary, as we do not need to use the spectral sequence
for $\Ext$ of filtered sheaves.)

\begin{prop}
  Let $M$ be a rigid torsion-free sheaf with $M|_Q$ indecomposable.  Then
  there is an induced filtration $0=F_n\subset F_{n-1}\subset\cdots\subset
  F_0=M$ with subquotients $E_i=F_{i-1}/F_i$ such that $(E_1,\dots,E_n)$ is
  an exceptional collection.
\end{prop}

\begin{proof}
  By assumption, $\Ext^1(M,M)=0$, and slope considerations imply
  $\Ext^2(M,M)=0$, so that $\End(M)\cong \End_Q(M|_Q)\cong k[t]/(t^n)$ for
  some $n$.  We may thus induce a filtration on $M$ by taking
  $F_i=\ker(t^{n-i})$, and see that the restriction of this filtration to
  $Q$ has the same description.  In particular,
  $\rank(F_i)=(1-i/n)\rank(M)$, and thus each $E_i|_Q$ is isomorphic to the
  unique stable constituent $V$ of $M|_Q$.  The action of $t$ certainly
  respects this filtration, and thus has an induced action on the
  associated graded.  Since $t^{n-1-i}$ induces an injection $E_i\to E_n$,
  we see that $t:E_i\to E_{i+1}$ is injective for all $i$, giving a chain
  $E_1\subset \cdots \subset E_n$ of sheaves restricting to the same stable
  bundle $V$ on $Q$.

  Slope considerations again tell us that $\Ext^2(E_i,E_j)=0$ for all
  $i,j$, and thus $\Hom(E_i,E_j)\to \End(V)\cong k$ is injective for all
  $i,j$.  Moreover, the natural map $\Ext^1(E_i,E_{i+1})\to
  \Ext^1(V,V)\cong k$ is nonzero (since the corresponding subquotient of
  $M|_Q$ is indecomposable), and thus surjective, implying that
  $\Ext^2(E_i,\theta E_{i+1})\to \Ext^2(E_i,E_{i+1})=0$ is injective.  We
  thus conclude that $\Hom(E_{i+1},E_i)=0$ so that $\Hom(E_i,E_j)=0$ for
  all $i>j$, and thus (reversing the above argument) $\Ext^1(E_i,E_j)\to
  \Ext^1(V,V)$ is surjective for $i<j$.

  We thus conclude that the $\Ext$ spaces satisfy the inequalities
  \begin{align}
  \dim\Hom(E_i,E_j)&=\delta_{i\le j}\notag\\
  \dim\Ext^1(E_i,E_j)&\ge \delta_{i<j}\notag\\
  \dim\Ext^2(E_i,E_j)&=0\notag
  \end{align}
  and thus
  \[
  \chi(E_i,E_j)\le \delta_{i\le j}-\delta_{i<j}=\delta_{ij}.
  \]
  Since
  \[
  \chi(M,M)=\sum_{i,j}\chi(E_i,E_j)\le n=\chi(M,M),
  \]
  the inequalities must be tight, and thus in particular $(E_1,\dots,E_n)$
  is an exceptional collection as required.
\end{proof}

\begin{prop}
  If $M$, $N$ are sheaves with $\mu_Q(M)=\mu_Q(N)\in \Q\cup\infty$, then
  $\chi(M,N)=\chi(N,M)$. In particular, if $E$ and $F$ are exceptional
  sheaves of the same slope, then $\chi(E,F)=\chi(F,E)=1+\rho^2/2\le 1$,
  where $\rho=c_1(E)-c_1(F)$.
\end{prop}

\begin{proof}
  From the distinguished triangle
  \[
  R\Hom(M,\theta N)\to R\Hom(M,N)\to R\Hom_Q(M|_Q,N|_Q)\to,
  \]
  we obtain the identity
  \[
  \chi_Q(M|_Q,N|_Q) = \chi(M,N)-\chi(M,\theta N).
  \]
  Serre duality gives $\chi(M,\theta N)=\chi(N,M)$, and the fact that
  $\chi_Q$ is a symplectic pairing on $K_0(Q)$ implies
  $\chi_Q(M|_Q,N|_Q)=0$ and thus $\chi(M,N)=\chi(N,M)$.

  In the exceptional case, we have
  \[
  -\rho^2 = \chi([M]-[N],[M]-[N]) = \chi(M,M)+\chi(N,N)-2\chi(M,N) =
  2-2\chi(M,N)
  \]
  and solving for $\chi(M,N)$ gives the desired result.  The inequality
  then follows from the fact that $\rho$ is in the negative definite even
  lattice of divisors orthogonal to $Q$.
\end{proof}

\begin{lem}\label{lem:exceptionals_comparable}
  Let $M_1\oplus M_2$ be a rigid torsion-free sheaf such that $M_1|_Q$,
  $M_2|_Q$ are indecomposable bundles with isomorphic stable constituents.
  If $\rank(M_1)\ge \rank(M_2)$, then (up to isomorphism) the exceptional
  collection associated to $M_1$ contains that associated to $M_2$.
\end{lem}

\begin{proof}
  As usual, slope considerations imply that $\Ext^2(M_i,M_j)=0$.  Let the
  two exceptional collections be $(E_1,\dots,E_{n_1})$ and
  $(E'_1,\dots,E'_{n_2})$, with $n_1\ge n_2$.  Then $\chi(E_i,E'_j)\le
  \delta_{E_i\cong E'_j}$, and thus
  \[
  n_2 = \dim\Hom(M_1,M_2) = \chi(M_1,M_2) = \sum_{i,j}\chi(E_i,E'_j)\le 
  \sum_{i,j} \delta_{E_i\cong E'_j}\le n_2,
  \]
  since each $E'_j$ is isomorphic to at most one $E_i$.  It follows that
  both inequalities are tight, and in particular that every $E'_j$ is
  isomorphic to some $E_i$ as required.
\end{proof}

\begin{lem}
  Let $M_1\oplus M_2$ be a rigid torsion-free sheaf such that $M_1|_Q$ and
  $M_2|_Q$ are indecomposable bundles with $\mu(M_1|_Q)\le
  \mu(M_2|_Q)<\mu(M_1|_Q)+Q^2$.  Then the exceptional collections
  associated to $M_1$ and $M_2$ are contained in a common exceptional
  collection.
\end{lem}

\begin{proof}
  If $M_1|_Q$, $M_2|_Q$ have isomorphic stable constituents, then this
  follows from the previous Lemma.  Otherwise, let the two exceptional
  collections again be $(E_1,\dots,E_{n_1})$ and $(E'_1,\dots,E'_{n_2})$.
  Then we claim that $(E_1,\dots,E_{n_1},E'_1,\dots,E'_{n_2})$ is an
  exceptional collection.
  
  The slope constraint implies that $\Hom(M_2,M_1)=\Ext^2(M_2,M_1)=0$ and
  (since they have the same slopes) $\Hom(E'_i,E_j)=\Ext^2(E'_i,E_j)=0$ for
  all $i,j$.  By hypothesis, $\Ext^1(M_2,M_1)=0$, so that
  $\chi(M_2,M_1)=0$.  We then find that
  \[
  0 = \chi(M_2,M_1)=\sum_{i,j} \chi(E'_i,E_j) \le 0
  \]
  and thus $\chi(E'_i,E_j)=0$ for all $i,j$, implying $\Ext^1(E'_i,E_j)=0$
  as required.
\end{proof}

\begin{cor}
  Let $M$ be a rigid torsion-free sheaf with slopes contained in a
  half-open interval of length $Q^2$.  Then there is an exceptional
  collection $(E_1,\dots,E_n)$ and a filtration of $M$ such that the $i$th
  subquotient has the form $n_iE_i$.
\end{cor}

\begin{proof}
  The two Lemmas tell us that any two indecomposable summands of $M$ have
  compatible exceptional collections, and thus we immediately conclude that
  there is a single exceptional collection containing all of them;
  moreover, that exceptional collection can be split into blocks associated
  to the different stable constitutents of $M|_Q$, and the existence of a
  filtration as described for each such block is straightforward.
\end{proof}

\begin{rem}
  Since the hypothesis ensures that $\Ext^2(M,M)=0$, this is (in the
  commutative case) a special case of \cite[Thm. 2.5.1(2)]{KuleshovSA:1997}.  The
  argument there can be extended to the noncommutative case, but we will
  not need this generalization, so omit the details.
\end{rem}

Of course, the above is in a sense the converse of the result we want: we
want to know when a given exceptional collection corresponds to a rigid,
unobstructed torsion-free sheaf.  In the indecomposable case, this is
straightforward.

\begin{prop}\label{prop:construct_fitered_indec}
  Suppose $(E_1,\dots,E_n)$ is an exceptional collection on a
  noncommutative del Pezzo surface, such that each $E_i$ is a torsion-free
  sheaf and the restrictions to $Q$ are all isomorphic.  Then
  $\dim\Hom(E_i,E_j)=\dim\Ext^1(E_i,E_j)=1$ for $i<j$, and there is (up to
  isomorphism) a unique filtered sheaf $0=F_n\subset \cdots \subset F_0=M$
  with $F_{i-1}/F_i\cong E_i$, $\Ext^1(M,M)=0$, and $M|_Q$ indecomposable.
\end{prop}

\begin{proof}
  For $i<j$, since $R\Hom(E_j,E_i)=0$, we have $R\Hom(E_i,E_j)\cong
  R\Hom(E_i|_Q,E_j|_Q)$, and the sheaves on $Q$ are isomorphic stable
  sheaves.  We then define a sequence $F_i$ by letting $F_{i-1}$ be the
  generic extension of $E_i$ by $F_i$, giving a sheaf $M$.  That this
  satisfies $\Ext^1(F_i,F_i)=0$ for all $i$ follows as in
  \cite{KuleshovSA:1997}, giving existence.  Moreover, the only way to
  obtain a nonisomorphic sheaf with the same subquotients is to take a
  nongeneric extension in some step, but since $\Ext^1(E_i,F_i)\cong
  \Ext^1(E_i|_Q,F_i|_Q)\cong k$ at the first such step, this would make
  $F_0|_Q$ decomposable.
\end{proof}

\begin{thm}
  Let $(E_1,\dots,E_n)$ be an exceptional collection of torsion-free
  sheaves, with $\mu_Q(E_i)$ constant.  Then for any sequence
  $(m_1,\dots,m_n)$ of nonnegative integers, there is a unique rigid
  sheaf $M$ having a filtration with subquotients $E_i^{m_i}$.
\end{thm}

\begin{proof}
  For existence, note that for any subset $S\subset \{1,\dots,n\}$, there
  is a corresponding sheaf $M_S$ obtained as in Proposition
  \ref{prop:construct_fitered_indec}: at each step, take the generic extension.
  We then find as above that
  \[
  \chi(M_S,M_T) = |S\cap T|.
  \]
  Since $\Ext^2(M_T,M_S)=0$, we find that $\dim\Hom(M_S,M_T)\le
  \dim\Hom(M_S|_Q,M_T|_Q)$.  If $S\subset T$ or $T\subset S$, then
  $\dim(M_S|_Q,M_T|_Q)=|S\cap T|$, and
  \[
  |S\cap T|=\chi(M_S,M_T)\le \dim\Hom(M_S,M_T)\le |S\cap T|,
  \]
  forcing $\dim\Ext^1(M_S,M_T)=0$ and $\dim\Hom(M_S,M_T)=|S\cap T|$.
  Existence then follows by taking the sum $\sum_i M_{S_i}^{c_i}$ for a
  suitable chain of subsets and multiplicities $c_i$.  (The chain of
  subsets is determined by a permutation that sorts the sequence $\vec{m}$,
  and the $c_i$ are then differences of consecutive entries in the sorted
  sequence.  This permutation need not be unique, but the sequence $c_i$ is
  independent of the choice of permutation, and when $c_i\ne 0$, $S_i$ is
  also uniquely determined.)

  For uniqueness, we may use the fact that $M$ decomposes as a direct sum
  over different isomorphism classes of stable constituents to reduce to
  the case that $E_1|_Q\cong\cdots\cong E_n|_Q$, at which point we conclude
  from Lemma \ref{lem:exceptionals_comparable} that if both $M_S$ and $M_T$
  appear, then $S$ and $T$ are comparable.
\end{proof}

\begin{rem}
  There is a slight additional technical issue if we want to perform this
  construction in families, namely the fact that over a field, we are
  taking {\em generic} extensions at various points in the process, and
  thus run into potential problems in constructing a global extension class
  that is nonzero whenever possible.  Note that since the relevant $\Ext^1$
  spaces are preserved by restriction to $Q$, it suffices to show that one
  can define $M_S|_Q$ globally for each $S$.  Since this sheaf is
  semistable, applying a suitable derived autoequivalence of
  $D^b_{\coh}(E)$ turns it into a torsion sheaf (\cite{AtiyahMF:1957}; see
  Proposition \ref{prop:atiyah} above), the fibers of which are isomorphic
  to structures sheaves of divisors.  But such a torsion sheaf can
  certainly be defined globally: simply view the family of divisors as a
  divisor on the total space of the family and take its structure sheaf.
  Inverting the derived autoequivalence gives a global version of $M_S|_Q$
  as required, and thus the required global extension classes needed to
  construct $M_S$ itself.
\end{rem}

More generally, if we are given an exceptional collection consisting of
torsion-free sheaves with nondecreasing $Q$-slopes contained in a half-open
interval of length $Q^2$, and for each $\mu$ we construct a rigid sheaf
$M_\mu$ from the corresponding subcollection, then $M=\bigoplus_\mu M_\mu$
is rigid and unobstructed.  Indeed, we have $R\Hom(M_\nu,M_\mu)=0$ for
$\mu<\nu$, and thus $R\Hom(M_\mu,M_\nu)\cong
R\Hom_Q(M_\mu|_Q,M_\nu|_Q)\cong \Hom(M_\mu|_Q,M_\nu|_Q)$, so that in
neither direction are there higher morphisms.

Thus the question of classifying filtered deformations arising from
Proposition \ref{prop:fildef_from_sheaf} reduces to one of classifying
exceptional collections.  The basic question here is: given a sequence of
classes in $K_0(X)$ satisfying the obvious numerical conditions, is there a
corresponding exceptional collection, and is it unique?  Unfortunately, the
main result of \cite{KuleshovSA:1997} on exceptional collections is not
particularly useful in this regard: it says that any exceptional collection
can be obtained from a certain standard collection via the braid group
action, but it is unclear how to turn that into an existence test, and
(worse yet) the arguments depend in a crucial way on the commutativity of
the surface.  (They involve restriction functors to $D^b_{\coh}(\P^1)$
coming from $-1$-curves.)

Although we do not have a complete answer, we can reduce the problem
considerably.

\begin{prop}
  Suppose $(E_1,\dots,E_n)$ is a sequence of exceptional torsion-free
  sheaves such that $\mu_Q(E_1)\le \mu_Q(E_2)\le \cdots\le \mu_Q(E_n)<
  \mu_Q(E_1)+Q^2$ and $\chi(E_j,E_i)=0$ for $i<j$; suppose moreover that if
  $\mu_Q(E_i)=\mu_Q(E_j)$ for $i<j$, then there is no sheaf of class
  $[E_i]-[E_j]$.  Then $(E_1,\dots,E_n)$ is an exceptional collection.
\end{prop}

\begin{proof}
  Let $i<j$.  Since $\mu_Q(E_j)<\mu_Q(E_i)+Q^2$, we find
  $\Ext^2(E_j,E_i)=0$.  Since $\chi(E_j,E_i)=0$, it remains only to show
  that $\Hom(E_j,E_i)=0$.  If $\mu_Q(E_i)<\mu_Q(E_j)$, this follows from
  stability, while if the slopes agree, any morphism must be injective, and
  thus the cokernel is a sheaf of the given class in $K_0(X)$.
\end{proof}

\begin{rem}
  Note that as observed above, when the slopes agree, $c_1(E_i)-c_1(E_j)$
  is a root of the lattice $Q^\perp\sim E_{9-Q^2}$.  It follows that there
  is a sheaf of class $[E_i]-[E_j]$ iff $c_1(E_i)-c_1(E_j)$ is a {\em
    positive} root and $[E_i]|_Q\cong [E_j]|_Q$.  In particular, if we
  arrange for $c_1(E_j)-c_1(E_i)$ to be a positive root for $i<j$, then
  $[E_i]-[E_j]$ can never be effective.
\end{rem}

Uniqueness is even easier.

\begin{prop}
  If $E$, $F$ are exceptional torsion-free sheaves with $[E]=[F]\in
  K_0(X)$, then $E\cong F$.
\end{prop}

\begin{proof}
  Since $\chi(E,E)=\chi(F,F)=1$, we have $\chi(E,F)=\chi(F,E)=1$.  Since
  $E$ is rigid and unobstructed, $k\cong \End(E)\cong \End_Q(E|_Q)$, and
  thus $E|_Q$ is stable, so that $E$ is $\mu_Q$-stable.  Slope considerations
  imply that $\Ext^2(E,F)=\Ext^2(F,E)=0$ and thus there are nonzero
  morphisms in either direction.  Stability forces the morphisms to be
  injective, and thus so is the composition in either direction, which must
  therefore be a nonzero multiple of the identity, implying $E\cong F$ as
  required.
\end{proof}

Thus to understand the exceptional collections corresponding to filtered
deformations, it suffices to understand (a) which sequences of classes in
$K_0(X)$ satisfy the numerical conditions (and positivity of roots, when
relevant) and (b) when a given class in $K_0(X)$ is represented by an
exceptional torsion-free sheaf.

\begin{lem}
  Let $M$ be a $\mu_Q$-stable torsion-free sheaf, and let $D\in \NS(X)$ be
  an ample divisor class.  Then for all sufficiently large $n$, $M$ is
  $nQ+D$-stable.
\end{lem}

\begin{proof}
  Suppose otherwise, and for each $n\ge 0$, let $B_n$ be the maximally
  destabilizing subsheaf of $M$.  Since $M$ is $\mu_Q$-stable, we must have
  $\mu_Q(B_n)<\mu_Q(M)$, and thus there is some $n'>n$ such that
  \[
  n'\mu_Q(B_n)+\mu_D(B_n)< n'\mu_Q(M)+\mu_D(M),
  \]
  and thus in particular
  \[
  n'\mu_Q(B_{n'})+\mu_D(B_{n'})>n'\mu_Q(B_n)+\mu_D(B_n)
  \]
  Since $B_n$ is maximally destabilizing for $nQ+D$, we also have
  \[
  n\mu_Q(B_{n'})+\mu_D(B_{n'})\le n\mu_Q(B_n)+\mu_D(B_n),
  \]
  and thus $\mu_Q(B_n)<\mu_Q(B_{n'})$.  We thus obtain an infinite sequence
  of strictly increasing slopes bounded by $\mu_Q(M)$, but since the
  denominators are bounded by $\rank(M)$, this is impossible.
\end{proof}

\begin{rem}
  If $M$ is only assumed to be $\mu_Q$-semistable, an analogous argument
  shows that either $M$ is still $nQ+D$-stable for $n\gg 0$ or $B_n$
  eventually stabilizes, and thus by induction the corresponding
  Harder-Narasimhan filtration of $M$ eventually stabilizes.
\end{rem}

\begin{prop}\label{prop:all_del_Pezzos_have_same_exceptionals}
  Let $R$ be a discrete valuation ring with residue field $k$ and field of
  fractions $K$, and let $X_R$ be a noncommutative del Pezzo surface over
  $R$.  An exceptional torsion-free sheaf on either fiber of $X_R$ extends
  to a torsion-free sheaf which is exceptional on the other fiber.
\end{prop}

\begin{proof}
  An exceptional sheaf $E_k$ on $X_k$ is simple, so induces a point of the
  moduli space of simple sheaves.  (It is shown in
  \cite[Cor. 11.3]{noncomm2} that this moduli problem is represented by an
  algebraic space, with Prop. 11.1 op.~cit. implying that the obstruction
  theory is as expected.)  Since $\Ext^2(E_k,E_k)=0$, that algebraic space
  is formally smooth (and thus smooth) at $E_k$, so that the point
  corresponding to $E_k$ extends to $R$.  Since $\gcd(c_1(E)\cdot
  Q,\rank(E))=\gcd(\deg(E|_Q),\rank(E|_Q))=1$, there is a universal family.
  (The moduli stack of simple sheaves is a $\G_m$-gerbe over the moduli
  space of simple sheaves, so that the obstruction to the existence of a
  universal family is a class in the Brauer group.  For any sufficiently
  ample $D$, this class is represented by the Azumaya algebra
  $\End(\Hom(\sO_X(-D),E))$, and the $\gcd$ condition ensures that the
  degrees of those Azumaya algebras generate $\Z$, so the obstruction is
  trivial.)  Thus this point corresponds to a simple sheaf $E_R$, which has
  exceptional generic fiber by semicontinuity of $\Ext^i$.

  For the other direction, fix an ample divisor $D$ on $X_k$, so that an
  exceptional torsion-free sheaf $E_K$ on $X_K$ is $nQ+rD$-stable for some
  $n$, with $r=\rank(E_K)$.  Since the semistable moduli functor is
  universally closed, there is thus an extension to a sheaf $E_R$ with
  $nQ+rD$-semistable special fiber $E_k$.  If $\gcd(n,r)=1$, then the
  corresponding slope of $E_K$ has denominator $r$, and thus cannot be the
  slope of a sheaf of smaller rank, so that $E_k$ is actually {\em
    stable}. It follows that $\End(E_k)\cong k$, and thus that
  $\Hom(E_k,\theta E_k)=0$, giving $\Ext^2(E_k,E_k)=0$.  Since
  $\chi(E_k,E_k)=\chi(E_K,E_K)=1$, it follows that $E_k$ is exceptional as
  required.
\end{proof}

It follows immediately that, given a class $[E]\in K_0(X)$ such that
$\rank([E])>0$ and $\chi([E],[E])=1$, if {\em any} noncommutative del Pezzo
surface (with the correct value of $Q^2$ and parity if $Q^2=8$; i.e., in a
given component of the moduli stack of noncommutative surfaces) has an
exceptional torsion-free sheaf of class $[E]$, then {\em every}
noncommutative del Pezzo surface has such a sheaf.  (Indeed, every $X$ can
be obtained from a finite base change of the generic surface by a
descending down a sequence of valuations of the corresponding field, and
the existence of a universal family lets us descend through finite field
extensions.)  We thus in particular see that the existence question reduces
to the case of commutative del Pezzo surfaces, which may even be assumed
nondegenerate.

In fact, it reduces to the case of (commutative, nondegenerate) del Pezzo
surfaces of degree 1.  Let $\pi:Y\to X$ be a birational morphism with $Y$ a
nondegenerate commutative del Pezzo surface of degree $1$.  If $[E]\in
K_0(X)$ is the class of an exceptional bundle $E$, then $\pi^*[E]$ is the
class of the exceptional bundle $\pi^*E$.  Conversely, if $\pi^*[E]$ is the
class of an exceptional bundle $E'$, then it descends to an exceptional
bundle on $X$ by \cite[Cor 3.2]{KuleshovSA/OrlovDO:1994}.

A fairly strong necessary condition for $[E]$ to be representable comes
from the following.

\begin{lem}
  If $E$ is an exceptional bundle on $X$, then for any exceptional bundle $E'$
  with $\mu(E')-Q^2<\mu(E)<\mu(E')$, $\chi(E',E)\le 0$.
\end{lem}

\begin{proof}
  Indeed, the conditions on the slopes ensure that
  $\Hom(E',E)=\Ext^2(E',E)=0$ and thus $\chi(E',E)=-\dim\Ext^1(E',E)\le 0$.
\end{proof}

Since the twist of an exceptional bundle by a line bundle remains
exceptional, any given $E'$ induces a whole family of necessary conditions:
$\chi(E'(-D),[E])\le 0$ for all $D\in Q^\perp=\{v\in \NS(X)|v\cdot Q=0\}$.
Since (letting $r=\rank([E])$, $r'=\rank(E')$)
\[
\chi(E'(-D),[E])
-
\frac{rr'}{2}
(D - (D_{E'}/r'-D_E/r))^2
\]
is independent of $D$, with $D_E:=c_1(E)-\frac{c_1(E)\cdot Q}{Q^2}Q$, we
see that $\chi(E'(-D),[E])$ is maximized when $D$ is the closest vector
(relative to the negative of the intersection form) in $Q^\perp$ to
$D_{E'}/r'-D_E/r$; if the maximum is positive, then $[E]$ is not
representable.  The resulting condition is particularly straightforward to
compute when $E'$ ranges over line bundles, and is already quite strong in
that case.

It turns out it is useful to reformulate the question in terms of one on
the elliptic surface $\tilde{X}$ associated to the nondegenerate degree 1
del Pezzo surface $X$.  (Note that nondegeneracy of $X$ corresponds to
smoothness of the Weierstrass model of $\tilde{X}$.)  Here, the
representability condition turns out to be essentially trivial.

\begin{prop}
  Let $\tilde{X}$ be a smooth rational Weierstrass elliptic surface, and
  let $[E]\in K_0(\tilde{X})$ be a class with $\rank([E])>0$,
  $\chi([E],[E])=1$, and $\gcd(\rank([E]),c_1([E])\cdot Q)=1$ (where $Q$
  denotes the divisor class of a fiber of $\tilde{X}$). Then $[E]$ is the
  class of an exceptional vector bundle on $\tilde{X}$, and any two bundles
  of class $[E]$ are isomorphic.
\end{prop}

\begin{proof}
  The derived autoequivalences of elliptic curves extend in a
  straightforward way to any Weierstrass curve, and thus in particular we
  can apply them to $\tilde{X}$, viewed as a family of such curves.  Those
  derived autoequivalences act as $\SL_2(\Z)$ on $(\rank(M),c_1(M)\cdot
  Q)$, and thus in particular there is a derived autoequivalence taking
  $[E]$ to a class with $\rank([E'])=0$, $c_1([E'])\cdot Q=1$,
  $\chi([E'],[E'])=1$.  Any such class is represented by a sheaf
  $\sO_e(\chi([E'])-1)$ where $e=c_1([E'])$ is a section of the elliptic
  fibration, and thus we can apply the inverse autoequivalence to obtain an
  object $E$ representing $[E]$.  Since a section $e$ cannot meet the
  singular point of any fiber, we find that the restriction of $E$ to any
  fiber is a shift of a vector bundle, and thus $E$ itself is a shift of a
  vector bundle.  This establishes existence, with uniqueness similarly
  following from the rank 0 case.
\end{proof}

\begin{prop}
  Let $\pi:\tilde{X}\to X$ be the birational morphism contracting the
  identity section of the smooth rational Weierstrass elliptic surface
  $\tilde{X}$ with marked smooth fiber $Q$, and let $[E]\in K_0(X)$ be a
  class with $\rank([E])>0$, $\chi([E],[E])=1$, and
  $\gcd(\rank([E]),c_1([E])\cdot Q)=1$.  Let $E'$ be the exceptional bundle
  on $\tilde{X}$ representing $\pi^*[E]$.  Then $[E]$ is represented by an
  exceptional bundle iff $R\Hom(E',\sO_e(-1))=0$, where $e$ is the identity
  section.
\end{prop}

\begin{proof}
  If $[E]$ is representable, then we must have $E'\cong \pi^*E$, and thus
  $R\Hom(E',\sO_e(-1))=0$ as required.  Conversely, if
  $R\Hom(E',\sO_e(-1))=0$, then $E'\cong \pi^*E$ for some sheaf $E$ of
  class $[E]$, necessarily exceptional.
\end{proof}

\begin{rem}
  This reformulation appears particularly amenable to a computational
  approach: rather than asking about whether an object exists, it asks
  about the properties of an explicit object $E'$.  In addition, the
  condition $R\Hom(E',\sO_e(-1))=0$ can be reformulated as $E'|_e\cong
  \sO_e(-1)^r$, so one need not even compute Ext spaces.  In addition, it
  suffices to consider {\em any} rational elliptic surface, and the
  reduction of Proposition \ref{prop:all_del_Pezzos_have_same_exceptionals}
  works equally well when the discrete valuation ring $R$ has mixed
  characteristic, so that one may feel free to compute in finite
  characteristic.  It is particularly tempting to work with the surface
  $y^2=x^3+z^6+w^6$ over a suitable finite field.  Not only is it easy to
  find prime fields over which this surface has rational Mordell-Weil group
  (any prime that splits completely in $\Q(2^{1/3})$ will do), but the
  Mordell-Weil group is the group of $C_6$-equivariant maps from the Fermat
  curve $t^6+u^6+v^6=0$ to the curve $E_0:y^2=x^3-1$ of $j$-invariant $0$,
  suggesting that it should be possible to reduce to a question about
  equivariant bundles on the Fermat curve.
\end{rem}

More generally, if $E$, $E'$ are two exceptional bundles on $\tilde{X}$
satisfying $R\Hom(E,E')=0$, then we can apply the standard derived
autoequivalences to put $E'$ in the form $\sO_{e'}(\chi-1)$, at which point
translating by the corresponding section and twisting by a line bundle
makes it $\sO_e(-1)$, and thus the image of $E$ under the same operations
gives rise to an exceptional bundle on $X$.  In fact, this works equally
well if $E$, $E'$ are bundles on some blowdown of $\tilde{X}$, since their
pullbacks to $\tilde{X}$ still form an exceptional pair.  This gives a
powerful technique for constructing exceptional bundles on del Pezzo
surfaces.

As an example, consider the case of two line bundles.  We may as well take
one of the bundles to be $\sO_X$, and thus ask when
$R\Hom(\sO_X,\sO_X(-D))=0$, with $D\cdot Q>0$.  In particular, we must have
$0=\chi(\sO_X(-D))=(D\cdot (D-Q))/2+1$, and since
$\chi(\sO_X(-D),\sO_X)>0$, we must have $\Hom(\sO_X(-D),\sO_X)\ne 0$, so
that $D$ is effective.  In other words, if $(\sO_X(-D),\sO_X)$ forms an
exceptional pair, then $D$ must be the class of a rational curve on $X$.
Conversely, if $D$ is the class of a rational curve, then $H^p(\sO_D)\cong
H^p(\sO_X)$, and thus $R\Gamma(\sO_X(-D))=0$ as required.

There is a {\em contravariant} derived autoequivalence $\Phi^*$ of
$\tilde{X}$ that swaps $\sO_{\tilde{X}}$ and $\sO_e(-1)$, and thus
preserves the triangulated subcategory of objects $M$ on $X$ with
$R\Gamma(M)=0$.  This is the composition of the autoequivalence $\Psi$ of
\cite[Lem. 12.5]{noncomm1} with Cohen-Macaulay duality, and thus we may
explicitly compute
\begin{align}
\rank(\Phi^* M) &= c_1(M)\cdot Q\\
c_1(\Phi^* M) &= -c_1(M)+(c_1(M)\cdot Q+\rank(M))(Q+e)+(c_1(M)\cdot e-\chi(M))Q\\
\chi(\Phi^* M) &= -c_1(M)\cdot e.
\end{align}
On the subcategory $\sO_X^\perp$ of $D^b_{\coh}(X)$, the action simplifies
to
\begin{align}
  \rank(\Phi^*M) &= c_1(M)\cdot Q  \\
  c_1(\Phi^*M) &= -c_1(M) + (c_1(M)\cdot Q+\rank(M))Q,\\
  \chi(\Phi^*M) &= 0,
\end{align}
where we note that the $Q$ here is the one on $X$, corresponding to
$\tilde{Q}+e$ on $\tilde{X}$.

In particular, since $\sO_X(-D)$ is an exceptional object on $X$ with
$R\Gamma(\sO_X(-D))=0$, applying $\Phi^*$ and shifting gives another such
object, with invariants
\[
\rank(E)=D\cdot Q,\quad
c_1(E)=-D+((D\cdot Q)-1)Q,\quad
\chi(E)=0.
\]
In this way, we obtain an exceptional bundle of slope $-1/r$ for all $r\ge
2$, and in two different ways for $r=6$ (corresponding to the two types of
rational curve of degree $6$: bidegree $(1,2)$ curves in $\P^1\times \P^1$
and conics in $\P^2$).

Conversely, suppose we are given an exceptional bundle $E$ of slope $-1/r$,
with Chern class of the form $-Q+v$ with $v\in \Lambda_{E_8}$.  This has
Euler characteristic $(r^2-r+2-v^2)/2r$ (taking the positive inner product
on $\Lambda_{E_8}$ rather than the intersection form), which we insist must
be negative for any representative of the coset $v+r\Lambda_{E_8}$.  We may
thus just as well assume that $v/r$ is in the Voronoi cell around $0$ (the
closure of the set of points closer to $0$ than to any other lattice
vector), so that $v$ is a minimal coset representative.  The maximum of
$w^2$ over the Voronoi cell is $1$, and thus $v^2\le r^2$, forcing
$v^2=r^2-r+2$ and $\chi=0$.  We thus obtain an exceptional pair $(E,\sO_X)$
on $X$, and thus an exceptional triple $(\sO_e(-1),\pi^*E,\sO_{\tilde{X}})$
on $\tilde{X}$.  Applying $\Phi^*$ turns $\pi^*E$ into (a shift of) a line
bundle, and thus $E$ must have come from the above construction.

The above example is somewhat ad hoc, but there is a somewhat more
systematic approach available.  Consider the case of slope $-2/r$ ($r$
odd).  Here $v^2=r^2-2r+5$, so we still must have $\chi=0$, but it is
difficult to control when $R\Hom(\sO_X,V)=0$ for $V$ an exceptional bundle
of rank 2 and very small negative slope.  However, since $v^2>(r-1)^2$, so
that $v/(r-1)$ must be outside the Voronoi cell, we find that our putative
bundle $E$ fits into an exceptional triple on $E$ of the form
\[
(E,\sO_X(-\alpha),\sO_X)
\]
where $\alpha$ is the highest root.  Applying the Cohen-Macaulay dual gives
a triple
\[
(\omega_X,\omega_X(\alpha),E^D),
\]
inducing a triple
\[
(E^D,\sO_X,\sO_X(\alpha)),
\]
in which we can swap $\sO_X$ and $\sO_X(\alpha)$.  Lifting to $\tilde{X}$
and applying $\Phi^*$ then gives a four-term exceptional collection of the
form
\[
(\sO_e(-1),\sO_{e'}(-1),M,\sO_X)
\]
where $M$ has rank $r-2$ and slope $-r/(r-2)$ and (up to the action of
$W(E_8)$) is the pullback of a sheaf on a degree 2 del Pezzo surface.

We could, of course, do this calculation entirely in the Grothendieck group
to get a class $[M]$.  Clearly, if $[E]$ is representable, then $[M]$ is
representable, but it turns out that the converse applies as well.  Indeed,
if $[M]$ is representable, then since $\mu(M)\in (-2,0)$, we find that
$R\Hom(\sO_X,M)=0$ and thus we indeed obtain a four-term exceptional
collection, and can then invert the derived equivalences to deduce that
$[E]$ is representable.  Moreover, the representability of $[M]$ is
equivalent to that of the class $[M](Q_X)$ of slope $-2/(r-2)$, and may
thus be shown by induction.  (To make this explicit, we must reduce
$v/(r-2)$ to the fundamental alcove, but this requires only a small number
of affine reflections for $r$ large since it is already quite close to the
fundamental alcove!)  We use the standard coordinatization of the
fundamental alcove via inner products with the simple roots.  (In terms of
the basis of $\Pic(X)$ coming from $\P^2$, the simple roots are
$(h-e_1-e_2-e_3,e_1-e_2,e_2-e_3,\dots,e_7-e_8)$.)

\begin{prop}
  An exceptional sheaf of slope $-2/r$ on a degree 1 del Pezzo surface has
  Chern class $-2Q+v$ where the image of $v/r$ in the fundamental alcove is
  either
  \[
  r^{-1}(0,(r-5)/2,0,0,1,0,0,0,0)\notag
  \]
  or one of the following sporadic cases:
  \begin{align}
  3^{-1}(1,0,0,0,0,0,0,0,0)&, \notag\\
  9^{-1}(1,0,0,1,0,0,0,0,0)&, \notag\\
  11^{-1}(1,0,2,0,0,0,0,0,0)&,\notag\\
  13^{-1}(3,2,0,0,0,0,0,0,0)&,\notag\\
  15^{-1}(5,0,0,0,0,0,0,0,0)& \notag
  \end{align}
\end{prop}

\begin{proof}
  Since $v^2$ is close to $r^2$, $v/r$ must be close to the fundamental
  weight $(0,1/2,0,0,0,0,0,0,0)$, from which it is straightforward to
  verify that for $r\gg 0$, the only possibility is as stated.  Moreover,
  one finds that for $r\ge 5$, the image of $(r/(r-2))v$ in the fundamental
  alcove has the same form, while for $r=5$, the image is
  $3^{-1}(1,0,0,0,0,0,0,0,0)$.  The only other possibilities for $r\le 17$
  are as stated, with the sporadic instances for $r\in \{11,13,15\}$ each
  reducing to the sporadic instance of rank $r-2$, and the sporadic
  instance for $r=9$ reducing to the regular instance of rank $7$.  On the
  other hand, for the regular instance of rank $r\ge 9$, there is a unique
  orbit of $-1$-curves on $X$ (consisting of four orthogonal $-1$-curves)
  such that the slope $-(r-2)/r$ bundle has $\chi(E,\sO_e(-1))=0$, and thus
  (by duality) at most one $W(E_7)$-orbit of numerical classes of slope
  $-(r+2)/r$ corresponding to exceptional bundles on the degree $2$ del
  Pezzo surface.  We thus see that the regular instance can account for at
  most one type of bundle of rank $r+2$, from which the claim follows by
  induction.
\end{proof}  

\begin{rem}
  In both the slope $-1/r$ and $-2/r$ cases, one can verify that there are
  no other classes with $v^2=r^2-r+2$ or $v^2=r^2-2r+5$, and thus that the
  above necessary condition is sufficient for bundles of degree $\pm 1$,
  $\pm 2$ modulo $r$.
\end{rem}

The interpretation via $\tilde{X}$ also gives rise to one further
construction.  If $E$ is an exceptional bundle on $X$ of slope $-d/r\in
(0,1)$, then in addition to the pair $(\pi^*E,\sO_e(-1))$ on $\tilde{X}$,
we also obtain a pair $(\sO_e,\pi^*E)$, which itself gives rise to an
exceptional bundle $E'$ on $X$.  This new bundle has slope $-d'/r\in (0,1)$
where $d'd=-1(r)$ and satisfies $v' = d'v$.  As before, this construction
not only gives a new bundle but a new necessary condition.  In this way,
the slope $-2/r$ construction gives rise to a construction for slope
$(r-1)/2r$ ($r$ odd).  Together with Cohen-Macaulay duality, this gives us
a complete classification for slopes $-1/r$, $-2/r$, $(r-1)/2r$,
$(r+1)/2r$, $(r-2)/r$, $(r-1)/r$.  This accounts for every possibility for
$r\le 9$ except $\{-3/8,-5/8\}$.  For slope $-3/8$, there are two cases
with $\chi\le 0$, but for one of those cases the putative bundle $E'$ of
slope $-d'/r=-5/8$ violates the necessary condition.  We thus see that the
only possibility for slope $-3/8$ corresponds to $v$ with coordinates
\[
8^{-1}(0,0,1,0,0,1,0,0),
\]
and can verify that this reduces to the known bundle of slope $-3/5$.

For rank $\ge 10$, the above techniques no longer suffice to obtain a
classification: there is a candidate of slope $-3/10$ with maximal Euler
characteristic $-1$.  (A direct computation on a suitable Weierstrass
surface shows that no such bundle exists.)

\section{Classifications}

We now return to our original question: given a bundle $V$ and
autoequivalence $\Psi$, understand the filtered deformations of
$B_{V,\Psi}$.  If as usual we constrain the slopes of $V$ to an interval of
length $d$, then we can produce such a deformation from the following data:
(1) a noncommutative del Pezzo surface $X$ with anticanonical curve $Q$
such that $\theta\otimes H^1(\sO_Q)\cong \Psi^{-1}$ on $\sO_Q$, so in
particular $Q^2=d$, (2) an abelian equivalence $\psi:\coh(Q)\cong \coh(E)$,
and (3) a rigid, unobstructed, torsion-free sheaf $M$ satisfying $\psi
M|_Q\cong V$.  The discussion above further allows us to replace $M$ by an
exceptional collection $(E_1,\dots,E_n)$ such that the bundles $(\psi
E_1|_Q,\dots,\psi E_n|_Q)$ are the stable constituents of $V$ (with
appropriate multiplicities).

Now, a pair $(X,\psi)$ with $X$ a noncommutative del Pezzo surface and
$\psi:\coh(Q)\cong \coh(E)$ induces a homomorphism $K_0(X)\to K_0(E)$ by
$[M]\mapsto [\psi M|_Q]$, and one can in fact recover the surface from this
homomorphism.  Indeed, $\psi$ as usual has the form $\psi(M)=f_*(M)\otimes
\psi(\sO_Q)$ where $f:Q\cong E$ is an automorphism of curves.  Since
$\psi(\sO_Q)\cong\psi(\sO_X|_Q)$, we can recover $[\psi(\sO_Q)]$ (and thus,
up to isomorphism, $\psi(\sO_Q)$) from the homomorphism $K_0(X)\to K_0(E)$.
Since
\[
[f_*M|_Q]\cong [\psi M|_Q]-\rank(M|_Q)[{\cal L}],
\]
we can thus compute the homomorphism $[M]\mapsto [f_*M|_Q]$ from the given
homomorphism.  But this homomorphism is precisely of the kind we considered
in Section 4, so lets us recover $X$.

This works for any homomorphism $\phi:K_0(X)\to K_0(E)$, as long as
$\rank(\phi([M]))=\rank([M])$ for all $[M]\in K_0(X)$, and
$\deg(\phi([M]))=c_1([M])\cdot Q+\delta \rank([M])$ for some $\delta\in
\Z$.  (It then follows that
$\phi'([M]):=\phi([M])-\rank([M])(\phi([\sO_X])-[\sO_Q])$ gives rise to a
noncommutative surface $(X,Q)$ with a marked isomorphism $Q\cong E$, and we
recover $\psi$ after twisting by a line bundle of class $\phi([\sO_X])$.)
We thus find that the moduli stack of pairs $(X,\psi)$ (with blowdown
structure) consists of countably many copies of $E^{12-d}$, one for each
value of $\deg(\phi([\sO_X]))\in \Z$.  (More precisely, if one allows $E$
to vary, one obtains the corresponding power of the universal curve over
the moduli stack of elliptic curves, since one can act on $\phi$ by
$\Aut(E)\ltimes E$.)

Note that the moduli space of possible $\psi$s for a given choice of
$(X,Q)$ is 2-dimensional: there is a degree of freedom corresponding to a
choice of marked point (to identify $Q$ with an elliptic curve) and a
degree of freedom corresponding to the choice of $\psi(\sO_Q)$.  Thus, for
instance, the moduli stack of pairs $(X,Q)$ with $X$ a noncommutative
$\P^2$ may be identified with the moduli space of cubic plane curves with a
choice of translation, or equivalently by the moduli stack of smooth genus
1 curves with marked line bundles of degrees 3 and 0.  Making $Q$ an
elliptic curve requires marking a point, at which point we may use
$\Pic^2(E)\cong \Pic^0(E)\cong E$ to write the family as $E^2$, where $E$
is the universal elliptic curve.  This determines how $\psi$ acts on
structure sheaves of points, but it remains to specify its image on
$\sO_Q$, which is a point of $\Pic(E)\cong \Z\times E$, giving countably
many copies of $E^3$.

There is a natural action of $W(E_{9-d})$ on $\NS(X)$, which extends to an
action of $W(E_{9-d})\ltimes \NS(X)$ on $K_0(X)$.  This arises from the
fact that line bundles span $K_0(X)$ (and, in fact, there is a basis of
line bundles), and thus we can define an action by
\[
  [\sO_X(D)]\mapsto [\sO_X(wD+D_0)].
\]
There is then an induced action of this group on our moduli stack, and we
obtain a derived equivalence between the categories of coherent sheaves on
any two surfaces in the same orbit.  (Indeed, $W(E_{9-d})$ acts on
commutative del Pezzo surfaces preserving the structure sheaf, while
$\NS(X)$ corresponds to twists by line bundles.  One caveat is that in the
noncommutative setting, twisting by a line bundle changes the surface!)
Similarly, there is an action of $\Aut(E)\ltimes E\times \Pic(E)$ on the
other side, and this preserves the surface (since it only changes $\psi$).

Note that the autoequivalence $\Psi_X$ of $\coh(E)$ corresponding to the
action of $\theta^{-1}$ depends only on $\phi([\text{pt}])$,
$\phi([\sO_Q])$ and the integer $\delta=\deg(\phi([\sO_X]))$.  Indeed, the
automorphism $\tau$ is translation by $q:=\det\phi([\text{pt}])$, so that
$\Psi_X(\sO_x)\cong \sO_{x+q}$, while the line bundle is determined by the
equation
\[
[{\cal L}]+\delta q+\phi([\sO_X]) = \Psi_X(\phi([\sO_X])) =
\phi([\sO_X(Q)]) = \phi([\sO_X])+\phi([\sO_Q])+d q,
\]
giving
\[
[{\cal L}] = (d-\delta)q + \phi([\sO_Q]).
\]
(One caveat, which we will discuss further below, is that the {\em
  isomorphism class} of ${\cal L}$ is fixed, but ${\cal L}$ itself is a
more subtle issue.)  Thus any given isomorphism class of $\Psi$ cuts out a
codimension 2 subscheme of the moduli stack on which $\phi([\text{pt}])$
and $\phi([\sO_Q])$ are fixed.  For $d<8$, the two classes $[\text{pt}]$
and $[\sO_Q]$ span a saturated sublattice of $K_0(X)$, so that this
subscheme is itself a power of $E$, namely $E^{10-d}$; this also holds on
the component for $d=8$ corresponding to odd Hirzebruch surfaces.  (For
$\P^2$, the subscheme is a coset of $E[3]\times E$, while for $F_0/F_2$ it
is a coset of $E[2]\times E^2$.)

It will be useful to give names for the coordinates in $E^{12-d}$.  Recall
that the expression of $X$ as an iterated blowup of $\P^2$ (i.e., such an
expression of the underlying commutative surface) induces a basis of
$\NS(X)$ of the form
\[
h,e_0,\dots,e_{8-d}
\]
which in turn induces a basis
\[
  [\sO_X],[\sO_h(-1)],[\sO_{e_0}(-1)],\dots,[\sO_{e_{8-d}}(-1)],[\text{pt}]
\]
of $K_0(X)$, so that $(X,\psi)$ is determined by
\[
u:=\det\phi[\sO_X],\quad h:=\det\phi[\sO_h(-1)],\quad q:=[\text{pt}]
\]
and
\[
x_i:=\det\phi[\sO_{e_i}(-1)].
\]
(An expression of $X$ as an iterated blowup of an even Hirzebruch surface
gives a similar parametrization, but with $h$ and $e_0$ replaced by $s$ and
$f$.)

Now, given a representable numerical exceptional collection
$([E_1],\dots,[E_n])$ in $K_0(X)$ (satisfying the positive root condition
where the slopes agree), there is an associated family $(E_1,\dots,E_n)$ of
exceptional collections on the moduli stack, and the bundles
$(E_1|_Q,\dots,E_n|_Q)$ can be read off (up to isomorphism) from the
classes $\phi([E_i])$, since any stable sheaf on $E$ is determined by its
class in $K_0(E)$.  Thus fixing the bundles is again tantamount to fixing
the values under the corresponding homomorphisms from $E^{12-d}$ to $E$.
An exceptional collection always forms a saturated sublattice (the matrix
of Mukai pairings has determinant 1), and thus fixing the bundles gives a
family isomorphic to $E^{12-n-d}$.

We thus see in general that the family of surfaces with blowdown structure
and choice of $\psi$ such that the given collection gives rise to a given
$B_{V,\Psi}$ is the appropriate fiber of a homomorphism $E^{12-d}\to
E^{n+2}$.  (Of course, the fiber will be empty unless the values
$\rank(\phi([E_i]))$, $\deg(\phi([E_i]))$ are all correct, a constraint on
the numerical collection along with the integer $\delta$.)  We must then
take the union over all numerical collections and quotient by the action of
$W(E_{9-d})\ltimes \NS(X)$ to obtain the ``true'' moduli space.
(Equivalently, we choose one representative from each orbit of valid
numerical collections, and for each one, take the quotient of the relevant
subscheme of $E^{12-d}$ by the action of the stabilizer, a subgroup of the
affine Weyl group $W(E_{9-d})\ltimes Q^\perp$.)

One caveat here is that reflections do not give {\em abelian}
autoequivalences; they do not in general preserve the effective cone.  In
particular, it is not a priori obvious that $W(E_{9-d})$ preserves the
graded algebra.  Certainly, where none of the roots of $E_{9-d}$ are
effective, $W(E_{9-d})$ acts as abelian equivalences, and thus there is no
difficulty.  In particular, for any coordinate on the appropriate weighted
projective space, the two algebras we want to compare induce sections of
the same line bundle on the relevant subscheme of $E^{12-d}$ that agree
where no roots are effective (or, more precisely, where the given $w$ does
not make any effective root negative).  So the only way those coordinates
can fail to agree is if the condition on roots {\em always} fails.  But in
that case (which since it holds for all $q$, must come from a case in which
the underlying commutative surface has a $-2$-curve in which we wish to
reflect), the reflection corresponding to such a root actually acts
trivially on the entire component.  It follows in either case that the map
to the deformation {\em scheme} factors through the quotient by
$W(E_{9-d})\ltimes \NS(X)$.  (Caveat: the image in the deformation {\em
  stack} is not the quotient stack: weighted projective spaces have cyclic
stabilizers, unlike the action of $W(E_{9-d})\ltimes \NS(X)$ on $E^{12-d}$.
This is presumably related to the fact that the quotient by the reflection
subgroup of the stabilizer of a point is regular at that point; quotienting
by that normal subgroup gives a cyclic group in the cases of interest.)

It remains to understand the morphism from this family to deformation
space.  The missing piece here is that there is a natural line bundle on
deformation space (it is a quotient by $\G_m$, after all), and thus we need
to understand the pullback of this line bundle to the relevant family of
surfaces.  (The image will then be the corresponding embedding in weighted
projective space.)  It turns out that this is in fact the restriction of a
line bundle on the ambient moduli space of pairs $(X,\psi)$.  The main
issue is the fact that when defining the moduli space, we only determined
the various exceptional sheaves and ${\cal L}$ up to isomorphism.  For the
exceptional sheaves, this is not a problem: rescaling the exceptional
sheaves by scalars simply acts on $B_{V,\Psi}$ by an inner automorphism.
The line bundle is thus entirely determined by the failure of the induced
autoequivalence $\Psi_X$ of $Q$ to be a pullback from the corresponding
family over $E^2$.

We thus need to understand this autoequivalence better.  For convenience,
we consider its inverse instead.  A pair $(X,\psi)$ determines an adjoint
pair of functors $j_*:\coh E\to \coh X$, $j^*:\coh X\to \coh E$, and there
is a functorial distinguished triangle
\[
\Psi_X^{-1}[1]\to Lj^*j_*\to \text{id}\to
\]
(actually a direct sum decomposition of $Lj^*j_*$).  Since
$\Psi_X^{-1}M\cong (\Psi_X^{-1}\sO_E)\otimes \tau_q^*M$, it remains only to
compute $\Psi_X^{-1}\sO_E$.  This is a line bundle on $E\times E^{12-d}$,
and is functorial in $E$, so is really a line bundle on the corresponding
moduli stack.  Such line bundles were classified in
\cite[Prop. 2.2]{elldaha}, and have the form ${\cal L}_{Q,w}$ with
``weight'' $w\in \Z/12\Z$ (the isomorphism class of the pullback of the
identity section) and ``polarization'' $Q$ (the class in the N\'eron-Severi
group of the generic fiber).  The polarization is a symmetric element of
$\End(E\times E^{12-d})\cong \Mat_{13-d}(\Z)$, and thus may be expressed as
a quadratic form $vQv^t/2$ in the coordinates, in our case $z$ on $E$ and
$u,h,q,x_0,\dots,x_{8-d}$ on $E^{12-d}$.  (Note that homomorphisms between
powers of $E$ are generically given by integer matrices, and the pullback
acts on the polarization in the obvious way.  Also, the line bundle is
relatively ample iff the polarization is positive definite.)

Composing $\psi$ with an equivalence $\coh(X)\cong \coh(X')$ does not
change $Lj^*j_*$, so leaves $\Psi_X^{-1}$ alone.  Not only does this imply
that $\Psi_X^{-1}$ is well-defined as a functor on the category of sheaves
on the universal curve over the parameter space, but that it is invariant
under the action of $\NS(X)$.  Similarly, if we compose with an
autoequivalence of $\coh(E)$, this simply conjugates $\Psi_X^{-1}$ by that
autoequivalence, and thus we know how $\Psi_X^{-1}\sO_E$ must transform
under that action on parameter space.

The action of $E\times \Pic(E)$ suffices to determine the dependence of
$\Psi_X^{-1}\sO_E$ on $u$ and $\delta$ (i.e., it tells us that the
polarization has the form $uq-\delta(qz+q^2/2)$ plus terms independent of
$u$ and $\delta$).  The twist-by-$h$ symmetry implies that the polarization
is invariant under $(h,u,\delta)\mapsto (h+q,u+2q+h,\delta+3)$, which
together with the known dependence on $u$ and $\delta$ tells us that the
$h$-dependent terms are precisely $-h^2/2+3hz$.  Similarly, invariance
under $(x_i-q,u+x_i,\delta+1)$ (i.e., twisting by $e_i$) implies that the
$x_i$-dependent terms are precisely $x_i^2/2-x_iz$.  Since the
$z$-dependence is determined by the {\em isomorphism classes} of the fibers
over $E^{12-d}$ (i.e., that they correspond to the point
$-3h+x_0+\cdots+x_{8-d}+\delta q$ of $\Pic^{-d}(E)$), we thus conclude that
\[
\Psi_X^{-1}\sO_E
\cong
{\cal L}_{-dz^2/2+(3h-x_0-\cdots-x_{8-d}-\delta q)z
-h^2/2+x_0^2/2+\cdots+x_{8-d}^2/2+uq-(\delta+a)q^2/2,w},
\]
for some unknown $a\in \Z$ and $w\in \Z/12\Z$.  (One can also easily verify
that this is invariant under $W(E_{9-d})$.)

Since changing $a$ and $w$ merely twists by a line bundle pulled back from
the curve parametrizing $q$, this is enough to tell us the following.
(Note that the dependence on $h,\vec{x}$ is essentially just the
intersection form, making the analogue for $F_0/F_2$ straightforward to
write down.)

\begin{prop}
  If we fix $E$, $V$, $\Psi$, then for any family of filtered deformations
  of $B_{V,\Psi}$ constructed via Proposition \ref{prop:fildef_from_sheaf},
  the pullback of the natural line bundle on deformation space is
  isomorphic to the restriction of ${\cal
    L}_{-h^2/2+x_0^2/2+\cdots+x_{8-d}^2/2+uq,0}$.
\end{prop}

\begin{proof}
  Twisting by the inverse of this bundle turns $\Psi_X^{-1}\sO_E$ into
\[
{\cal L}_{-dz^2/2+(3h-x_0-\cdots-x_{8-d}-\delta q)z-(\delta+a)q^2/2,w},
\]
which is the pullback through the map $(z,3h-x_0-\cdots-x_{8-d},q)$ of a
family of line bundles depending only on the isomorphism class of the
autoequivalence $\Psi$.  In particular, for fixed $\Psi$, twisting by the
given line bundle makes the associated graded of the universal filtered
deformation globally isomorphic to $B_{V,\Psi}$ as required.
\end{proof}

%\begin{rem}
%  One can in fact show that
%  \[
%  \Psi^{-1} \sO_E
%  \cong
%  {\cal L}_{-dz^2/2 + (3h-x_0-\cdots-x_{9-d}-\delta q)z - h^2/2+x_0^2/2+\cdots
%  +x_{9-d}^2/2+uq-(\delta+1)q^2/2,d-5},
%  \]
%  but the full calculation requires additional machinery not otherwise
%  relevant to the present work.
%\end{rem}

Since the polarization is not positive definite, this line bundle is not
ample on all of $E^{12-d}$.  On the fibers of $(q,3h-x_0-\cdots-x_{8-d})$,
the polarization becomes positive {\em semi}definite, and becomes positive
definite once we fix the image under $\psi$ of some exceptional sheaf of
positive rank (as this fixes a linear combination of $u$ and the other
parameters).  We thus see that we indeed obtain an ample bundle on our
family of deformations in this way, and thus a map to a weighted projective
space, in such a way that the deformation parameters are homogeneous in the
coordinates on weighted projective space.

As an example, let us consider the case that $V$ is semistable of rank $r$
and slope $m\in \Z$ (divisorial, for simplicity) while $\Psi$ increases
degrees by $d$.  Since changing $\deg\phi([\sO_X])$ adds the same integer
to all slopes in the exceptional collection, the corresponding exceptional
collection must consist of line bundles of the same slope.  It must
therefore have the form
\[
(\sO_X(D-\alpha_1-\cdots-\alpha_{r-1}),
 \sO_X(D-\alpha_1-\cdots-\alpha_{r-2}),
 \dots
 \sO_X(D))
\]
where $D$ is some divisor class and $\alpha_1,\dots,\alpha_{r-1}$ are the
simple roots of a subsystem of type $A_{r-1}$ inside $E_{9-d}$.  The orbits
of such collections under $W(E_{9-d})\ltimes \NS(X)$ are in natural
correspondence with the conjugacy classes of $A_{r-1}$ subsystems, and for
each orbit, we may take the representative with $D=0$.  (This forces us to
take $\delta=m$ in the map $K_0(X)\to K_0(E)$, though of course in this
case we can tensor $V$ by a line bundle to obtain an equivalent problem
with $m=0$.)  

For $r=1$, there is (trivially) a unique orbit of $A_0$ subsystems, and we
find that, for $d\le 6$, the base of the corresponding family is the
quotient of (a torsor over) $E^{9-d}$ by $W(E_{9-d})$, with line bundle
corresponding to the minimal invariant polarization.  (Minimality can be
read off from the relation of the above polarization to the intersection
form, so that any root evaluates to 1.)  Per Looijenga
\cite{LooijengaE:1976}, such a quotient for an irreducible Weyl group is
isomorphic to a weighted projective space, in which the degrees of the
generators are the labels on the affine Dynkin diagram.  We thus obtain the
following descriptions of the quotient subschemes (corresponding to $E_8$,
$E_7$, $E_6$, $D_5$, $A_4$, $A_2A_1$:
\[
\P^{[122334456]},\P^{[11222334]},\P^{[1112223]},\P^{[111122]},\P^4,
\P^2\times \P^1.
\]
For $d=7$, the Weyl group is too small, so we obtain a quotient of $E^2$ by
$\Sym_2$, while for $d=8$, there are two components (corresponding to the
two parities of Hirzebruch surface), one isomorphic to $E$ (the odd case),
and one isomorphic to $E[2]\times E/A_1$ (the even case).  Finally, for
$d=9$, we obtain $E[3]$.

For $r=2$, $d\le 5$, there is again a unique orbit of $A_1$ subsystems,
and thus we must replace $E_{9-d}$ by the stabilizer of a root in
$E_{9-d}$.  This gives ($E_7$, $D_6$, $A_5$, $A_3A_1$)
\[
\P^{[11222334]},\P^{[1111222]},\P^5,\P^3\times \P^1
\]
for $d\le 4$.  For $d=5$, the rank is too small, so we have an elliptic
factor, and for $d=6$, there are two orbits of roots, one of which has an
elliptic factor.  For $d=7$, we have a unique orbit of roots, with an
elliptic factor, while for $d=8$ only the even case has a root, making the
moduli space in that case a torsor over $E[2]$, and for $d=9$ no
configuration exists.

For $r=3$, we similarly obtain
\[
\P^{[1112223]},\P^5,\P^2\times \P^2
\]
for $d\le 3$; for $d=4$, $d=5$, there are unique orbits but with a rank
deficit.  For $d=6$, there is a unique orbit, and the space has the form
$E[3]\times \P^1$.  For $d>6$, no configuration exists.

For $r=4$, we have
\[
\P^{[111122]},\P^3\times \P^1,
\]
for $d=1,2$, a unique orbit with a rank deficit for $d=3$, two orbits for
$d=4$ (one with a rank deficit, and one with a torsion factor), and two
orbits for $d=5$ (the maximum), again with a rank deficit.

For $r=5$, we have $\P^4$ for $d=1$, a unique orbit with a rank deficit for
$d=2,3$, two orbits with a rank deficit for $d=4$, and for $d=5$ obtain
$E[5]$.

For $r=6$, we have $\P^2\times \P^1$ for $d=1$, two orbits (one with a rank
deficit, the other with torsion) for $d=2$, two orbits (with torsion) for
$d=3$, and nothing for $d>3$.  For $r=7$, we have a rank deficit for
$d=1,2$ and nothing for $d>2$, while for $r=8$ we have two components (one
with a rank deficit and one with torsion) for $d=1$ and $E[2]$ for $d=2$,
and for $r=9$ we have $E[3]$ for $d=1$.  (In general, if $r+d>10$, then no
configuration exists, since the rank of the subsystem is too large.)

The astute reader will have noticed a symmetry above: the description of
the given portion of deformation space is invariant under swapping $r$ and
$d$.  Indeed, if we blow up $d$ further points, then we obtain an
exceptional collection on a (deformed!) elliptic surface, and a derived
equivalence swaps the roles of the $r$ line bundles and the $d$ exceptional
curves.  Such a symmetry holds in general for the semistable case, with the
one caveat that it will in general change the slope from $a/r$ to $b/r$
with $ab=1(r)$.

In many of the above cases, the image of the family in weighted projective
space contains an orbifold point, which leads to additional interesting
families.  In general, if the space for $B_{V,\Psi}$ passes through an
orbifold point of degree $g$, then the corresponding deformation extends to
one of $B_{V',\Psi'}$ with
\[
V' = V\oplus \Psi(V)\oplus \cdots \oplus \Psi^{g-1}(V),
\quad
\Psi' = \Psi^g.
\]
We can then further deform $V'$, $\Psi'$.

For instance, since the family for $d=r=1$ is $\P^{[122334456]}$, we have
orbifold points of order $2\le g\le 6$.  The result corresponds to taking
$V$ to be a sum of line bundles of slopes $0,1,\dots,g-1$, with $\Psi$
increasing degrees by $g$.  The differences between the Chern classes of
the corresponding line bundles on $X$ (now a del Pezzo surface of degree
$g$) must then be rational curves.  It is straightforward to classify such
configurations, and for $g=2,3,4$, we find that the family is the
Looijenga-style quotient associated to $E_6$, $D_4$, $A_2$ respectively.
For $g=5$, there is a unique filtered deformation, while for $g=6$, the
{\em generic} such $V$ has no filtered deformations, but there is a
codimension 2 subfamily of $V$s with unique filtered deformations.  Indeed,
in this case, we find that the exceptional collection is actually a {\em
  full} exceptional collection, and thus $\Psi$ is determined by $V$.

Of course, the cases corresponding to orbifold points are special cases of
this scenario, and have an additional symmetry corresponding to shifting
the slopes by $1$ (and applying $\Psi$ or $\Psi^{-1}$ as needed).  This
relates to the fact that for $d\le 6$, the group $W(E_{9-d})\ltimes \NS(X)$
has the affine Weyl group (or product of two affine Weyl groups, for $d=6$)
as a normal subgroup, and thus we have an isomorphism
\[
W(E_{9-d})\ltimes \NS(X)\cong \Z\ltimes \tilde{W}(E_{9-d}),
\]
where the splitting is induced by the requirement that $\Z$ preserves the
set of simple roots.  In other words, the group $W(E_{9-d})\ltimes \NS(X)$
contains an element $\pi$ corresponding to a diagram automorphism.  Since
$\tilde{W}(E_{9-d})$ does not change slopes, $\pi$ must change the slope
(increasing it by 1, say).  Moreover, $\pi^d$ commutes with the affine Weyl
group, and thus (since it increases slopes by $d$) is translation by $Q$.
We can then verify that the divisor classes $\pi(0),\dots,\pi^{d-1}(0)$ all
correspond to rational curves, and thus induce a configuration as in the
previous paragraph.  The orbifold case then corresponds to the case that
the point in the moduli space is invariant under $\pi$.

For $d=2$, $r=1$, we have decimations for $g\in \{2,3,4\}$.  In each case,
there is a single configuration (up to the Weyl group action and choice of
parity for $2g=8$).  For $g=2$, $g=3$, the result is a Looijenga quotient,
with Weyl groups $D_4$ and $A_1$ respectively, while for $g=4$, the result
is again a full exceptional collection (in an even Hirzebruch surface), and
thus there is a unique deformation on the appropriate codimension 2
subscheme.

For $d=1$, $r=2$, we again have $D_4$ and $A_1$ for $g\in \{2,3\}$ and a
unique deformation over a codimension 2 subscheme of the base for $g=4$,
with the only difference from the $(d,r)=(2,1)$ case being that the
surfaces have degree $g$ rather than $2g$.  (In particular, for $g=4$, we
obtain a full exceptional collection of line bundles on a quartic del Pezzo
surface.)

For $d=3$, $r=1$, the $g=2$ case is the Looijenga quotient of type $A_2$,
while the $g=3$ case is a full exceptional collection on $\P^2$.  The $g=3$
case corresponds to the fact that the standard noncommutative deformation
(the three generator Sklyanin algebra
\cite{ArtinM/TateJ/VandenBerghM:1990}) of $k[x,y,z]$ is obtained by
removing the degree $3$ equation from the $d=3$, $r=1$ elliptic algebra.

For $d=1$, $r=3$, the $g=2$ case is again the Looijenga quotient of type
$A_2$, while the $g=3$ case is a full exceptional collection on a
noncommutative cubic surface.

For $d=4$, $r=1$, $g=2$, the exceptional collection lives on an even
Hirzebruch surface, and $V$ imposes a single (saturated!) condition on
$\Psi$.  So deformations only exist on a codimension $1$ subscheme of
parameter space, but when they exist are classified by $\P^1\cong E/A_1$.
Each such deformation corresponds to removing a quadratic relation from the
algebra with $d=4$, $r=1$; since the scheme is $\P^1$ embedded by $O(1)$,
the possibilities span a $2$-dimensional subspace of the dual.  Removing
both relations gives Sklyanin's deformation \cite{SklyaninEK:1982} of
$k[x,y,z,w]$ (with the two removed relations corresponding to the
$2$-dimensional space of central elements of degree 2.)

For $d=1$, $r=4$, $g=2$, we again have a $\P^1$ of deformations on a
codimension 1 subscheme.  Again, we could also remove both relations, which
gives something best described as an Azumaya algebra on a noncommutative
Fano 3-fold; it is unclear how to characterize that 3-fold, however.

Finally, for $d=r=2$, $g=2$, deformations only exist on a codimension 1
subscheme, while the general nonempty fiber is the Looijenga quotient of
type $A_3$, a.k.a. $\P^3$.

\medskip

We now turn to semisimple examples with non-integer slope.  Note that in
the case of a single bundle, the process of testing representability of the
candidate exceptional classes involves first putting them in the
fundamental alcove.  It is of course quite straightforward to find all
candidates in the fundamental alcove, and for each candidate surviving the
Euler characteristic test, we can read off the stabilizer from its standard
coordinates.  This works particularly well for $r=d=1$; for larger values
of $r$, we need only classify configurations of $r-1$ roots that have
intersection $1$ (modulo the denominator) with the fundamental
representative and $-1$ with each other, while for larger values of $d$, we
need only find $d-1$ exceptional curves that are orthogonal (modulo the
denominator) and to each other.  (And when both are larger, the roots
must be orthogonal to the $-1$-curves.)

For $\mu=-1/2$, we obtain the following nice cases (omitting cases where there
is either a rank deficit or torsion) as $r$, $d$ vary:
\[
\begin{pmatrix}
  D_8 & A_7 & A_5A_1\\
  A_7 & A_3A_3 \\
  A_5A_1
\end{pmatrix}
\]
In the case $r=9$, $d=1$ (or vice versa), there is a unique combinatorial
configuration, with the corresponding subscheme of parameter space a torsor
over $E[10]$.  This case is primarily of note because there actually are
{\em no} bundles of slope $-1/2$ on a noncommutative $\P^2$, and thus this
case can only arise by taking $\delta\ne 0$.

In the $r=d=1$ case, there are orbifold points, so we can again take a
decimation.  It turns out that for bundles of slope $-3/2$, $-1/2$ (on a
degree 2 del Pezzo surface) to form an exceptional collection, it is
necessary for their Chern classes to differ by $Q$.  This means that a
deformation will only exist on a codimension 1 subscheme of parameter
space, and when it does will be classified by the quotient corresponding to
$A_7$ (i.e., the same as for $r=1$, $d=2$).

For slope $-1/3$ (or equivalently $-2/3$), we obtain the nice cases
\[
\begin{pmatrix}
  A_8 & A_5A_2\\
  A_5A_2
\end{pmatrix}.
\]
For slope $-1/4$, the $r=d=1$ case is nice, with corresponding Weyl group
$A_7A_1$.  For slope $-1/5$, it appears that no cases are nice, but for
$r=d=1$, slope $-2/5$, the $r=d=1$ case is nice, with group $A_4A_4$.

One caveat here is that we have not shown that every invariant section of a
power of the ample bundle can be expressed as a polynomial in the
coefficients of the associated filtered deformation.  For sections of the
ample bundle itself, this would reduce to a calculation of the
infinitesimal deformation of $\Proj(B[t])$ associated to the $\G_m$ action
on the general filtered deformation in the family; this, together with
analogous calculations for orbifold points, should be enough to establish
the result in general.  (Note that this would also establish that the map
from the quotient scheme to the deformation scheme is an embedding.)

\medskip

The fact that precisely these weighted projective spaces arise when
computing moduli spaces of filtered deformations of random elliptic
algebras suggests the following.

\begin{conj}\label{conj:main}
  Any nontrivial filtered deformation of any elliptic algebra $B_{V,\Psi}$
  has Rees algebra of the form $A^+_M$ for some rigid, unobstructed,
  torsion-free sheaf $M$ on a noncommutative del Pezzo surface.
\end{conj}

\begin{rem}
  Since the moduli stack is closed in a weighted projective space, the
  Hilbert series of the moduli stack is upper semicontinuous.  As a result,
  the argument of Lemma \ref{lem:no_neg_derivs_general} would let us reduce
  to the case that $\tau$ has finite order (so the construction involves
  Azumaya algebras on del Pezzo surfaces), at the cost of having to work in
  finite characteristic.  We also immediately deduce that if the claim
  holds for a given algebra $B$, then it holds for the generic algebra of
  the same shape.
\end{rem}

In the cases where this would make the deformation space a weighted
projective space, there is a possible strategy for proving this: show that
for each degree, the space of (nontrivial) infinitesimal deformations
agrees with the putative number of generators.  This implies that the
deformation space is {\em contained} in the conjectured weighted projective
space, and therefore equal.

Although one could hope to compute the desired Hochschild cohomology spaces
by hand, this does at the very least simplify the calculation when done by
computer.  In particular, for the line bundle cases with $r=1$ and $d\in
\{1,2\}$, we could carry out the calculation not just for generic
parameters, but for {\em general} parameters: that is, for any $E$ and any
$\Psi$ that increases degrees by $1$ or $2$, every filtered deformation of
the algebra $B_{\sO_E,\Psi}$ arises from the surface construction.  We have
also checked the conjecture for {\em random} instances (over a moderately
large finite field) in every slope 0 case (including those where the family
is not nice), as well as for a few of the cases with noninteger slope; as
mentioned, this implies the conjecture holds generically in those cases as
well.

Although the general twisted case appears out of reach at the moment, there
is some hope for the untwisted case, since any deformation has a large
center.  The first step would be to show that the filtered deformation
arises from a sheaf on a {\em singular} del Pezzo surface, which would
reduce to the following (which would follow from the main conjecture, as
would the analogous statement when twisted by a torsion point, since the
center of a noncommutative del Pezzo surface with $q$ torsion is a
commutative del Pezzo surface).

\begin{conj}
  Any filtered deformation of an untwisted elliptic algebra $B_{V,\Psi}$
  which induces a trivial deformation of the center is trivial.
\end{conj}

\begin{rem}
  This holds for $d=1$ if $V$ is stable of rank 2 or a sum of two degree 0
  line bundles, as in both cases we showed above (Examples \ref{eg:d1stab}
  and \ref{eg:d1r2} respectively) that the moduli space is the correct
  weighted projective space.
\end{rem}

The other missing ingredients are showing that the resulting algebra on a
singular del Pezzo surface lifts to a flat sheaf of algebras on the minimal
desingularization, and showing that the result is still the endomorphism
algebra of a rigid, unobstructed object.  (In the twisted case, one would
also need to show that one obtains the correct class in the Brauer group of
the affine del Pezzo surface; this should not be an issue in the untwisted
case where the algebra remains Azumaya on the compactification, since the
compactification has trivial Brauer group.)

\bibliographystyle{plain}

\end{document}